\PassOptionsToPackage{prologue,dvipsnames}{xcolor} 
\documentclass[hidelinks,onefignum,onetabnum]{siamart250211}

\usepackage{lipsum}
\usepackage{amsfonts}
\usepackage{graphicx}
\usepackage{epstopdf}
\usepackage{algorithmic}
\ifpdf
  \DeclareGraphicsExtensions{.eps,.pdf,.png,.jpg}
\else
  \DeclareGraphicsExtensions{.eps}
\fi

\usepackage{arydshln} 
\usepackage{tikz}
\usepackage{tikz-cd}
\usepackage{stmaryrd} 
\usepackage{tablefootnote} 
\usepackage{graphbox}
\usepackage{dashrule}
\usepackage{mathtools} 
\usepackage[shortlabels]{enumitem}
\usepackage{subcaption}

\newsiamremark{remark}{Remark}
\newsiamremark{hypothesis}{Hypothesis}
\crefname{hypothesis}{Hypothesis}{Hypotheses}
\newsiamthm{claim}{Claim}
\newsiamremark{fact}{Fact}
\crefname{fact}{Fact}{Facts}

\newsiamremark{example}{Example}
\crefname{example}{Example}{Examples}
\newsiamremark{assumption}{Assumption}
\crefname{assumption}{Assumption}{Assumptions}

\headers{Multirate NPRK Methods}{T. Buvoli, B.K. Tran, and B. S. Southworth}

\title{Multirate Runge-Kutta for Nonlinearly Partitioned Systems\thanks{Submitted to the editors DATE.
\funding{T.B. was funded by the National Science Foundation under grant NSF-OIA-2327484. B.S.S. was supported by the DOE Office of Advanced Scientific Computing Research Applied Mathematics program through Contract No. 89233218CNA000001. B.K.T. was supported by the Marc Kac Postdoctoral Fellowship at the Center for Nonlinear Studies at Los Alamos National Laboratory. Los Alamos National Laboratory report number LA-UR-25-23138.}}}

\author{
    Tommaso Buvoli\thanks{Tulane University, Department of Mathematics, New Orleans, Louisiana 70118}
    \and
    Brian K. Tran\thanks{Los Alamos National Laboratory, Theoretical Division, Los Alamos, New Mexico 87545 \\ \quad E-mail address: \texttt{tbuvoli@tulane.edu, btran@lanl.gov, southworth@lanl.gov } }
    \and
    Ben S. Southworth\footnotemark[3]
}

\usepackage{amsopn}

\newcommand{\cmp}[2]{#1^{\{#2\}}}
\newcommand{\Mcmp}[1]{\cmp{M}{#1}}
\newcommand{\Acmp}[1]{\cmp{A}{#1}}
\newcommand{\acmp}[1]{\cmp{a}{#1}}
\newcommand{\bcmp}[1]{\cmp{b}{#1}}
\newcommand{\scmp}[1]{\cmp{s}{#1}}

\newcommand{\Mhcmp}[1]{\cmp{\widehat{M}}{#1}}
\newcommand{\Ahcmp}[1]{\cmp{\widehat{A}}{#1}}
\newcommand{\AhcmpB}[1]{\cmp{\widehat{\mathbf{A}}}{#1}}
\newcommand{\bhcmpB}[1]{\cmp{\widehat{\mathbf{b}}}{#1}}
\newcommand{\ahcmp}[1]{\cmp{\widehat{a}}{#1}}
\newcommand{\bhcmp}[1]{\cmp{\widehat{b}}{#1}}
\newcommand{\chcmp}[1]{\cmp{\widehat{c}}{#1}}

\newcommand{\AcmpB}[1]{\cmp{\mathbf{A}}{#1}}
\newcommand{\bcmpB}[1]{\cmp{\mathbf{b}}{#1}}
\newcommand{\Scmp}[1]{\cmp{\mathcal{S}}{#1}}

\def\treetype{2} 

\if\treetype1 

	\def\T#1#2{\scalebox{1}{\enspace\raisebox{-1mm}{\csname Tree#1#2\endcsname}\enspace}}
	\tikzstyle{nodeB}=[circle,fill=black,draw=black,scale=0.20,line width=0.8pt,radius=0.7pt]
	\tikzstyle{nodeW}=[circle,fill=black,draw=black,scale=0.24,line width=0.8pt,radius=0.7pt,fill=white]
	
	\tikzstyle{treeX}=[rotate=180,scale=0.15,line width=1.1pt]
	\expandafter\newcommand\csname Tree10\endcsname{\begin{tikzpicture}[treeX]\node[nodeB]{};\end{tikzpicture}}
	
	\expandafter\newcommand\csname Tree11\endcsname{\begin{tikzpicture}[treeX]\node[nodeW]{};\end{tikzpicture}}
	
	\expandafter\newcommand\csname Tree20\endcsname{\begin{tikzpicture}[treeX]\node[nodeB]{}child{node[nodeB]{}};\end{tikzpicture}}
	\expandafter\newcommand\csname Tree21\endcsname{\begin{tikzpicture}[treeX]\node[nodeB]{}child{node[nodeW]{}};\end{tikzpicture}}
	
	\expandafter\newcommand\csname Tree22\endcsname{\begin{tikzpicture}[treeX]\node[nodeW]{}child{node[nodeB]{}};\end{tikzpicture}}
	\expandafter\newcommand\csname Tree23\endcsname{\begin{tikzpicture}[treeX]\node[nodeW]{}child{node[nodeW]{}};\end{tikzpicture}}

	\expandafter\newcommand\csname Tree30\endcsname{\begin{tikzpicture}[treeX]\node[nodeB]{}child{node[nodeB]{}}child{node[nodeB]{}};\end{tikzpicture}}
	\expandafter\newcommand\csname Tree31\endcsname{\begin{tikzpicture}[treeX]\node[nodeB]{}child{node[nodeW]{}}child{node[nodeB]{}};\end{tikzpicture}}
	\expandafter\newcommand\csname Tree32\endcsname{\begin{tikzpicture}[treeX]\node[nodeB]{}child{node[nodeB]{}}child{node[nodeW]{}};\end{tikzpicture}}
	\expandafter\newcommand\csname Tree33\endcsname{\begin{tikzpicture}[treeX]\node[nodeB]{}child{node[nodeW]{}}child{node[nodeW]{}};\end{tikzpicture}}
	
	\expandafter\newcommand\csname Tree34\endcsname{\begin{tikzpicture}[treeX]\node[nodeW]{}child{node[nodeB]{}}child{node[nodeB]{}};\end{tikzpicture}}
	\expandafter\newcommand\csname Tree35\endcsname{\begin{tikzpicture}[treeX]\node[nodeW]{}child{node[nodeW]{}}child{node[nodeB]{}};\end{tikzpicture}}
	\expandafter\newcommand\csname Tree36\endcsname{\begin{tikzpicture}[treeX]\node[nodeW]{}child{node[nodeB]{}}child{node[nodeW]{}};\end{tikzpicture}}
	\expandafter\newcommand\csname Tree37\endcsname{\begin{tikzpicture}[treeX]\node[nodeW]{}child{node[nodeW]{}}child{node[nodeW]{}};\end{tikzpicture}}

	\expandafter\newcommand\csname Tree40\endcsname{\begin{tikzpicture}[treeX]\node[nodeB]{}child{node[nodeB]{}child{node[nodeB]{}}};\end{tikzpicture}}
	\expandafter\newcommand\csname Tree41\endcsname{\begin{tikzpicture}[treeX]\node[nodeB]{}child{node[nodeB]{}child{node[nodeW]{}}};\end{tikzpicture}}
	\expandafter\newcommand\csname Tree42\endcsname{\begin{tikzpicture}[treeX]\node[nodeB]{}child{node[nodeW]{}child{node[nodeB]{}}};\end{tikzpicture}}
	\expandafter\newcommand\csname Tree43\endcsname{\begin{tikzpicture}[treeX]\node[nodeB]{}child{node[nodeW]{}child{node[nodeW]{}}};\end{tikzpicture}}
	
	\expandafter\newcommand\csname Tree44\endcsname{\begin{tikzpicture}[treeX]\node[nodeW]{}child{node[nodeB]{}child{node[nodeB]{}}};\end{tikzpicture}}
	\expandafter\newcommand\csname Tree45\endcsname{\begin{tikzpicture}[treeX]\node[nodeW]{}child{node[nodeB]{}child{node[nodeW]{}}};\end{tikzpicture}}
	\expandafter\newcommand\csname Tree46\endcsname{\begin{tikzpicture}[treeX]\node[nodeW]{}child{node[nodeW]{}child{node[nodeB]{}}};\end{tikzpicture}}
	\expandafter\newcommand\csname Tree47\endcsname{\begin{tikzpicture}[treeX]\node[nodeW]{}child{node[nodeW]{}child{node[nodeW]{}}};\end{tikzpicture}}

\fi

\if\treetype2 

	\def\T#1#2{\scalebox{1.5}{\enspace\raisebox{-1mm}{\csname Tree#1#2\endcsname}\enspace}}
	\tikzstyle{dns}=[circle,fill=black,draw=black,scale=0.20,line width=0.5pt,solid] 
	\tikzstyle{edgeB}=[black,line width=1pt,solid]
	\tikzstyle{edgeW}=[black,line width=0.5pt,densely dotted]
	
	\tikzstyle{treeX}=[rotate=180,scale=0.15,line width=1.1pt]
	\expandafter\newcommand\csname Tree10\endcsname{\begin{tikzpicture}[treeX]\node[dns]{};\end{tikzpicture}}
	
	\expandafter\newcommand\csname Tree11\endcsname{\begin{tikzpicture}[treeX]\node[dns]{};\end{tikzpicture}}
	
	\expandafter\newcommand\csname Tree20\endcsname{\begin{tikzpicture}[treeX]\node[dns]{}child[edgeB]{node[dns]{}};\end{tikzpicture}}
	\expandafter\newcommand\csname Tree21\endcsname{\begin{tikzpicture}[treeX]\node[dns]{}child[edgeW]{node[dns]{}};\end{tikzpicture}}
	
	\expandafter\newcommand\csname Tree22\endcsname{\begin{tikzpicture}[treeX]\node[dns]{}child{node[dns]{}};\end{tikzpicture}}
	\expandafter\newcommand\csname Tree23\endcsname{\begin{tikzpicture}[treeX]\node[dns]{}child{node[dns]{}};\end{tikzpicture}}

	\expandafter\newcommand\csname Tree30\endcsname{\begin{tikzpicture}[treeX]\node[dns]{}child{node[dns]{}}child{node[dns]{}};\end{tikzpicture}}
	\expandafter\newcommand\csname Tree31\endcsname{\begin{tikzpicture}[treeX]\node[dns]{}child[edgeW]{node[dns]{}}child{node[dns]{}};\end{tikzpicture}}
	\expandafter\newcommand\csname Tree32\endcsname{\begin{tikzpicture}[treeX]\node[dns]{}child{node[dns]{}}child[edgeW]{node[dns]{}};\end{tikzpicture}}
	\expandafter\newcommand\csname Tree33\endcsname{\begin{tikzpicture}[treeX]\node[dns]{}child[edgeW]{node[dns]{}}child[edgeW]{node[dns]{}};\end{tikzpicture}}
	
	\expandafter\newcommand\csname Tree34\endcsname{\begin{tikzpicture}[treeX]\node[dns]{}child{node[dns]{}}child{node[dns]{}};\end{tikzpicture}}
	\expandafter\newcommand\csname Tree35\endcsname{\begin{tikzpicture}[treeX]\node[dns]{}child{node[dns]{}}child{node[dns]{}};\end{tikzpicture}}
	\expandafter\newcommand\csname Tree36\endcsname{\begin{tikzpicture}[treeX]\node[dns]{}child{node[dns]{}}child{node[dns]{}};\end{tikzpicture}}
	\expandafter\newcommand\csname Tree37\endcsname{\begin{tikzpicture}[treeX]\node[dns]{}child{node[dns]{}}child{node[dns]{}};\end{tikzpicture}}

	\expandafter\newcommand\csname Tree40\endcsname{\begin{tikzpicture}[treeX]\node[dns]{}child{node[dns]{}child{node[dns]{}}};\end{tikzpicture}}
	\expandafter\newcommand\csname Tree41\endcsname{\begin{tikzpicture}[treeX]\node[dns]{}child{node[dns]{}child[edgeW]{node[dns]{}}};\end{tikzpicture}}
	\expandafter\newcommand\csname Tree42\endcsname{\begin{tikzpicture}[treeX]\node[dns]{}child[edgeW]{node[dns]{}child[edgeB]{node[dns]{}}};\end{tikzpicture}}
	\expandafter\newcommand\csname Tree43\endcsname{\begin{tikzpicture}[treeX]\node[dns]{}child[edgeW]{node[dns]{}child[edgeW]{node[dns]{}}};\end{tikzpicture}}
	
	\expandafter\newcommand\csname Tree44\endcsname{\begin{tikzpicture}[treeX]\node[dns]{}child{node[dns]{}child{node[dns]{}}};\end{tikzpicture}}
	\expandafter\newcommand\csname Tree45\endcsname{\begin{tikzpicture}[treeX]\node[dns]{}child{node[dns]{}child{node[dns]{}}};\end{tikzpicture}}
	\expandafter\newcommand\csname Tree46\endcsname{\begin{tikzpicture}[treeX]\node[dns]{}child{node[dns]{}child{node[dns]{}}};\end{tikzpicture}}
	\expandafter\newcommand\csname Tree47\endcsname{\begin{tikzpicture}[treeX]\node[dns]{}child{node[dns]{}child{node[dns]{}}};\end{tikzpicture}}

\fi

\definecolor{mathematica_plot_blue}{RGB}{93.9463, 129.229, 180.998}
\definecolor{mathematica_plot_orange}{RGB}{224.584, 155.815, 36.223}
\definecolor{mathematica_plot_green}{RGB}{142.846, 176.35, 49.6957}

\makeatletter
\newcommand{\subalign}[1]{%
  \vcenter{%
    \Let@ \restore@math@cr \default@tag
    \baselineskip\fontdimen10 \scriptfont\tw@
    \advance\baselineskip\fontdimen12 \scriptfont\tw@
    \lineskip\thr@@\fontdimen8 \scriptfont\thr@@
    \lineskiplimit\lineskip
    \ialign{\hfil$\m@th\scriptstyle##$&$\m@th\scriptstyle{}##$\hfil\crcr
      #1\crcr
    }%
  }%
}
\makeatother

\ifpdf
\hypersetup{
  pdftitle={Multirate Runge-Kutta for Nonlinearly Partitioned Systems},
  pdfauthor={T. Buvoli, B. Tran, and B. S. Southworth}
}
\fi

\begin{document}

\maketitle

\begin{abstract}
Multirate integration is an increasingly relevant tool that enables scientists to simulate multiphysics systems. Existing multirate methods are designed for equations whose fast and slow variables can be linearly separated using additive or component-wise partitions. However, in realistic applications, this assumption is not always valid. Building on the recently developed class of nonlinearly partitioned Runge--Kutta (NPRK) methods, we develop a framework for multirate NPRK (MR-NPRK) that allows for arbitrary nonlinear splittings of the evolution operator.  We discuss order conditions, formalize different types of coupling between timescales, and analyze joint linear stability of MR-NPRK methods. We then introduce a class of 2nd- and 3rd-order methods, referred to as ``implicitly-wrapped'' multirate methods, that combine a user-specified explicit method for integrating the fast timescale with several slow implicit stages. These methods are designed to be algorithmically simple with low memory costs and minimal operator evaluations. Lastly, we conduct numerical experiments to validate our proposed methods and show the benefits of multirating a nonlinear partition.
\end{abstract}




Multirate time integration is a methodology for numerically solving coupled processes that evolve on different timescales. Multirate integrators differentiate themselves from single-rate integrators by employing different timestep sizes to capture fast and slow  dynamics. The general concept dates back over 60 years to the early works  \cite{Rice.1960,Andrus.1979,Gear.1984} that develop multirate  Runge--Kutta (RK) methods and multirate linear multistep (LMM) methods. 

In recent years, there has been a resurgence in literature on multiscale time integration techniques. These works can be broadly split into ``multirate'' (MR) \cite{Constantinescu.2022,Schafers.2023,Constantinescu.2007,Schlegel.2010,Schlegel.2012,Messahel.2021,Gunther.2016,Sarshar.2019,Roberts.2021}
 and ``multirate infinitesimal'' step (denoted as MRI or MIS) \cite{Wensch.2009,Wensch.2009,Sexton.2018,Sandu.2019,Luan.2020,Bauer.2021,Chinomona.2021nt,Fish.2024,Knoth.1998} methods. In contrast to the commonly used operator-split sub-cycling strategy, the broad idea behind these methods is that a modified fast equation is integrated by including ``slow tendency terms'' that reflect  dynamic coupling. Unlike MRI, MR methods have a predefined structure for integrating the fast dynamical timescale (e.g. a specific RK method run with a local time step that is smaller than the macro timestep). In contrast, MRI methods allow the user to integrate the fast component using any sufficiently accurate method.

The aforementioned multirate literature targets linearly split evolution operators such as the additive equation ${\dot{y} = \widehat{F}_f(y) + \widehat{F}_s(y)}$ and the partitioned equation $\dot{u} = \widehat{G}_f(u,v) , \dot{v} = \widehat{G}_s(u,v)$, where $\widehat{F}_f$ and $\widehat{G}_f$ govern the fast (``micro'') dynamic timescale while $\widehat{F}_s$ and $\widehat{G}_s$ govern a slow (``macro'') dynamic timescale. 
Both additive and component partitioned methods assume a linear separation of scales and stiffnesses. 
In realistic applications, this assumption is not always valid. This is particularly motivated by our recent and ongoing work in simulation of thermal radiation transport (TRT) \cite{imex-trt} and radiation hydrodynamics \cite{rad-hydro}, complex models fundamental to multiscale physics such as inertial confinement fusion. In these models, material coefficients and/or closure terms that nonlinearly interact with stiff processes depend directly on nonstiff variables, prohibiting the natural use of additive partitioning. Similarly, respecting energy balance relations prohibits the direct use of component partitioning. Additionally, it can be advantageous from the perspective of both computational efficiency and ease of implementation to have integrators that do not require treating nonlinearities  monolithically or assembling exact Jacobians. These types of considerations were a primary motivation for the development of semi-implicit methods \cite{Boscarino.2016,Boscarino.2015} and the more general nonlinearly partitioned Runge--Kutta (NPRK) methods \cite{nprk1,nprk2}, which naturally facilitate the high-order integration of an arbitrary nonlinear partitioning of equations.

 The primary objective of this paper is to develop a framework that enables multirate integration for both linear and nonlinearly partitioned systems. This new framework, named  multirate NPRK (MR-NPRK) is a generalization of our previous work on single-rate NPRK \cite{nprk1} methods. Since NPRK methods satisfy additive and partitioned order conditions \cite{nprk2}, the resulting methods can be naturally applied to any additive or component partitioned system. 
 
Another area of focus for this paper is the construction and analysis of implicit-explicit (IMEX) MR-NPRK. This class of integrators is motivated by multiphysics equations in which it is common that at least one variable or process needs to be treated implicitly. The development of IMEX multirate methods introduces important questions of joint stability; namely, if one can achieve damping of the fast process via implicit solves of the slow process. To this end, we introduce a method ansatz that: (1) guarantees joint L-stability in the limit of an infinitely stiff implicit variable, and (2) allows the user to select any explicit RK method for integrating the explicit component of the nonlinear partition. 
 
This paper is organized as follows. In \Cref{sec:nprk,sec:mr-nprk} we review NPRK methods, introduce MR-NPRK, and discuss order conditions and linear stability. The second half of this work is focused on IMEX MR-NPRK. In \Cref{sec:constructing} we develop a methodology for constructing simple and practical MR-NPRK methods for use in large-scale codes. Key objectives include (i) algorithmic simplicity, (ii) low storage requirements, (iii) minimal operator evaluations, and (iv) robust 2nd- and 3rd-order methods for linear or nonlinear partitions. Lastly, in \Cref{sec:numerical-experiments} we validate our MR-NPRK methods and highlight the advantages of nonlinear partitioning over additive or component partitioning.


\section{Nonlinearly Partitioned Runge--Kutta}\label{sec:nprk}
\sloppy This section briefly introduces the concept of nonlinear partitioning and nonlinearly partitioned Runge--Kutta (NPRK) methods, and directly relates them to additive and partitioned Runge--Kutta methods. Given a differential equation 
\begin{align}
	\dot{y} = G(y),
	\label{eq:classical_ode}
\end{align}
we say that a function $F(u,v$) is a {\em nonlinear partition} of $G(y)$ \cite{nprk1,nprk2} if it satisfies
\begin{align}
	F(y,y) = G(y), \text{ for all } y.	
	\label{eq:np_ode}
\end{align}
This allows us to rewrite the ODE \cref{eq:classical_ode} as
\begin{align}\label{eq:multivariable}
	\dot{y} = F(y,y).	
\end{align}
For example, $F(u,v) = uv$ is a nonlinear partition of $G(y) = y^2$; this nonlinear partition factors the nonlinearity as a product of two terms.

This reformulation facilitates the derivation of numerical methods that treat each argument of $F(y,y)$ with different levels of implicitness. One of the simplest example methods is the first-order nonlinearly partitioned IMEX Euler method,
\begin{align}
	y_{n+1} = y_n + h F(y_{n+1}, y_{n}),
\end{align}
which treats the first argument of $F(y,y)$ implicitly and the second argument explicitly. For example, this is the standard approach in thermal radiation transport dating back many decades \cite{Larsen.1988}. Nonlinear partitions have a straightforward extension to the case where $F$ has $M$ arguments \cite{nprk1,nprk2}; however, we will not consider this generalization in this paper.\\

\subsection{The NPRK Framework} NPRK methods \cite{nprk1,nprk2} are multistage integrators for solving the nonlinearly partitioned equation \eqref{eq:np_ode}, that can be written as
\begin{subequations} 
	\label{eq:nprk-general}
	\begin{align}
	    	Y_i &= y_n + h\sum_{j=1}^{s} \sum_{k=1}^{s} a_{ijk} F(Y_j,Y_k), \quad i=1,\ldots, s, \label{eq:nprk-general-a} \\
			y_{n+1} &= y_n + h\sum_{i=1}^s \sum_{j=1}^s b_{ij} F(Y_i,Y_j).	\label{eq:nprk-general-b}
	\end{align}
\end{subequations}
The classical RK matrix $a_{ij}$ and weight vector $b_j$ has been replaced with a third-order tensor $a_{ijk}$ and weight matrix $b_{ij}$. This additional generality allows one to substitute any combination of stage values ($Y_j$, $Y_k$) into the arguments of $F$. 

Every NPRK method can be interpreted as a nonlinear combination of two classical Runge--Kutta methods, called {\em underlying methods}.
\begin{definition}
	A classical RK integrator $(a_{ij}, b_i, c_i)$ is an \emph{underlying RK method} of an NPRK integrator, if \eqref{eq:nprk-general} reduces to the RK method with coefficients $(a_{ij}, b_i, c_i)$ when the partition $F(u,v)$ depends only on a single argument.
\end{definition}
If ${F(u,v) = G(u)}$, then \eqref{eq:nprk-general} reduces to the classical RK method
\begin{align}
        a^{\{1\}}_{ij} = \sum_{k=1}^s a_{ijk}, \quad  
        b^{\{1\}}_{i} = \sum_{j=1}^s b_{ij}, \quad
        c^{\{1\}}_i = \sum_{j=1}^{s} a^{\{1\}}_{ij}.
    \label{eq:nprk-underlying-rk1}
\end{align}
Similarly, if ${F(u,v) = G(v)}$, then \eqref{eq:nprk-general} reduces to
\begin{align}
        a^{\{2\}}_{ik} = \sum_{j=1}^s a_{ijk}, \quad
        b^{\{2\}}_{j} = \sum_{i=1}^s b_{ij}, \quad
        c^{\{2\}}_i = \sum_{k=1}^{s} a^{\{2\}}_{ik} = c^{\{1\}}_i.
    \label{eq:nprk-underlying-rk2}
\end{align}
One possible approach for constructing NPRK methods is to first specify the underlying RK methods and then enforce the additional linear and nonlinear coupling conditions in the NPRK tensor. \\[1em]

\subsection{Connections with additive RK and partitioned RK} For appropriately selected partitions, NPRK methods reduce to additive Runge-Kutta (ARK) \cite{Kennedy.2003tv4} and partitioned Runge-Kutta (PRK) methods \cite{Hairer.1993}[Section II.15]. To obtain an ARK method, one must select a nonlinear partition of the form
\begin{align}
	F(u,v) &= G^{\{1\}}(u) + G^{\{2\}}(v).
\end{align}
Similarly, to obtain a PRK method we first express $y$ in components as $y = [u~ v]^T$ and similarly express $G(y)$ in components as $ G(y) = [G_1(u,v) ~ G_2(u,v)]^T$. Then, if one selects the nonlinear partition
\begin{align}
	F\left(y_1, y_2\right) &= \left[ \begin{array}{l} G_1(u_1,v_1) \\ G_2(u_2,v_2) \end{array} \right] 
		\quad \text{with} \quad
		\begin{array}{l}
			y_1 = (u_1\ v_1)^T \\
			y_2 = (u_2\ v_2)^T
		\end{array}.
\end{align}
The NPRK method reduces to a PRK method. The two coefficient matrices of the PRK or ARK methods are those of the underlying methods \cref{eq:nprk-underlying-rk1,eq:nprk-underlying-rk2}. The connection between NPRK, PRK, and ARK implies that all the multirate NPRK methods developed in this work immediately double as multirate PRK and ARK methods.


\section{Multirate NPRK methods}\label{sec:mr-nprk}

Before defining multirate NPRK methods, we first discuss the reducibility of NPRK integrators and their underlying methods, since these concepts are essential for defining multirate NPRK integrators. 
\begin{definition}
	A Runge--Kutta method with tableau $M$ is irreducible if it contains no equivalent stages and if no set of stages can be removed without altering the output of the method \cite{butcher2016numerical}[Def 381E]. Otherwise, the RK method is said to be reducible; we denote the corresponding reduced tableau with $\widehat{M}$.
\end{definition}

This definition applies to NPRK methods, ARK methods, and classical RK methods. We remark on an important  point regarding NPRK methods that is essential for multirate methods

\begin{remark}
	An irreducible NPRK method may have reducible underlying integrators; \Cref{ex:first-order-multirate} below is one such simple example.
\end{remark}

This fact allows us to construct an irreducible NPRK method whose reduced underlying methods have a different number of stages; as we will see, this is precisely the starting point for constructing a multirate NPRK integrator. For the remainder of this paper, we use the following notation for an NPRK method and its underlying methods:
\begin{align}
		\begin{aligned}
            &\text{NPRK method:} &  M &:= (a_{ijk}, b_{ij}), \\
			&\text{Unreduced underlying integrators:} & M^{\{r\}} &:=(a_{ij}^{\{r\}}, b_{i}^{\{r\}}, c_{i}^{\{r\}}) , \\
			&\text{Reduced underlying integrators:}  & \widehat{M}^{\{r\}} &:=(\widehat{a}_{ij}^{\{r\}}, \widehat{b}_{i}^{\{r\}}, \widehat{c}_{i}^{\{r\}})	,	
		\end{aligned}
		&& 
		r \in \{1,2\}.
	\end{align}

To introduce multirate NPRK (MR-NPRK) methods, we do not require any additional generalization; instead, we simply impose conditions on the underlying methods of an NPRK integrator. 
\begin{definition}[Multirate NPRK Method]\label{def:multirate-nprk}
    A multirate NPRK method is an NPRK method whose {\bf reduced} underlying methods $\widehat{M}^{\{1\}}$ and $\widehat{M}^{\{2\}}$ possess a differing number of stages.
\end{definition}

As a preliminary example, we present a first-order, IMEX MR-NPRK method that integrates the first argument with one step of backward Euler with stepsize $h$ and the second argument with three steps of explicit Euler with stepsize $h/3$.

\begin{example}\label{ex:first-order-multirate} Consider the MR-NPRK method
\begin{align}
	\begin{aligned}
		Y_{1} &= y_n, \\
		Y_{2} &= y_n + \tfrac{h}{3} F(Y_1, Y_1), \\
		Y_{3} &= y_n + \tfrac{h}{3} F(Y_1, Y_1) + \tfrac{h}{3} F(Y_1, Y_2), \\
		Y_{4} &= y_n + \tfrac{h}{3} F(Y_1, Y_1) + \tfrac{h}{3} F(Y_1, Y_2) - \tfrac{2h}{3}F(Y_1, Y_3) + h F(Y_4, Y_3), \\
		y_{n+1} &= Y_{4}.
	\end{aligned}
	\label{eq:first-order-multirate-example}
\end{align}
whose unreduced and reduced underlying methods are
\begin{align}
	\begin{tabular}{cccc}
		$M^{\{1\}}$ & $M^{\{2\}}$ & $\widehat{M}^{\{1\}}$ & $\widehat{M}^{\{2\}}$ \\[0.5em]
		\renewcommand*{\arraystretch}{1.25}
		\begin{tabular}{c|cccc}
			$0$ & $0$ \\
			$\tfrac{1}{3}$ & $\tfrac{1}{3}$ & $0$ \\
			$\tfrac{2}{3}$ & $\tfrac{2}{3}$ & $0$ & $0$ \\ 
			$1$ & $0$ & $0$ & $0$ & $1$ \\ \hline
			  & $0$ & $0$ & $0$ & $1$
		\end{tabular}
		&
		\renewcommand*{\arraystretch}{1.25}
		\begin{tabular}{c|cccc}
			$0$ & $0$ \\
			$\tfrac{1}{3}$ & $\tfrac{1}{3}$ & $0$ \\
			$\tfrac{2}{3}$ & $\tfrac{1}{3}$ & $\tfrac{1}{3}$ & $0$ & \\ 
			$1$ & $\tfrac{1}{3}$ & $\tfrac{1}{3}$ & $\tfrac{1}{3}$ & $0$ \\ \hline
		 	& $\tfrac{1}{3}$ & $\tfrac{1}{3}$ & $\tfrac{1}{3}$ & $0$
		\end{tabular}
		&
		\begin{tabular}{c|cccc}
		$0$ & $1$ \\ \hline
		  & $1$
		\end{tabular}
		&
		\renewcommand*{\arraystretch}{1.25}
		\begin{tabular}{c|cccc}
			$0$ & $0$ \\
			$\tfrac{1}{3}$ & $\tfrac{1}{3}$ & $0$ \\
			$\tfrac{2}{3}$ & $\tfrac{1}{3}$ & $\tfrac{1}{3}$ & $0$ \\ \hline
			 & $\tfrac{1}{3}$ & $\tfrac{1}{3}$ & $\tfrac{1}{3}$
		\end{tabular}
	\end{tabular}	
\end{align}
It follows that \cref{eq:first-order-multirate-example} reduces to backward Euler for the partition $F(u,v)=G(u)$, and three steps of forward Euler for the partition $F(u,v)=G(v)$.
\end{example}

To simplify the presentation of multirate methods, we also introduce two sets that identify the irreducible stages of an MR-NPRK method and its underlying methods. In \Cref{sec:time-scale-coupling,sec:mrnprk-order-conditions}, we will show how these sets simplify method classification and order conditions for MR-NPRK methods. 

\begin{definition}[Stage Index Set]
    The stage index set of an $s$-stage MR-NPRK method is $\mathcal{S} := \{1,\dots,s\}. $
\end{definition}

\begin{definition}[Irreducible Stage Index Sets]
	For an MR-NPRK method $M$, its irreducible stage index sets $\Scmp{r}$, $r=1,2$, contain the indices corresponding to the irreducible stages of the $r^{th}$ unreduced, underlying method $\Mcmp{r}$:
	\begin{align}
	\Scmp{r} := \{ \nu \in \mathcal{S} : \textup{$Y_\nu$ is an irreducible stage of $\Mcmp{r}$} \}.	
	\end{align}
	The cardinality of each set is denoted as $\scmp{r} := |\Scmp{r}|$. Using this definition, a multirate NPRK method (\Cref{def:multirate-nprk}) can be equivalently defined as an NPRK method with $\scmp{1} \neq \scmp{2}$.
\end{definition}

\begin{example}
	The stage sets for the method \eqref{eq:first-order-multirate-example} are $\mathcal{S} = \{1, 2, 3, 4\}$, $\Scmp{1}=\{4\}$ and $\Scmp{2}=\{1,2,3\}$ with $\scmp{1}=1$ and $\scmp{2}=3$.	
\end{example}

In general, it is possible to construct irreducible NPRK methods with stages that are simultaneously reducible in both $M^{\{1\}}$ and $M^{\{2\}}$ (an example is shown in \Cref{sup:nprk-example-violating-stage-union}). For simplicity, we will not consider such methods and make the following assumption.

\begin{assumption}\label{assump:stage-union}
    We only consider methods where $\mathcal{S}^{\{1\}} \cup \mathcal{S}^{\{2\}} = \mathcal{S}$. 
\end{assumption}

\subsection{Timescale coupling}
\label{sec:time-scale-coupling}

The MR-NPRK framework naturally incorporates different timescales for the first and second arguments of the nonlinear partition $F(u,v)$. Here, we describe four ways information can be exchanged between these two scales. 

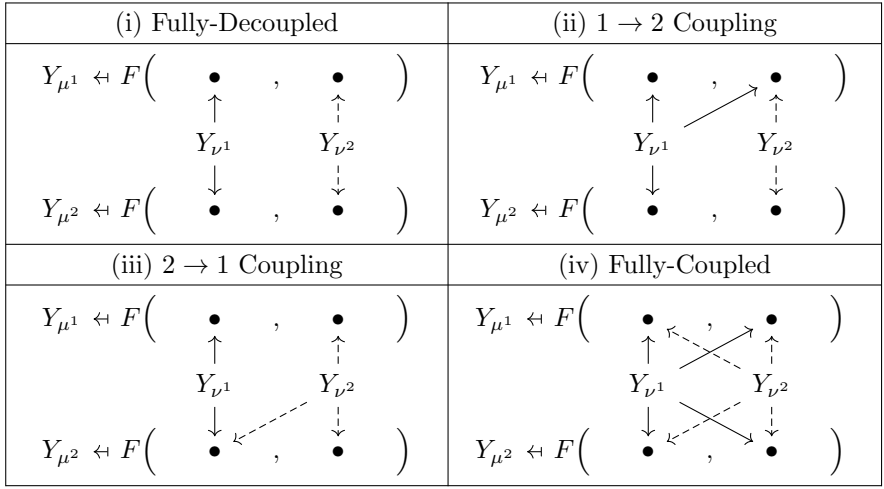
\begin{figure}[H]

\usetikzlibrary{decorations.pathmorphing}

\begin{center}

\renewcommand*{\arraystretch}{1.25}
\begin{tabular}{ |c|c| } 
 \hline
 (i) Fully-Decoupled & (ii) $1 \to 2$ Coupling \\ \hline 
 $$ \begin{tikzcd}[sep=tiny]
	{Y_{\mu^1}} & {F\Big(} & \bullet & {,} & \bullet & {\Big)} \\
	&& {Y_{\nu^1}} && {Y_{\nu^2}} \\
	{Y_{\mu^2}} & {F\Big(} & \bullet & {,} & \bullet & {\Big)}
	\arrow[maps to, from=1-2, to=1-1]
	\arrow[from=2-3, to=1-3]
	\arrow[from=2-3, to=3-3]
	\arrow[dashed, from=2-5, to=1-5]
	\arrow[dashed, from=2-5, to=3-5]
	\arrow[maps to, from=3-2, to=3-1]
\end{tikzcd} $$ & $$\begin{tikzcd}[sep=tiny]
	{Y_{\mu^1}} & {F\Big(} & \bullet & {,} & \bullet & {\Big)} \\
	&& {Y_{\nu^1}} && {Y_{\nu^2}} \\
	{Y_{\mu^2}} & {F\Big(} & \bullet & {,} & \bullet & {\Big)}
	\arrow[maps to, from=1-2, to=1-1]
	\arrow[from=2-3, to=1-3]
	\arrow[from=2-3, to=1-5]
	\arrow[from=2-3, to=3-3]
	\arrow[dashed, from=2-5, to=1-5]
	\arrow[dashed, from=2-5, to=3-5]
	\arrow[maps to, from=3-2, to=3-1] \end{tikzcd}$$  \\ 
 \hline
 (iii) $2 \to 1$ Coupling & (iv) Fully-Coupled \\ \hline 
$$\begin{tikzcd}[sep=tiny]
	{Y_{\mu^1}} & {F\Big(} & \bullet & {,} & \bullet & {\Big)} \\
	&& {Y_{\nu^1}} && {Y_{\nu^2}} \\
	{Y_{\mu^2}} & {F\Big(} & \bullet & {,} & \bullet & {\Big)}
	\arrow[maps to, from=1-2, to=1-1]
	\arrow[from=2-3, to=1-3]
	\arrow[from=2-3, to=3-3]
	\arrow[dashed, from=2-5, to=1-5]
	\arrow[dashed, from=2-5, to=3-3]
	\arrow[dashed, from=2-5, to=3-5]
	\arrow[maps to, from=3-2, to=3-1]
\end{tikzcd}$$ & \begin{tikzcd}[sep=tiny]
	{Y_{\mu^1}} & {F\Big(} & \bullet & {,} & \bullet & {\Big)} \\
	&& {Y_{\nu^1}} && {Y_{\nu^2}} \\
	{Y_{\mu^2}} & {F\Big(} & \bullet & {,} & \bullet & {\Big)}
	\arrow[maps to, from=1-2, to=1-1]
	\arrow[from=2-3, to=1-3]
	\arrow[from=2-3, to=3-3]
    \arrow[from=2-3, to=1-5]
    \arrow[from=2-3, to=3-5]
	\arrow[dashed, from=2-5, to=1-5]
	\arrow[dashed, from=2-5, to=3-3]
	\arrow[dashed, from=2-5, to=3-5]
    \arrow[dashed, from=2-5, to=1-3]
	\arrow[maps to, from=3-2, to=3-1]
\end{tikzcd} $$ \\ 
 \hline
\end{tabular}
\end{center}

\caption{Schematic of the allowed function argument evaluations for the four NPRK timescale couplings. The variables $\mu^1,\nu^1 \in \Scmp{1}$, $\mu^2, \nu^2 \in \Scmp{2}$, and $ Y_\mu \mapsfrom F(\cdot,\cdot)$ denotes function evaluations used in the stage update for $Y_\mu$.}\label{fig:timescale-coupling-diagrams}
\end{figure}

\begin{enumerate}[label = (\roman*)]
	\item MR-NPRK with {\bf fully-decoupled timescales}. These methods have non-intersecting irreducible stage sets that keep the two timescales separated. Specifically, only stages with index in $\mathcal{S}^{\{1\}}$ and stages with index in $\mathcal{S}^{\{2\}}$ can be evaluated in the first and second arguments of $F(\cdot,\cdot)$, respectively. This requirement applies to function evaluations in the stages and the output.
	\item MR-NPRK with {\bf 1 $\to$ 2 coupling}. These methods relax one of the requirements of fully-decoupled MR-NPRK. Specifically, there is at least one stage $Y_i$, $i \in \Scmp{1}$, with a non-zero coefficient  $a_{iuv}$ for $\mu, \nu \in \mathcal{S}^{\{1\}}$ -- this implies that there is an evaluation $F(Y_\mu, Y_\nu)$ in the $i^{th}$ stage computation. Intuitively, in stages pertaining to the first timescale, information from the first timescale may be transferred into the second.
	
	\item MR-NPRK with {\bf 2 $\to$ 1 coupling}. These are obtained by relaxing a similar requirement of fully-decoupled MR-NPRK methods in the opposite direction. Specifically, there  is at least one stage $Y_i$, $i \in \Scmp{2}$, with a non-zero coefficient $a_{iuv}$ for $\mu, \nu \in \mathcal{S}^{\{2\}}$. Intuitively, in stages pertaining to the second timescale, information from the second timescale may be transferred into the first. 
 		
	\item MR-NPRK with {\bf fully-coupled timescales.} This family of methods contains all of the remaining possibilities. Intuitively, the information between the two timescales is freely exchanged. 
		
\end{enumerate}
The allowed function argument evaluations for the above four types of MR-NPRK methods are shown schematically in \Cref{fig:timescale-coupling-diagrams}. 

\begin{example}
The MR-NPRK method \eqref{eq:first-order-multirate-example} is fully-coupled since the stage $Y_1$ belongs to $\Scmp{2}$ but it is substituted into the first argument of $F(\cdot,\cdot)$ in the fourth stage. Since $4 \in \Scmp{1}$, this violates the allowed function argument evaluations for fully-decoupled, $1 \to 2$, and $2 \to 1$ coupled methods.
\end{example}

\subsection{Order conditions}
\label{sec:mrnprk-order-conditions}

Necessary and sufficient order conditions for NPRK methods were derived in \cite{nprk2} using edge-colored rooted trees. Order conditions for MR-NPRK methods are identical since MR-NPRK methods are a subclass of NPRK methods. Nevertheless, the reducibility of the underlying methods of MR-NPRK often leads to simplifications in the summed indices found in each order condition. Moreover, the timescale coupling of a method influences the degree of simplifications, with fully-decoupled methods experiencing the greatest simplifications.

\Cref{fig:timescale-coupling-diagrams} described the allowable non-zero stage coefficients for MR-NPRK methods with different types of timescale coupling. This leads to implied sparsity in the MR-NPRK coefficient tensors that can be used to simplify order conditions; this sparsity is described in \Cref{tab:sparsity}.

\begin{table}[h]
	\caption{Sparsity of multirate NPRK methods with specified timescale couplings - proof in \Cref{subsec:implied-sparsity}.}
	\label{tab:sparsity}
	\begin{center}
	\renewcommand*{\arraystretch}{1.25}
	\begin{tabular}{l|ccc}
		 	& Implied Coefficient Conditions  \\
	   & ($\nu \in \mathcal{S}^{\{1\}}$, $\gamma \in \mathcal{S}^{\{2\}}$, $i,j,k \in \mathcal{S}$) \\ \hline
		 	Fully-decoupled & $a_{i \gamma k} = 0$, $a_{ij\nu}=0$, and $b_{\gamma k} = 0 = b_{j\nu}$ \\
		 $1 \to 2$ coupling & $a_{i \gamma k} = 0$, $a_{\gamma j \nu} = 0$, and $b_{\gamma k} = 0 = b_{j\nu}$ \\
			$2 \to 1$ coupling & $a_{\nu \gamma k} = 0$, $a_{i j \nu} = 0$, and $b_{\gamma k} = 0 = b_{j\nu}$
		 \end{tabular}
	\end{center}
\end{table}

Since order conditions involve sums over a method's coefficients, the sparsity described in \Cref{tab:sparsity} enables one to replace sums over $\mathcal{S}$ with sums over smaller subsets. \Cref{tab:mr-nprk-order-conditions} contains order conditions for MR-NPRK methods up to order three, although a similar reduction holds for arbitrary order (see \Cref{sup:order-conditions}).

\begin{table}[h!]
    
	\begin{center}
				
		{\bf Index Sets $I_1(i)$, $I_2(i)$ for $a_{ijk}$} \\[1em]
		\renewcommand*{\arraystretch}{1.25}	
            \scalebox{0.9} {
		\begin{tabular}{r|llll}
    		& Fully-decoupled & $1 \to 2$ Coupled & $2 \to 1$ Coupled & Fully-coupled \\ \hline
    		$I_1(i)$ & $\Scmp{1}$ & $\Scmp{1}$ & $\begin{cases} \mathcal{S} \text{ if } i \in \Scmp{2} \\ \Scmp{1} \text{ if } i \in \Scmp{1} \end{cases}$ & $\mathcal{S}$ \\
    		$I_2(i)$ & $\Scmp{2}$ & $\begin{cases} \mathcal{S} \text{ if } i \in \Scmp{1} \\ \Scmp{2} \text{ if } i \in \Scmp{2} \end{cases}$  & $\Scmp{2}$	& $\mathcal{S}$		
    	\end{tabular}
            }
     
    	\vspace{1em}
    	{\bf Order Conditions} \\ [1em]
    	
		\renewcommand*{\arraystretch}{1.5}	
            \scalebox{0.8}{
		\begin{tabular}{c|c|ll|l}
			Tree 	& Elementary 	& Elementary 			& $\Phi(\tau)$ in terms of 								& $\gamma(\tau)$ \\[-0.5em]
			$\tau$ 	& Differential 	& Weight $\Phi(\tau)$ 	& \eqref{eq:nprk-underlying-rk1} and \eqref{eq:nprk-underlying-rk2}  &  \\[0.3em] \hline
			$\T10$ & $F$ 				& $\sum_{\substack{i \in \Scmp{1} \\ j \in \Scmp{2}}} b_{ij}$ 					& $= \sum_{i \in \Scmp{1}} b_i^{\{1\}} = \sum_{j \in \Scmp{2}} b_j^{\{2\}}$ & 1 \\ 
            [1.0em]\hdashline
			$\T21$ & $F_1[F]$  			& $\sum_{\substack{i \in \Scmp{1}, k \in I_1(i) \\ j \in \Scmp{2}, l \in I_2(i)}} b_{ij}a_{ikl}$ 			& $= \sum_{i \in \Scmp{1}} b_i^{\{1\}}c_i$  & 2 \\
			$\T20$ & $F_2[F]$  			& $\sum_{\substack{i \in \Scmp{1}, k \in I_1(j) \\ j \in \Scmp{2}, l \in I_2(j)}} b_{ij}a_{jkl}$ 			& $= \sum_{j \in \Scmp{2}} b_j^{\{2\}}c_j$ & 2 \\ 
            [1.0em]\hdashline
			$\T33$ & $F_{11}[F,F]$ 		& $\sum_{\substack{i \in \Scmp{1}, k,u \in I_1(i) \\ j \in \Scmp{2}, l,v \in I_2(i)}} b_{ij} a_{ikl} a_{iuv}$ 	& $= \sum_{i \in \Scmp{1}} b_i^{\{1\}} c_i c_i$ & 3 \\
			$\T32$ & $F_{12}[F,F]$ 		& $\sum_{\substack{i \in \Scmp{1}, k\in I_1(i), u \in I_1(j) \\ j \in \Scmp{2}, l \in I_2(i), v \in I_2(j)}} b_{ij} a_{ikl} a_{juv}$ 	& not possible & 3 \\
			$\T30$ & $F_{22}[F,F]$ 		& $\sum_{\substack{i \in \Scmp{1}, k,u \in I_1(j) \\ j \in \Scmp{2}, l,v \in I_2(j)}}  b_{ij} a_{jkl} a_{juv}$ 	& $= \sum_{j \in \Scmp{2}} b_j^{\{2\}} c_j c_j$ & 3 \\[1em]
			$\T43$ & $F_{1}[F_{1}[F]]$ 	& $\sum_{\substack{i \in \Scmp{1}, k \in I_1(i), u\in I_1(k) \\ j \in \Scmp{2}, l \in I_2(i) ,v \in I_2(k)}} b_{ij} a_{ikl} a_{kuv}$ 	& $= \sum_{\substack{i \in \Scmp{1} \\ k \in I_1(1) }} b_i^{\{1\}} a^{\{1\}}_{ik} c_k$ & 6\\[1em]
			$\T42$ & $F_{1}[F_{2}[F]]$ 	& $\sum_{\substack{i \in \Scmp{1}, k \in I_1(i), u \in I_1(l) \\ j \in \Scmp{2}, l \in I_2(i), v \in I_2(l)}} b_{ij} a_{ikl} a_{luv}$ 	& $= \sum_{\substack{i \in \Scmp{1} \\ l \in I_2(i)}} b_i^{\{1\}} a^{\{2\}}_{il} c_l$ & 6\\[1em]
			$\T41$ & $F_{2}[F_{1}[F]]$ 	& $\sum_{\substack{i \in \Scmp{1}, k \in I_1(j), u \in I_1(k) \\ j \in \Scmp{2}, l \in I_2(j), v \in I_2(k)}} b_{ij} a_{jkl} a_{kuv}$ 	& $= \sum_{\substack{j \in \Scmp{2} \\ k \in I_1(j)}} b_j^{\{2\}} a^{\{1\}}_{jk} c_k$ & 6\\[1em]
			$\T40$ & $F_{2}[F_{2}[F]]$ 	& $\sum_{\substack{i \in \Scmp{1}, k \in I_1(j), u \in I_1(l) \\ j \in \Scmp{2}, l \in I_2(j), v \in I_2(l)}} b_{ij} a_{jkl} a_{luv}$ 	& $= \sum_{\substack{j \in \Scmp{2} \\ l \in I_2(j) }} b_j^{\{2\}} a^{\{2\}}_{jl} c_l$ & 6\\
		\end{tabular}
            }
	\end{center}

	\caption{Order conditions are $\Phi(\tau) = 1/\gamma(\tau)$. The horizontal dashed lines separate conditions of different orders. The order conditions consist of the well-known ARK order conditions, along with an additional condition at third order (and higher) corresponding to nonlinear coupling between the arguments of $F$.}
		\label{tab:mr-nprk-order-conditions}
\end{table}

\subsection{Linear stability}\label{sec:linear-stability}

Linear stability analysis \cite{Hairer.1993}[Section IV.2] studies the stability of a time integration method applied to the Dahlquist test problem $y'=\lambda y$, $\lambda \in \mathbb{C}$, and is useful for determining the types of equations (diffusive or oscillatory) for which an integrator is stable. In this short section, we review the linear stability analysis for NPRK methods proposed in \cite{nprk1}, then discuss general stability properties of MR-NPRK methods. In \Cref{sec:constructing}, we will discuss how certain method ansatze lead to beneficial linear stability properties.

\subsubsection{The partitioned Dalqhuist equation}
 We study linear stability of NPRK methods using the partitioned Dahlquist equation,
	\begin{align}
		y' = F(y,y) 
		\quad \text{for} \quad
		F(u,v) = \lambda_1 u + \lambda_2 v.
		\label{eq:dahlquist-partitioned}	
	\end{align}

	This equation has been used to study the stability of ARK methods \cite{Ascher.1997,Izzo.2017,Sandu.2015}; however, it also approximates the nonlinearly partitioned system  ${y'=F(y,y)}$ if one linearizes in each argument and assumes that the resulting component Jacobians are simultaneously diagonalizable. 

	\subsubsection{Stability functions}
	Since \eqref{eq:dahlquist-partitioned} is an additive partition, an NPRK method applied to this problem reduces to its underlying ARK method with tableaux \eqref{eq:nprk-underlying-rk1} and \eqref{eq:nprk-underlying-rk2}. Therefore, the linear stability function of an NPRK method is the well-known linear stability function for ARK methods \cite{Kennedy.2003tv4}. 
	Specifically, an NPRK method applied to \eqref{eq:dahlquist-partitioned} results in the iteration 
\begin{align}
	y_{n+1} = R(z_1,z_2), \quad z_1 = h \lambda_1, \quad z_2 = h \lambda_2,	
	\label{eq:stability-iteration}
\end{align}
with stability function
\begin{align}
	\begin{aligned}
		R(z_1,z_2) &= \frac{P(z_1,z_2)}{Q(z_1,z_2)} = \frac{\det(\mathbf{I} - (z_1 \mathbf{A}^{\{1\}} + z_2 \mathbf{A}^{\{2\}}) + \mathbf{e}(z_1 \mathbf{b}^{\{1\}} + z_2 \mathbf{b}^{\{2\}})^\text{T} )}{\det(\mathbf{I} - z_1 \mathbf{A}^{\{1\}} - z_2 \mathbf{A}^{\{2\}})},
	\end{aligned}
	\label{eq:ark-stability-function}
\end{align}
where
$\mathbf{A}^{\{1\}}=[a^{\{1\}}_{ij}]$, $\mathbf{b}^{\{1\}}=[b^{\{1\}}_{i}]$, $\mathbf{A}^{\{2\}}=[a^{\{2\}}_{ij}]$, and $\mathbf{b}^{\{1\}}=[b^{\{2\}}_{i}]$ are the coefficient matrices and vectors for the underlying ARK methods defined in \eqref{eq:nprk-underlying-rk1} and \eqref{eq:nprk-underlying-rk2}. When $z_1=0$ or $z_2=0$, the stability function \eqref{eq:ark-stability-function} simplifies to the stability function of the first or second underlying method, respectively.

\subsubsection{Stability region definitions and visualization}
The traditional linear stability region for an ARK method is defined as 
\begin{align}
	\Gamma := \left\{ (z_1,z_2) \in \mathbb{C} \times \mathbb{C} : \left| R(z_1,z_2) \right| \le 1 \right\},
	\label{eq:ark-stability-region}	
\end{align}
which contains all $(z_1,z_2)$ pairs such that \eqref{eq:stability-iteration} remains bounded. When considering stability for PDE discretizations, the values $z_1$ and $z_2$ represent the eigenvalues of the linearized component Jacobians. Since real-valued PDE discretizations have spectrums that are symmetric about the imaginary axis, one desires stability for both $z_1$ and $z_1^*$. Therefore, it is more realistic to consider the more restrictive stability region
\begin{align}
	\mathcal{P} := \left\{ (z_1,z_2) \in \mathbb{C} \times \mathbb{C} : \max\left(\left| R(z_1,z_2) \right|,\left| R(z_1^*,z_2) \right|\right) \le 1 \right\} \subset \Gamma.
	\label{eq:ark-stability-region-symmetric}	
\end{align}

The regions \cref{eq:ark-stability-region,eq:ark-stability-region-symmetric} are both four-dimensional and cannot be visualized directly. In this paper, we will overlay two-dimensional slices of $\mathcal{P}$ pertaining to a fixed $z_1$ value,
	\begin{align}
		\mathcal{P}(z_1) &:= \left\{ z_2 \in \mathbb{C} : \max 
	\left( 
		\left| R(z_1, z_2) \right|,
		\left| R(z_1^*, z_2) \right|
	\right) \le 1 \right\}.
	\label{eq:stability-region-symmetric-slice} 
	\end{align}
	  When $\arg(z_1) = \pi/2$, these regions approximate stability for a skew-symmetric advection discretization, while $\arg(z_1)=\pi$ approximately represents a symmetric positive-definite diffusion discretization. An intermediate value, $\arg(z_1) \in (\pi/2, \pi)$, approximates a mix of advection and diffusion.

\section{Constructing IMEX MR-NPRK methods}
\label{sec:constructing}

For the remainder of this paper, we focus our attention on diagonally-implicit IMEX MR-NPRK methods that treat the first component implicitly and the second component explicitly. This implies the following sparsity on the NPRK tensor
\begin{align}
	a_{ijk} = 0 \text{ for } j > i, k \ge i.
\end{align}
We also assume that $\scmp{1} < \scmp{2}$ such that the first argument of $F(y,y)$ will be integrated using fewer stages than the second. We begin by discussing the linear stability of such methods and then introduce MR-NPRK methods of order two and three with favorable stability properties.

\subsection{Linear stability considerations}\label{sec:constructing:stab}

The stability function of an A-stable NPRK method satisfies ${|R(z_1,z_2)|\le 1}$ for all $z_1,z_2 \in \mathbb{C}^-$. All IMEX NPRK methods are only implicit in the first variable and therefore cannot be A-stable since $R(0,z_2)$ is the stability function of an explicit RK method. However, we can still seek methods that guarantee stability as the implicit variable $z_1$ becomes increasingly stiff.
\begin{definition} An NPRK method is L-stable in the stiff $z_1$ limit if
\begin{align}	
	\label{eq:l-stability-z1-limit}
	\lim_{|z_1|\to \infty} R(z_1,z_2) = 0.
\end{align}
\end{definition}

Intuitively, IMEX NPRK methods that are L-stable in the stiff $z_1$ limit have larger stability regions since $\lim_{|z_1| \to \infty} \mathcal{P}(z_1)$ from \eqref{eq:stability-region-symmetric-slice} is the entire  $z_2$ plane. However, we note that \eqref{eq:l-stability-z1-limit} does not guarantee that the area of $P(z_1)$ increases monotonically as one transitions from the finite stability domain $\mathcal{P}(0)$ to the infinite stability domain $\mathcal{P}(\infty)$. 

To facilitate the construction of methods that are L-stable in the stiff $z_1$ limit, we present the following theorem.
\begin{theorem}
	\label{thm:z1-L-stable-ansatz}
	An IMEX MR-NPRK method is L-stable in the stiff $z_1$ limit if it satisfies the following the properties:
	\begin{enumerate}[(i)]
		\item fully-decoupled or $1 \to 2$ coupled,
		\item stiffly-accurate,
		\item $\cmp{\widehat{M}}{1}$ (i.e., the reduced $\Mcmp{1}$) is a DIRK method with non-zero diagonals.		
	\end{enumerate}	
\end{theorem}
\begin{proof}
    See \Cref{appendix:stability-analysis}. \qed
\end{proof}

\subsubsection{The practical utility of L-stability for IMEX MR-NPRK}
\label{subsec:practical-utility-of-l-stability}

\Cref{thm:z1-L-stable-ansatz} provides a simple ansatz for constructing an L-stable method in the stiff $z_1$ limit. For single-rate IMEX methods, this is a highly desirable property that helps prevent numerically-induced oscillations on stiff problems. For IMEX MR-NPRK methods, the practical utility is more subtle. The stability function numerator is a polynomial of degree at least $\scmp{2}$ in $z_2$ (See \Cref{prop:stability-numerator}) that is being damped by a polynomial of degree $\scmp{1}$ in $z_1$ (See \Cref{prop:stability-denominator}). Thus, asymptotically speaking, even a fully-decoupled method will only have damping if $|z_1| \gg |z_2|^{\scmp{2}/\scmp{1}}.$
In practice, this implies that, barring an extreme separation of scales between the implicit and explicit processes, methods with $\scmp{2} \gg \scmp{1}$ must integrate the explicit process using a timestep that ensures stability for the underlying explicit method. Of course, taking $\scmp{2}$ large also implies a larger explicit stability region; however, this only scales linearly in $\scmp{2}$. In summary, the practical utility of L-stability decreases as $\scmp{2}$ increases.

We demonstrate these phenomena by considering two simple first-order IMEX MR-NPRK methods. The first generalizes \eqref{eq:first-order-multirate-example}, such that $\widehat{M
}^{\{2\}}$ is now $\scmp{2}$ steps of forward Euler with stepsize $h/\scmp{2}$ and $\cmp{\widehat{M}}{1}$ remains one step of backward Euler:
\begin{align}
	\begin{aligned}
		Y_{1} &= y_n, \\
		Y_{i} &= y_n + \tfrac{h}{\scmp{2}} \sum_{j=1}^{i-1} F(Y_1, Y_j), \quad i = 2, \ldots, \scmp{2}, \\
		Y_{\scmp{2}+1} &= y_n + \tfrac{h}{\scmp{2}} \sum_{j=1}^{\scmp{2}-1} F(Y_1, Y_j) -\tfrac{h(\scmp{2} -1)}{\scmp{2}} F(Y_1,Y_j) + hF(Y_{\scmp{2}+1}, Y_{\scmp{2}}), \\			
		y_{n+1} &= Y_{\scmp{2}+1}.
	\end{aligned}
	\label{eq:mr-nprk-example-z1-unstable}
\end{align}
This fully-coupled method is mildly unstable in the stiff $z_1$ limit, in the sense that
\begin{align}
	\lim_{z_1 \to \infty} R(z_1,z_2) = \sum_{j=1}^{\scmp{2}} c_j z_2^j 
	\quad \text{with} \quad \lim_{z_2 \to \infty} p(z_2) = \infty .
\end{align}

Using \Cref{thm:z1-L-stable-ansatz}, we construct the following first-order fully-decoupled method that is L-stable in the stiff-$z_1$ limit:
\begin{align}
	\begin{aligned}
		Y_{1} &= y_n, \\
		Y_{i} &= y_n + \tfrac{h}{\scmp{2}} \sum_{j=1}^{i-1} F(Y_2, Y_j), \quad i = 2, \ldots, \scmp{2}, \\
		Y_{\scmp{2}+1} &= y_n + \tfrac{h}{\scmp{2}} \sum_{j=1}^{\scmp{2}-1} F(Y_2, Y_j) + h(\tfrac{1}{\scmp{2}} - \tfrac{3}{4} ) F(Y_2,Y_{\scmp{2}}) + h \tfrac{3}{4} F(Y_{\scmp{2}+1}, Y_{\scmp{2}}), \\	
		y_{n+1} &= Y_{s},
	\end{aligned}
	\label{eq:mr-nprk-example-z1-lstable}
\end{align}
where $\widehat{M}^{\{2\}}$ remains unchanged while $\cmp{\widehat{M}}{1}$ is now a two-stage L-stable DIRK method. Subsequently, in \Cref{subsec:method-ansatze}, we construct second and third order IMEX MR-NPRK methods that are also L-stable in the stiff-$z_1$ limit.

\begin{figure}[h]

	\centering
	
	\begin{tabular}{cccc}
		& $\scmp{2}=4$ & $\scmp{2}=8$ & $\scmp{2}=32$ \\[1em]
		\rotatebox{90}{\hspace{3em} method \eqref{eq:mr-nprk-example-z1-unstable}}&
		\includegraphics[align=b,width=0.28\textwidth]{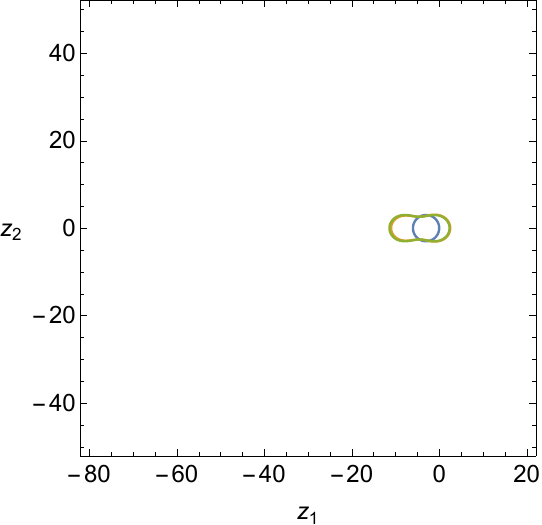} &
		\includegraphics[align=b,width=0.28\textwidth]{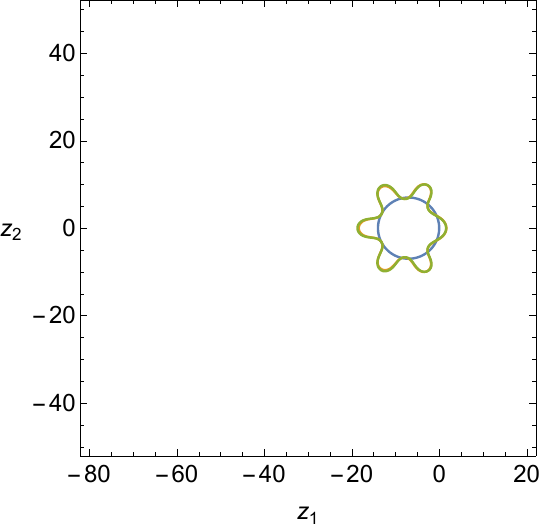} &
		\includegraphics[align=b,width=0.28\textwidth]{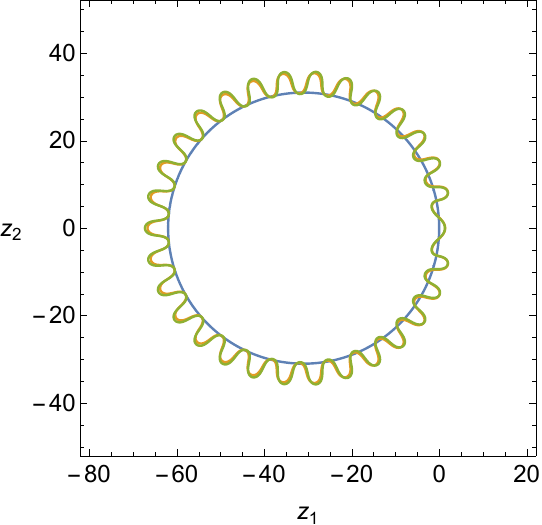} \\		
		\rotatebox{90}{\hspace{3em} method \eqref{eq:mr-nprk-example-z1-lstable}} & 
		\includegraphics[align=b,width=0.28\textwidth]{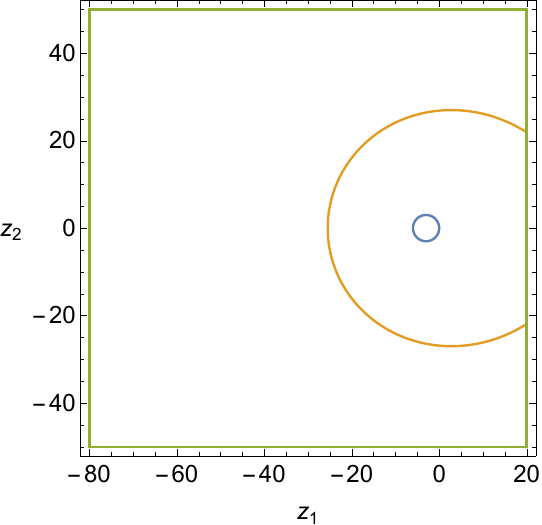} &
		\includegraphics[align=b,width=0.28\textwidth]{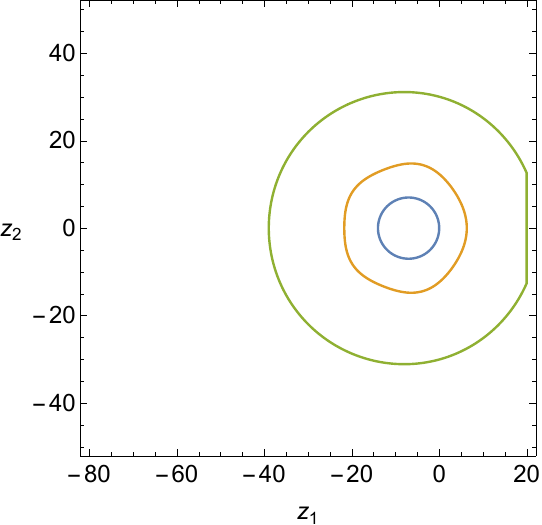} &
		\includegraphics[align=b,width=0.28\textwidth]{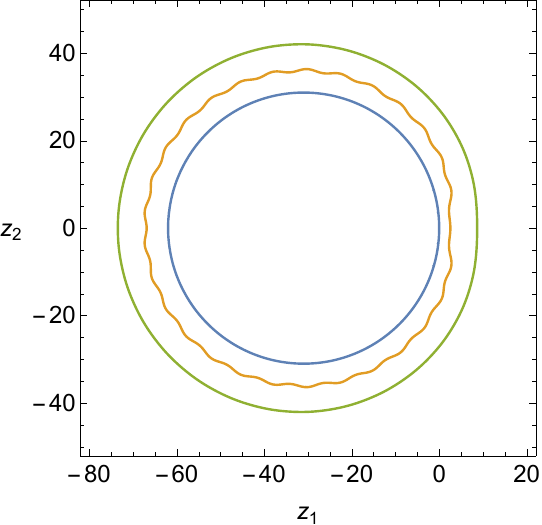}
	\end{tabular}
	
	\begin{tabular}{c}
		{\tiny \textcolor{mathematica_plot_blue}{\hdashrule[0.2ex]{2em}{2pt}{}} $z_1 = 0$ \hspace{1em}}
		{\tiny \textcolor{mathematica_plot_orange}{\hdashrule[0.2ex]{2em}{2pt}{}} $z_1 = -100$ \hspace{1em}}
		{\tiny \textcolor{mathematica_plot_green}{\hdashrule[0.2ex]{2em}{2pt}{}} $z_1 = -10000$}
	\end{tabular}

	\caption{Stability slices $\mathcal{P}(z_1)$, defined in \eqref{eq:stability-region-symmetric-slice}, for the MR-NPRK methods \eqref{eq:mr-nprk-example-z1-unstable} and \eqref{eq:mr-nprk-example-z1-lstable}. The first method is unstable in the stiff $z_1$ limit while the second is L-stable. Stability improvements due to L-stability in the stiff $z_1$ limit diminish as $\scmp{2}$ increases.}
	\label{fig:mr-lstab-importance}
	
\end{figure}

In \Cref{fig:mr-lstab-importance}, we plot slices $\mathcal{P}(z_1)$, defined in \eqref{eq:stability-region-symmetric-slice}, pertaining to the stability regions of the two methods  \cref{eq:mr-nprk-example-z1-unstable,eq:mr-nprk-example-z1-lstable}. We take $z_1 \in \mathbb{R}^-$ to represent the case where the first argument is a diffusive operator. For $\scmp{2} = 4$, \eqref{eq:mr-nprk-example-z1-lstable} provides significantly larger stability regions, particularly for large $|z_1|$ that are typical of even-order discretized differential operators. As $\scmp{2}$ increases, we see that the explicit stability region $\mathcal{P}(z_1=0)$ increases linearly in size. However, notice that there is also a decrease in damping for the method \eqref{eq:mr-nprk-example-z1-lstable} when $z_1$ is large. For example, when $\scmp{2}=32$ and $z_1=-10000$, the improvement in stability compared to \cref{eq:mr-nprk-example-z1-unstable} is minor despite $z_1$ being more than 100 times stiffer than the explicit component. Of course, for the method \eqref{eq:mr-nprk-example-z1-lstable}, $\lim_{z_1 \to \infty}\mathcal{P}(z_1)$ is the entire $z_2$ plane; however, in practice, a problem may not present sufficiently large timescale separations for this limit to represent a realistic stability region.

\subsection{Proposed IMEX MR-NPRK method ansatze}
\label{subsec:method-ansatze}

Due to the generality of the NPRK framework there are many ways to derive IMEX MR-NPRK methods. In this paper, we propose a method construction that: (i) ensures stability in the stiff $z_1$ limit following \Cref{thm:z1-L-stable-ansatz}, and (ii) places all implicit stages before and after all the  explicit stages -- we name this approach {\em implicit wrapping}. Implicit wrapping allows us to easily construct IMEX MR-NPRK methods based on any existing explicit RK method. Specifically, the proposed methods allow the user to specify any explicit method for integrating the second component of $F(y,y)$. Moreover, excluding implicit solves, the timestep computation is equivalent to integrating a modified equation using the specified explicit integrator. This final fact is difficult to see from the method formulas; we therefore also include a ``practical implementation'' and pseudocode section after each method that clearly reveal this structure. The result is a set of multirate methods with a simple, efficient, and low-memory structure, that can be easily incorporated into existing codes. Specifically, we simply require that a code can perform implicit solves involving the first argument of $F(\cdot,\cdot)$, and explicit integration of the system $w'=G(c,w)$ for fixed $c$.

\subsubsection{Second-order IMEX MR-NPRK}\label{sec:second-order-schemes}

We propose an $s$-stage second-order ansatz that pads an $\scmp{2}$-stage explicit RK method ($\ahcmp{2}_{ij}$, $\bhcmp{2}_i$) with one implicit stage at the beginning and one at end of the method. The method ansatz is
\begin{align}
	\begin{aligned}
		Y_1 &= y_n, \\
		{\color{magenta} Y_2} &= y_n + h a_{2,2,1} F({\color{magenta} Y_2},Y_1), \\
		Y_i &= y_n + h \acmp{2}_{i-1,1} F({\color{magenta} Y_2},Y_1) + h \sum_{k=3}^{i-1} \acmp{2}_{i-1,k-1} F({\color{magenta} Y_2},Y_k), \hspace{2em} i = 3,\ldots, \scmp{2}+1,\\
		{\color{Cerulean} Y_{s}} &= y_n 
			+ h \bhcmp{2}_1  F({\color{magenta} Y_2},Y_1) 
			+ h\sum^{\scmp{2}}_{k=3} \bhcmp{2}_{k-1} F({\color{magenta} Y_2},Y_k) \\ & \hspace{1em}+h\left[ (\bhcmp{2}_{\scmp{2}} - a_{s,s,s-1} ) F({\color{magenta} Y_2},Y_{\scmp{2}+1}) + a_{s,s,s-1}  F({\color{Cerulean} Y_{s}},Y_{s-1}) \right], \\
		y_{n+1} &= {\color{Cerulean} Y_{s}},
	\end{aligned}
	\label{eq:mr-nprk-second-order-minus}
\end{align}
with $s=\scmp{2}+2$. If we want the underlying implicit method to be singly-diagonally implicit (SDIRK), then second-order conditions require that
\begin{align}
	a_{2,2,1} = a_{s,s,s-1} = \frac{2 \pm \sqrt{2}}{2}.
	\label{eq:second_order_mr_nprk_gamma_conditions}
\end{align}
Selecting the negative case $(2 - \sqrt{2}) / 2$ leads to improved error while selecting the positive case $(2 + \sqrt{2}) / 2$ increases damping in the implicit method. A detailed derivation of this method can be found in  \Cref{sup:second-order-derivation}.

\subsubsection*{Practical implementation for MR-NPRK method \eqref{eq:mr-nprk-second-order-minus}}
 Given the initial condition $y_n$ and an explicit RK method $(\AhcmpB{2}, \bhcmpB{2})$ of at least second order:
\begin{enumerate}
	\item Select $a_{2,2,1}$ and $a_{s,s,s-1}$ according to \cref{eq:second_order_mr_nprk_gamma_conditions}.
	\item Solve the implicit system ${\color{magenta} Y_2} = y_n + h a_{2,2,1} F({\color{magenta}Y_2}, y_n)$.
	\item Compute one timestep of the ODE $w' = F({\color{magenta} Y_2}, w)$, $w_0=y_n$  using the explicit method $(\AhcmpB{2}_{ij}, \bhcmpB{2}_i)$. In addition to the output $w_1$, also store the final stage value $W_{\scmp{2}}$.
	\item Solve the implicit system ${\color{Cerulean} Y_s} = w_{1} + h a_{s,s,s-1} \left[ -F({\color{magenta} Y_2},W_{\scmp{2}}) + F({\color{Cerulean}Y_s}, W_{\scmp{2}}) \right]$.
\end{enumerate} 

\vspace{1em}
\begin{center}
	\begin{tabular}{|p{0.95\textwidth}|}
	  \hline
	  \vspace{-0.25em}
	  {\bf Pseudocode for MR-NPRK method \eqref{eq:mr-nprk-second-order-minus}} \hfill Select $\gamma = (\sqrt{2} \pm 2)/2$ \\[.25em] \hline         
	        
	        \vspace{-0.25em}
	        Solve ${\color{magenta} Y_2} = y_n + h \gamma F({\color{magenta} Y_2},y_n)$ \\[.5em]
	        
	        [$w_1$, $W_{\scmp{2}}$] = RK\_step($w' = F({\color{magenta} Y_2},w)$, $w_0=y_n$,  $\AhcmpB{2}$, $\bhcmpB{2}$, $h$) \\[.5em]
	        
	        Solve ${\color{Cerulean} Y_s} = w_{1} + h \gamma \left[ -F({\color{magenta} Y_2},W_{\scmp{2}}) + F({\color{Cerulean}Y_s}, W_{\scmp{2}}) \right]$ \\[.5em] \hline 
	\end{tabular}
\end{center}

\subsubsection{Third-order IMEX MR-NPRK}\label{sec:third-order-schemes}

In \Cref{tab:mr-nprk-third-order-v1}, we propose a third-order method ansatz that pads an $\scmp{2}$-stage explicit RK method ($\AhcmpB{2}$, $\bhcmpB{2}$) with two implicit stages at the beginning and one at end of the method. 
A detailed derivation of this method is contained in \Cref{sup:third-order-derivation}. The underlying implicit method, whose coefficients are denoted below using $\ahcmp{1}_{ij}$, is the L-stable, three-stage, stiffly-accurate SDIRK method \cite{butcher2016numerical}[p.274] 
	\begin{align}
			\renewcommand*{\arraystretch}{1.25}
			\begin{tabular}{l|lll}
				$\lambda$					& $\lambda$ \\
				$\tfrac{1}{2}(1 + \lambda)$	& $\tfrac{1}{2}(1 - \lambda)$ & $\lambda$ \\
				$1$  						& $\tfrac{1}{4}(-6 \lambda^2 + 16 \lambda - 1)$ & $\tfrac{1}{4}(6 \lambda^2 - 20 \lambda + 5)$ & $\lambda$ \\ \hline
				 							& $\tfrac{1}{4}(-6 \lambda^2 + 16 \lambda - 1)$ & $\tfrac{1}{4}(6 \lambda^2 - 20 \lambda + 5)$ & $\lambda$
			\end{tabular}
			\label{eq:stiffly-accurate-three-stage-SDIRK}
	\end{align}
	where $\lambda$ is the root of the polynomial $p(z) = \tfrac{1}{6} - \tfrac{3}{2} \lambda + 3 \lambda^2 - \lambda^3$, which to 32 digits of precision is $\lambda=0.43586652150845899941601945119356$. 
		
		\begin{remark}
			\label{remark:mr3-v1-unbounded-coeff}
			 If one selects an explicit method for which $\lim_{\scmp{2} \to \infty} |\bcmp{2}_{\scmp{2}}| = 0$ (true for any composite RK tableau \cref{eq:rk-subcycled}), then the coefficient $a_{\scmp{2}+2,3,\scmp{2}+1}$ in \cref{tab:mr-nprk-third-order-v1} will become unbounded as $\scmp{2} \to \infty$. Since we only considered moderately sized $\scmp{2}$ in our numerical experiments, this never posed a problem. However, for those who may want to multirate with very large $\scmp{2}$ we also propose an alternative method with stable coefficients in \cref{subsec:mr3-alternative}, but with twice the number of explicit function evaluations per step.
		\end{remark}
		
		\begin{table}[h]
		
		\begin{tabular}{|p{0.95\textwidth}|}
		\hline \\[-0.5em]
		{ 
			Select either $\omega=1$ or $\omega = 2$. The third-order MR-NPRK method ansatz is:
		\begin{align}
			\begin{aligned}
				Y_1 &= y_n, \\
				{\color{magenta} Y_2} &= y_n + h\ahcmp{1}_{11} F({\color{magenta} Y_2},Y_1), \\
				{\color{ForestGreen} Y_3} &= y_n + h \left[ (\ahcmp{1}_{21} - a_{3,2,2}) F({\color{magenta} Y_2},Y_1) + a_{3,2,2} F({\color{magenta} Y_2},{\color{magenta} Y_2}) + \ahcmp{1}_{22} F({\color{ForestGreen} Y_3},Y_\omega)  \right], \\
				Y_i &= y_n + h \sum_{\substack{k=1 \\ k \ne 2,3}}^{i-1} \ahcmp{2}_{i-2,\max(k-2,1)} F({\color{magenta} Y_2},Y_k),\quad \text{for } i = 4,\ldots, \scmp{2}+1, \\
				Y_{\scmp{2}+2} &= y_n + h \sum_{\substack{k=1 \\ k \ne 2,3}}^{\scmp{2}} \ahcmp{2}_{i-2,\max(k-2,1)} F({\color{magenta} Y_2},Y_k) \\
				& \hspace{2.4em} + h \Big[ (\ahcmp{2}_{\scmp{2},\scmp{2}-1} - a_{\scmp{2}+2,3,\scmp{2}+1}) F({\color{magenta} Y_2},Y_{\scmp{2}+1}) \\
                & \qquad \hspace{2.4em}  + a_{\scmp{2}+2,3,\scmp{2}+1} F({\color{ForestGreen} Y_3},Y_{\scmp{2}+1}) \Big], \\
				{\color{Cerulean} Y_{s}} &= y_n +  
					h \left[ \bhcmp{2}_{1} - a_{s,3,1} \right]F({\color{magenta} Y_2}, Y_1) + h a_{s,3,1} F({\color{ForestGreen} Y_3}, Y_1) \\				
					& \hspace{2.4em}
					+ h\sum^{s-2}_{k=4} \bhcmp{2}_{k-2}  F({\color{magenta} Y_2},Y_k) 
					+ h \left[ \bhcmp{2}_{\scmp{2}} - a_{s,3,s-1} - \ahcmp{1}_{33} \right]F({\color{magenta} Y_2}, Y_{s-1}) \\
                    & \hspace{2.4em} + h a_{s,3,s-1} F({\color{ForestGreen} Y_3}, Y_{s-1}) + h\ahcmp{1}_{33} F({\color{Cerulean} Y_s},Y_{s-1}), \\	
				y_{n+1} &= {\color{Cerulean} Y_{s}},
			\end{aligned}
			\label{eq:mr-nprk-third-order-coupling-coeff-v1}
		\end{align}

		where $s = \scmp{2}+3$. 	The coupling coefficients are
			\begin{align*}
				&a_{3,2,2} = 
					\begin{cases}
				 		\frac{2-6 \lambda }{18 \lambda ^3-60 \lambda ^2+15 \lambda } \approx 0.1825367204468751798095174531383 & \omega = 1 \\
				 		\frac{2-6 \lambda }{18 \lambda ^3-60 \lambda ^2+15 \lambda }-\lambda \approx -0.2533298010615838196065019980553 & \omega = 2
				 	\end{cases} \\
				&a_{\scmp{2}+2,3,\scmp{2}+1} = \frac{1 - 3\chcmp{1}_1}{6	\bhcmp{2}_{\scmp{2}}(\chcmp{1}_2 - \chcmp{1}_1)} \\
				&a_{s,3,1} = \frac{\tfrac{1}{3} - \ahcmp{1}_{33} \chcmp{2}_{\scmp{2}} + \chcmp{1}_{1}(\ahcmp{1}_{33} \chcmp{2}_{\scmp{2}} - \tfrac{1}{2}) - \ahcmp{1}_{32}(\chcmp{1}_2 - \chcmp{1}_1)\chcmp{2}_{\scmp{2}}}{(\chcmp{1}_1-\chcmp{1}_2)\chcmp{2}_{\scmp{2}}} \\
				&a_{s,3,s-1} = \ahcmp{1}_{32} - a_{s,3,1}
			\end{align*}
			} \\[-1em] \hline
		\end{tabular}
		
		\caption{Third-order MR-NPRK \eqref{eq:mr-nprk-third-order-coupling-coeff-v1} with $\omega \in \{1,2\}$.}
		\label{tab:mr-nprk-third-order-v1}
		\end{table}
				
		\subsubsection*{Practical implementation for MR-NPRK method \eqref{eq:mr-nprk-third-order-coupling-coeff-v1}}		
		Given the initial condition $y_n$ and an explicit RK method $(\AhcmpB{2}, \bhcmpB{2})$ of at least third order:
		
		\begin{enumerate}
			\item Select either $\omega = 1$ or $\omega = 2$ -- this determines which stage is used as the explicit argument for computing stage ${\color{ForestGreen} Y_3}$.
			\item Solve the implicit equations for ${\color{magenta} Y_2}$ and ${\color{ForestGreen} Y_3}$ listed in line two and three of \cref{eq:mr-nprk-third-order-coupling-coeff-v1}.

			\item Compute one timestep of the ODE $w' = F({\color{magenta} Y_2}, w)$, $w_0=y_n$ using the explicit method $(\acmp{2}_{ij}, \bcmp{2}_i)$ with $\scmp{2}$ stages. In addition to the output $w_1$, also store the final two stage values $W_{\scmp{2}}$ and $W_{\scmp{2}-1}$.
			\item Define the function 
				\begin{align}
					H(w) := F({\color{ForestGreen} Y_3},w) - F({\color{magenta} Y_2},w),
					\label{eq:h-fun-mr-nprk-third-order-coupling-coeff-v1}
				\end{align} 
                then, $Y_{\scmp{2}+2}$ is given by
			    \begin{align}
				    Y_{\scmp{2}+2} = W_{\scmp{2}} + h a_{\scmp{2}+2,3,k} H(W_{\scmp{2}-1}).
			    \end{align}
			\item To obtain {\color{Cerulean} $Y_s$}, solve the implicit system ${\color{Cerulean} Y_{s}} = \xi + h\ahcmp{1}_{33} F({\color{Cerulean} Y_s},Y_{s-1})$ with
			\begin{align}
					\begin{aligned}
						\xi \coloneqq w_1 & + h\big[ b_{\scmp{2}} (F({\color{magenta} Y_2},Y_{s-1}) -F({\color{magenta} Y_2},W_{\scmp{2}})) \\
						&+ a_{s,3,1} H(y_n) + a_{s,3,s-1} H(Y_{s-1}) - \ahcmp{1}_{33}  F({\color{magenta} Y_2},Y_{s-1}) \big].
					\end{aligned}
				\label{eq:b-vector-mr-nprk-third-order-coupling-coeff-v1}
			\end{align}
		\end{enumerate}
	
	\noindent
	\begin{center}
	\begin{tabular}{|p{0.95\textwidth}|}
	  \hline
	  \vspace{-0.25em}
	  {\bf Pseudocode for MR-NPRK \eqref{eq:mr-nprk-third-order-coupling-coeff-v1}} \hfill $\omega \in \{1,2\}$ and $Y_1=y_n$  \\[.25em] \hline         
	        
	        \vspace{-0.25em}
	        Solve ${\color{magenta} Y_2} = y_n + h\ahcmp{1}_{11} F({\color{magenta} Y_2},y_n)$ \\[0.5em]
	        Solve ${\color{ForestGreen} Y_3} = y_n + h \left[ (\ahcmp{1}_{21} - a_{3,2,2}) F({\color{magenta} Y_2},y_n) + a_{3,2,2} F({\color{magenta} Y_2},{\color{magenta} Y_2}) + \ahcmp{1}_{22} F({\color{ForestGreen} Y_3},Y_\omega)  \right]$   \\[0.5em]

	        [$w_1$, $W_{\scmp{2}}$, $W_{\scmp{2}-1}$] = RK\_step($w' = F({\color{magenta} Y_2},w)$, $w_0 = y_n$, $\AhcmpB{2}$, $\bhcmpB{2}$, $h$) \\[.5em]
	        
	        $Y_{s-1} = W_{\scmp{2}} + ha_{\scmp{2}+2,3,\scmp{2}+1} H(W_{\scmp{2}-1})$ \hfill ($H(w)$ defined in \eqref{eq:h-fun-mr-nprk-third-order-coupling-coeff-v1}) \\[.5em]
	        
	        Solve ${\color{Cerulean} Y_s} = \xi + h F({\color{Cerulean}Y_s}, Y_{s-1})$ \hfill ($\xi$ defined in \eqref{eq:b-vector-mr-nprk-third-order-coupling-coeff-v1}) \\[.5em] \hline 
	\end{tabular}
	\end{center}

\subsubsection{Creating high-stage explicit integrators via composition}
\label{sec:high-stage-composition}
	To create an implicitly wrapped MR-NPRK with arbitrary $\scmp{2}$ (i.e., the number of stages used to integrate the explicit component), we require an explicit method with an arbitrary number of stages. A simple way to accomplish this is to compose a known $s$-stage method $m$ times to create an explicit method with $s \times m$ stages; specifically, given an $s$ stage explicit method $(\mathbf{A},\mathbf{b},\mathbf{c})$, we consider the composite method
	\begin{small}
	\begin{align}
		\renewcommand*{\arraystretch}{1.25}
		\begin{tabular}{c|ccccc}
				$\tfrac{\mathbf{c}}{m}$ & $\tfrac{\mathbf{A}}{m}$ \\
				$\tfrac{1}{m}+\tfrac{\mathbf{c}}{m}$ & $\tfrac{\mathbf{B}}{m}$ & $\tfrac{\mathbf{A}}{m}$ \\
				$\vdots$ & $\vdots$ & $\ddots$& $\ddots$ \\
				$\tfrac{m-1}{m} + \tfrac{\mathbf{c}}{m}$ & $\tfrac{\mathbf{B}}{m}$ & $\hdots$ & $\tfrac{\mathbf{B}}{m}$ & $\tfrac{\mathbf{A}}{m}$ \\ \hline
				& $\tfrac{\mathbf{b}}{m}$ & $\hdots$ & $\hdots$ & $\tfrac{\mathbf{b}}{m}$
		\end{tabular} && 
		\begin{aligned}
			\mathbf{B} &:= \mathbf{b}^T\mathbf{e}, \\
			\mathbf{e} &:= [1, \ldots, 1] \in \mathbb{R}^{1 \times s}.
		\end{aligned}
		\label{eq:rk-subcycled}
	\end{align}
	\end{small}%
	One significant benefit of this approach is that the storage cost for the composite integrator is identical to that of the base integrator $(\mathbf{A},\mathbf{b},\mathbf{c})$.

\section{Numerical experiments}\label{sec:numerical-experiments}
	
	An application for our proposed IMEX MR-NPRK methods are PDEs that couple hyperbolic and parabolic equations, such as those found in radiation hydrodynamics \cite{rad-hydro}. In this introductory work, we validate the new methods using a one-dimensional Burgers' equation with nonlinear diffusion and spatially varying coefficients:
	\begin{align}
		u_t = a(x)(|u|^{1/2}u_x)_x + b(x)(u^2)_x && 
		\left\{\begin{aligned}
			a(x) &= \tfrac{1}{2} + 2 \exp \left( - \tfrac{1}{25}(x - 1)^2 \right),  \\
		b(x) &= \exp \left( - \tfrac{1}{25} (x + 1)^2 \right)	.
		\end{aligned}\right.
\label{eq:burgers-nonlinear-diffusion-equation}
	\end{align}
	The spatial domain is $x\in [-2,2]$ and the boundary conditions are periodic. We discretize the domain using $300$ spatial grid points and integrate the PDE to time $t=5$. The coefficients $a(x)$ and $b(x)$ are selected so that: (1) nonlinear advection is present throughout the spatial domain, albeit with higher amount in the left-half of the domain and (2) nonlinear diffusion is only dominant in the right-half of the domain. The nonlinear diffusion term is discretized using 
    \begin{align}
		(k(u)u_x)_x =& \frac{q_{i+1/2} - q_{i-1/2}}{\Delta x}, \nonumber \\
			q_{i + 1/2} = \frac{k_{i+1/2} (u_{i+1} - u_{i+1})}{\Delta x}, 
			\quad & \quad q_{i - 1/2} = \frac{k_{i-1/2} (u_{i} - u_{i-1})}{\Delta x}, \label{eq:diffusion-discretization} \\
			k_{i+1/2} = \frac{k(u_i) + k(u_{i+1})}{2}, \quad  & \quad
			k_{i-1/2} = \frac{k(u_i) + k(u_{i-1})}{2} \nonumber
	\end{align}
	while the Burgers' nonlinearity $b(x)(u^2)_x$ is discretized using the fifth-order finite difference WENO method described in \cite{jiang1996efficient}. Code for reproducing the numerical experiments can be found in \cite{githubMRNPRK2025}.
	
	\subsection{Advantages of nonlinear partitioning}
	
	Due to the relationship between NPRK, ARK, and PRK methods (see \Cref{sec:nprk}), MR-NPRK methods can be applied to both additive and component partitioned equations. However, for this experiment, a nonlinear partition has advantages compared to additive or component partitions.
		\begin{itemize}[leftmargin=*]
			\item {\em Component partitioning of \eqref{eq:burgers-nonlinear-diffusion-equation}.} Since diffusion is only dominant in the right portion of the domain, a natural component partitioning is to separate the left and right portions of the spatial domain and then treat the left-half explicitly and the right-half implicitly. However, this requires us to treat both advection and diffusion implicitly in the right-half, which is non-trivial for a WENO discretization.
			\item {\em Additive partitioning of \eqref{eq:burgers-nonlinear-diffusion-equation}.} The natural additive partitioning separates the nonlinear diffusion and advection terms, and treats only the nonlinear diffusion implicitly. The primary disadvantage of this approach is that the nonlinear diffusion term requires a nonlinear solve at each timestep.

			\item {\em Nonlinear partitioning of \eqref{eq:burgers-nonlinear-diffusion-equation}}.	Nonlinear partitioning allows us to treat advection explicitly, and nonlinear diffusion semi-implicitly, such that only a linear solve is required at each timestep. Specifically, we consider the splitting
				\begin{align}
					F(u,v) = a(x)(k(v)u_x)_x + b(x) \left(v^2\right)_x
				\end{align}
				that after discretization this takes the form
				\begin{align}
					F(\mathbf{u},\mathbf{v}) =  \cdot \text{diag}(a(\mathbf{x})) * \mathbf{D}[\mathbf{v}]\mathbf{u} +  \text{diag}(b(\mathbf{x})) *  \mathbf{W}[\mathbf{v}^2]\mathbf{v}^2
				\end{align}
				where $\mathbf{D}[\mathbf{v}]$ is the discrete diffusion operator pertaining to \eqref{eq:diffusion-discretization} where $k_{i+1/2}$ and $k_{i-1/2}$ are computed using $\mathbf{v}$, and $\mathbf{W}[\mathbf{v}^2]$ is discrete advection operator arising from the WENO-5 method.

		\end{itemize}

	\subsection{MR-NPRK method selection}
	
	It is natural to treat advection WENO discretizations with strong stability preserving (SSP) methods that minimize unwanted oscillations \cite{jiang1996efficient,Gottlieb.2001}. We therefore consider implicitly wrapped MR-NPRK based on the well-known second and third-order SSP methods
	\renewcommand*{\arraystretch}{1.25}
	\begin{align}
		\text{SSP2:} \quad 
		\begin{tabular}{c|cc}
			0 & \\
			1 & $1$ \\ \hline
			  & $\tfrac{1}{2}$ & $\tfrac{1}{2}$
		\end{tabular} &&		
		\text{SSP3:} \quad 
		\begin{tabular}{c|ccc}
			0 & \\
			1 & $1$ \\ 
			$\tfrac{1}{2}$ & $\tfrac{1}{4}$ & $\tfrac{1}{4}$  \\ \hline
			  & $\tfrac{1}{6}$ & $\tfrac{1}{6}$ & $\tfrac{2}{3}$
		\end{tabular}
		\label{eq:ssp-rk-methods}
	\end{align}
	To create multirate methods with large $\scmp{2}$, we simply compose these methods $m$ times to create an explicit method with $2m$ or $3m$ stages (see \Cref{sec:high-stage-composition}).

	Our numerical experiment compares explicit RK, IMEX-NPRK, and IMEX MR-NPRK of orders two and three. The methods we consider are:
	
	\vspace{1em}
	\begin{center}
		\begin{tabular}{|rp{.61\textwidth}|} \hline
			{\bf MR-NPRK2-[ssp2-$m$x]} & Second-order MR-NPRK \eqref{eq:mr-nprk-second-order-minus} with $a_{2,2,1} = \frac{2 - \sqrt{2}}{2}$, $a_{s,s,s-1} = \frac{2 - \sqrt{2}}{2}$. The underlying explicit integrator is the SSP2 method from \eqref{eq:ssp-rk-methods} composed $m$ times. 	\\
			{\bf MR-NPRK3-1[ssp3-$m$x]} & Third-order MR-NPRK \eqref{eq:mr-nprk-third-order-coupling-coeff-v1} with $\omega = 2$. The underlying explicit integrator is the SSP3 method from \eqref{eq:ssp-rk-methods} composed $m$ times. \\
			{\bf ssp2-[$m$x]} & the SSP2 method from \eqref{eq:ssp-rk-methods} composed $m$ times.	\\	
		   {\bf ssp3-[$m$x]} &	 the SSP3 method from \eqref{eq:ssp-rk-methods} composed $m$ times.	\\
		  {\bf NPRK2[32]b} & second-order, stiffly-accurate NPRK from \cite{nprk1}. \\
		 {\bf NPRK3[54]-Sa} & third-order, stiffly-accurate  NPRK from \cite{nprk1}. \\ \hline
		\end{tabular}
	\end{center}
	
	\subsection{Evaluating convergence and efficiency}
	
	We validate the efficiency and convergence of our methods using two initial conditions (IC):
	\begin{align}
		\text{\bf Two Gaussian IC:}& & u_0^{[1]}(x) &= \tfrac{1}{100} + \exp\left(-60 \left(x + \tfrac{3}{2}\right)^2\right) +  \exp\left(-60 x^2\right)
		\label{eq:burgers-nonlinear-diffusion-ic2} \\
		\text{\bf Three Gaussian IC:}& &  u_0^{[2]}(x) &= u_0^{[1]}(x) + \exp\left(-60 \left(x - \tfrac{3}{2}\right)^2\right).
		\label{eq:burgers-nonlinear-diffusion-ic3}
	\end{align}
	The two Gaussian pulses start in the left-half of the domain and advect leftwards forming shocks that eventually travel through the diffusive region. This implies that advection is the dominant term until the shockwave enters the diffusive region. As we will see, the three Gaussian pulse initial condition is more challenging to multirate because the additional Gaussian is located in the diffusive region of the domain, thus advection and diffusion both drive the dynamics from the start.
		
\begin{figure}[h]
		
		\begin{center}
			
			\begin{minipage}{.48\textwidth}
				\centering
				{\em Initial Condition \eqref{eq:burgers-nonlinear-diffusion-ic2}}
				
				\includegraphics[width=\textwidth]{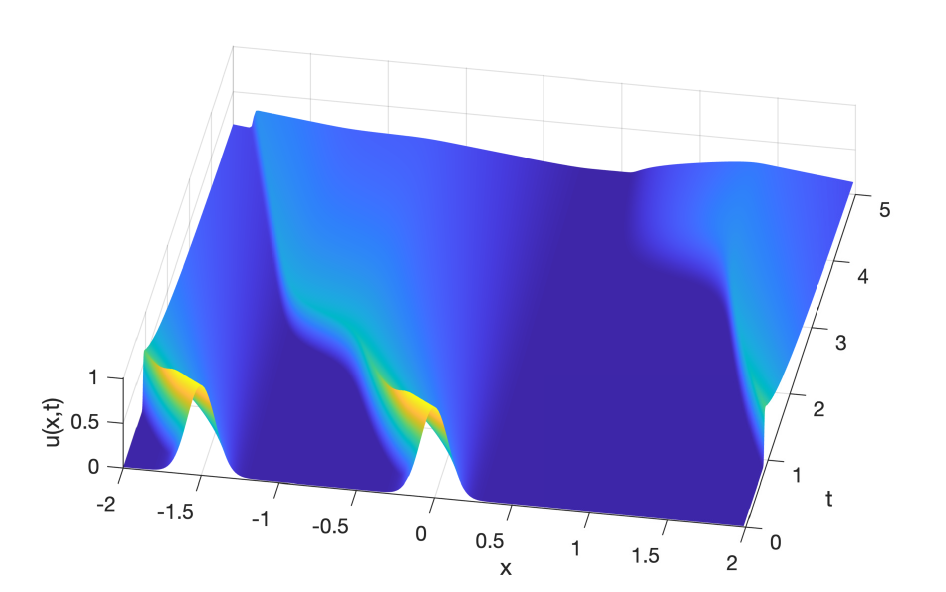}
				
			\end{minipage}
			\begin{minipage}{.48\textwidth}
				\centering
				{\em Initial Condition \eqref{eq:burgers-nonlinear-diffusion-ic3}}

				\includegraphics[width=\textwidth]{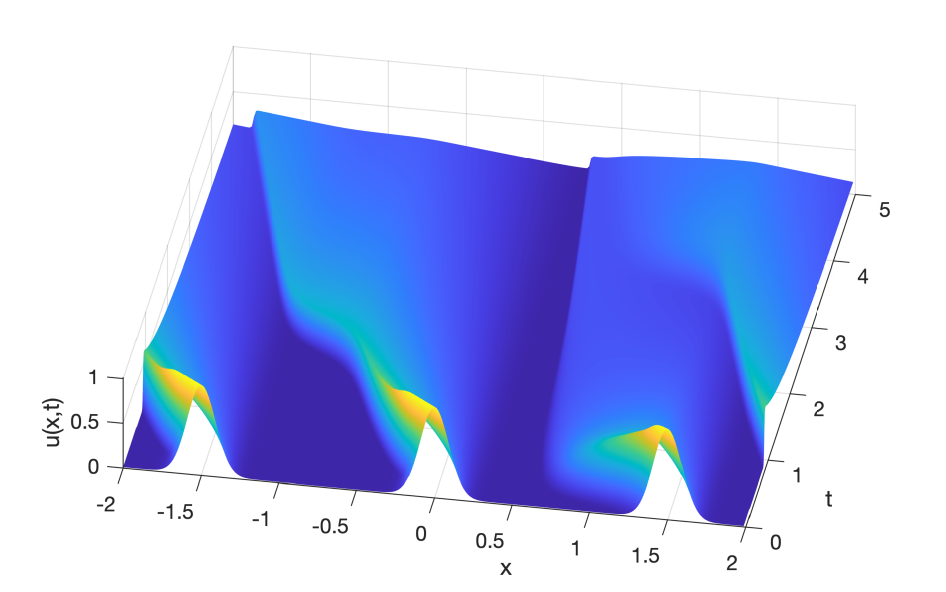}
				
			\end{minipage}	
			
		\end{center}
		
		\caption{Solution of the PDE \eqref{eq:burgers-nonlinear-diffusion-equation} with the two initial conditions \cref{eq:burgers-nonlinear-diffusion-ic2,eq:burgers-nonlinear-diffusion-ic3}.}	
		\label{fig:burgers-nld-solution-plot}
		
	\end{figure}

\begin{figure}
		
		\begin{center}			
			\begin{tabular}{|@{}c@{}|} 
				\hline
				\begin{tabular}{c|cc}
					& {\footnotesize {\bf Convergence }} & {\footnotesize {\bf Efficiency}}  \\[-.5em] 
					& {\footnotesize Second-Order Methods} & {\footnotesize Second-Order Methods} \\ 
					\rotatebox{90}{\hspace{4em} IC \cref{eq:burgers-nonlinear-diffusion-ic2}} & \includegraphics[width=0.42\textwidth,align=b]{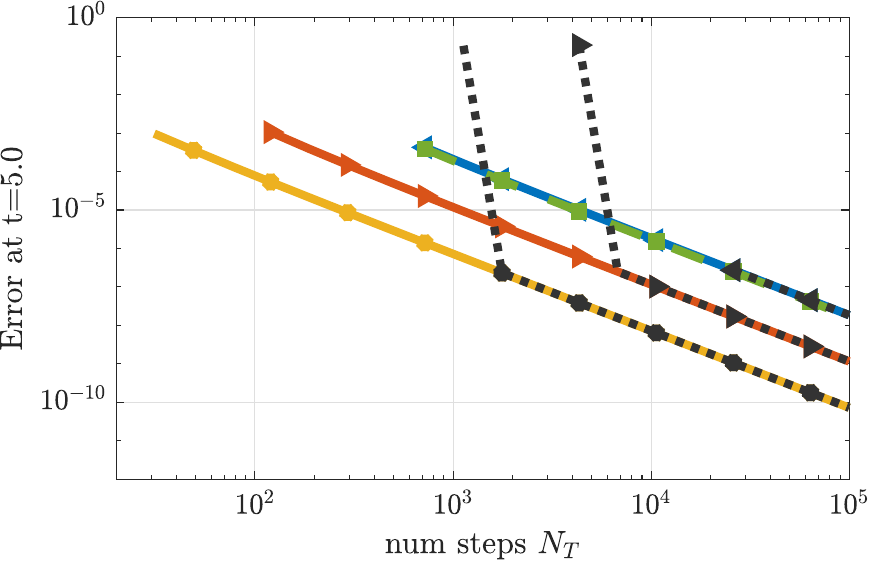} & \includegraphics[width=0.42\textwidth,align=b]{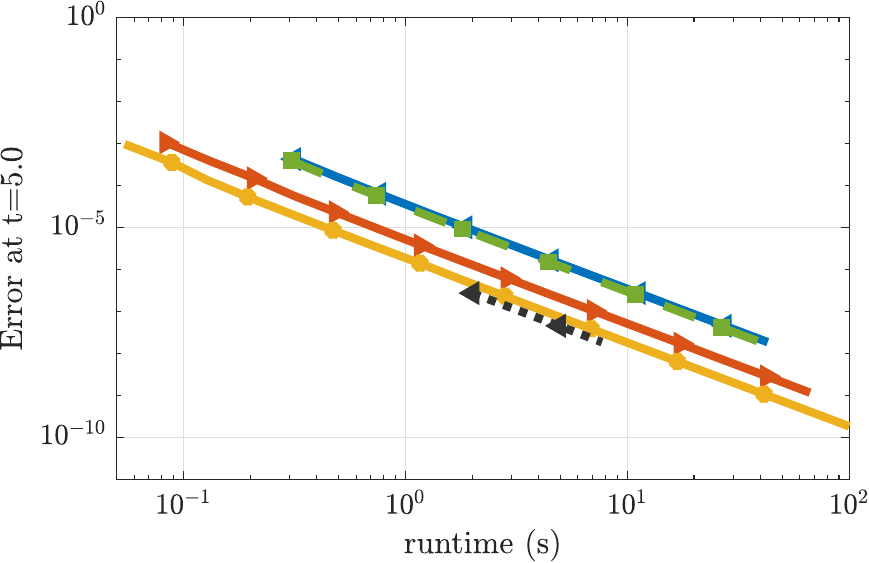}	 \\
					\rotatebox{90}{\hspace{4em} IC \cref{eq:burgers-nonlinear-diffusion-ic3}} &
					\includegraphics[width=0.42\textwidth,align=b]{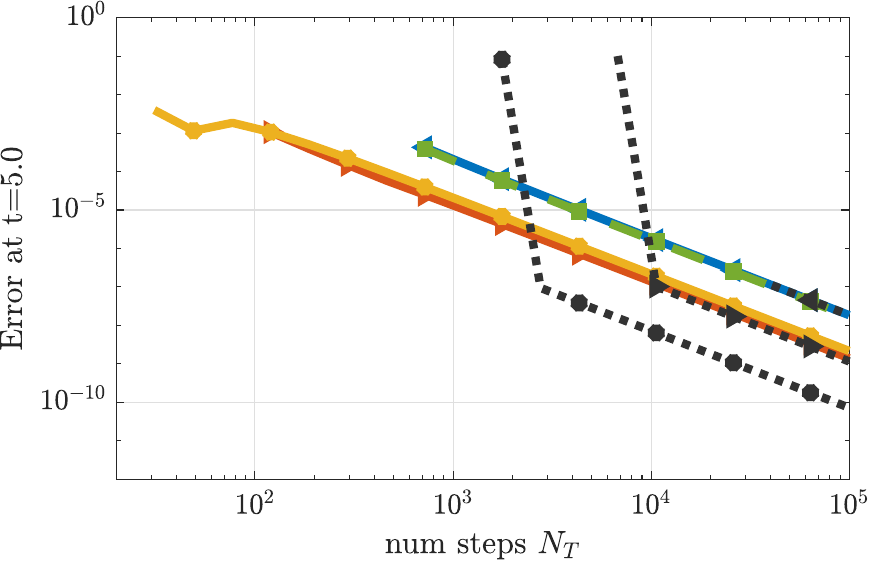} & 
					\includegraphics[width=0.42\textwidth,align=b]{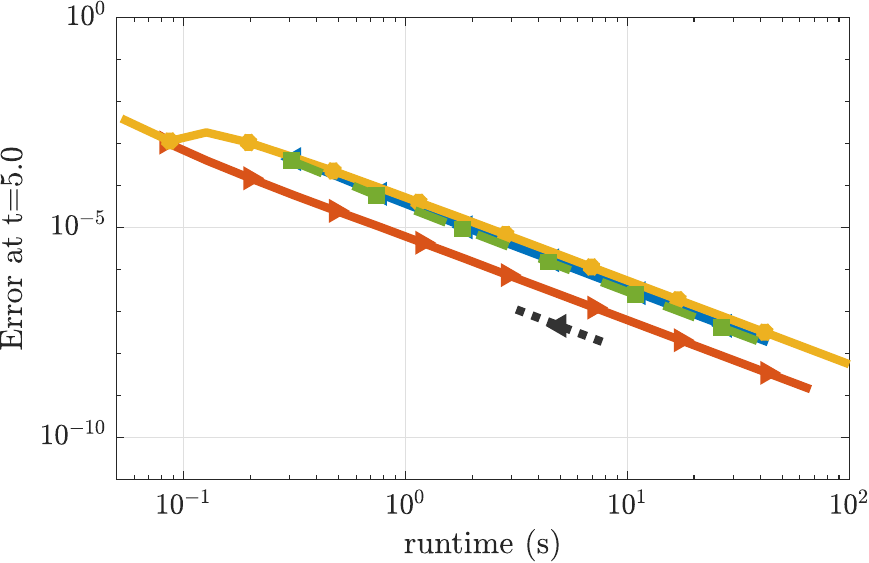} 
				\end{tabular} \\ \hline \\[-.75em]
				
				\includegraphics[width=.93\textwidth,trim={6.5cm 7.5cm .45cm .5cm},clip]{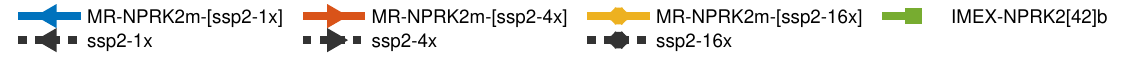}	\\[.25em] \hline	
							
			\end{tabular}			
		\end{center}
		
		\vspace{1em}
		\begin{center}			
			\begin{tabular}{|@{}c@{}|} 
				\hline
				\begin{tabular}{c|cc}
					& {\footnotesize {\bf Convergence }} & {\footnotesize {\bf Efficiency}}  \\[-.5em] 
					& {\footnotesize Third-Order Methods} & {\footnotesize Third-Order Methods} \\ 
					\rotatebox{90}{\hspace{4em} IC \cref{eq:burgers-nonlinear-diffusion-ic2}} & \includegraphics[width=0.42\textwidth,align=b]{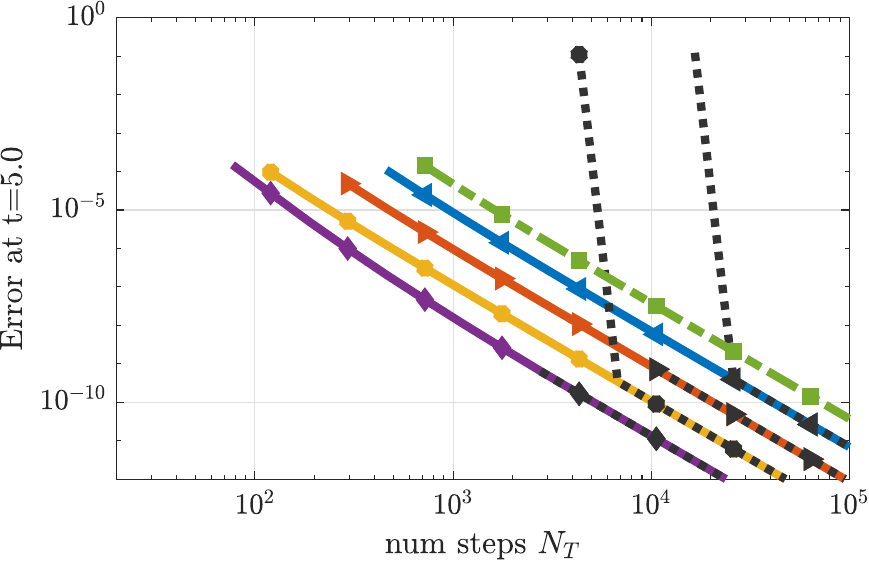} & \includegraphics[width=0.42\textwidth,align=b]{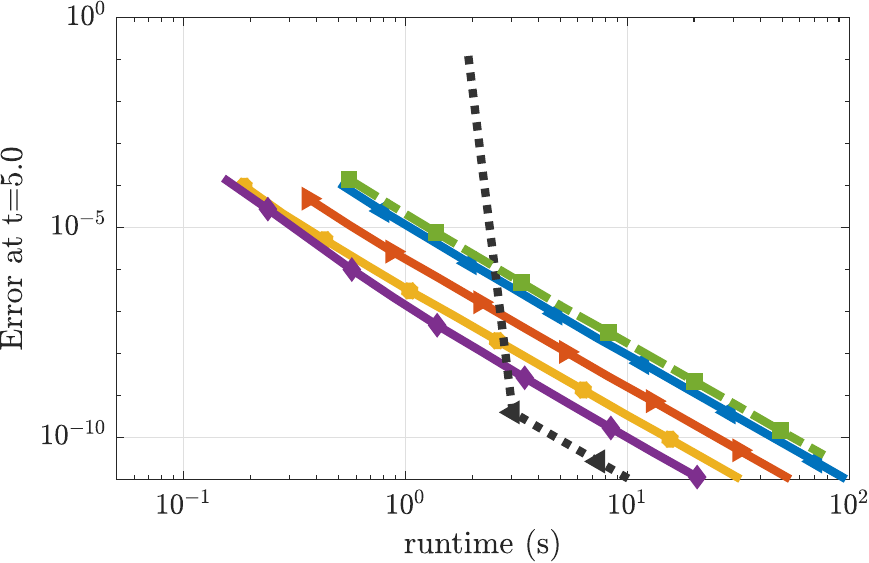}	 \\
					\rotatebox{90}{\hspace{4em} IC \cref{eq:burgers-nonlinear-diffusion-ic3}} &
					\includegraphics[width=0.42\textwidth,align=b]{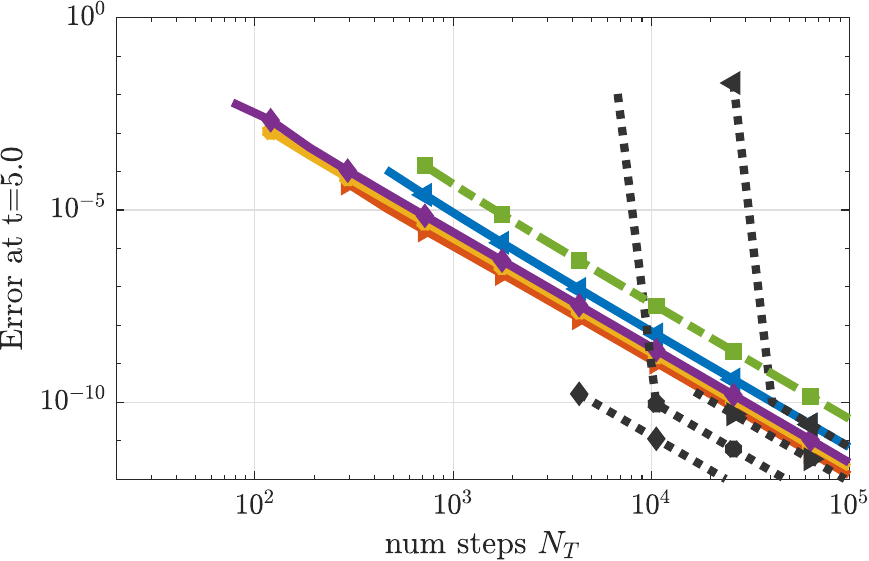} & 
					\includegraphics[width=0.42\textwidth,align=b]{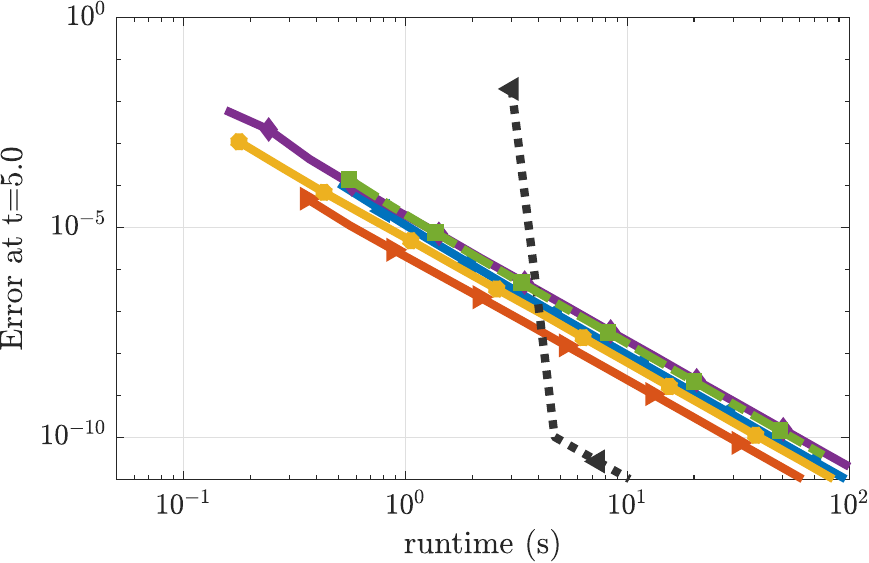} 
				\end{tabular} \\ \hline \\[-.75em]
				
				\includegraphics[width=.93\textwidth,trim={6.5cm 7.3cm .6cm .6cm},clip]{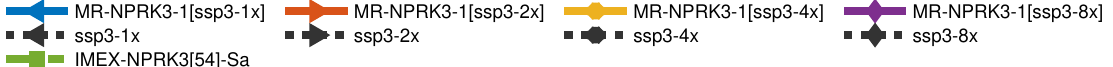}	\\[.25em] \hline							
			\end{tabular}			
		\end{center}
		
		\caption{Numerical solution of \cref{eq:burgers-nonlinear-diffusion-equation} with initial conditions (IC) \cref{eq:burgers-nonlinear-diffusion-ic2,eq:burgers-nonlinear-diffusion-ic3}. Curves for second-order methods are in the top box while those for third-order methods are in the bottom box. We only show the efficiency curves for ssp2-[1x] and ssp3-[1x] since all efficiency curves for ssp2-[$m$x] and ssp3-[$m$x], respectively, overlap.}	
		\label{fig:results}
		
	\end{figure}

	Visualizations of the solution for the initial conditions \eqref{eq:burgers-nonlinear-diffusion-ic2} and \eqref{eq:burgers-nonlinear-diffusion-ic3} are shown in \Cref{fig:burgers-nld-solution-plot}, and convergence and efficiency plots for second and third order methods are shown in \Cref{fig:results}.
	
	
	\subsection{Discussion} 
	
	Our numerical experiment reveals that the newly proposed MR-NPRK methods: (1) converge at their expected order of accuracy and (2) are able to achieve improved computational efficiency compared to both single-rate NPRK methods and explicit SSP methods. We discuss the efficiency gains in more detail:
	\begin{itemize}
		
		\item {\em Efficiency gains over explicit RK.} The cost per timestep of an implicitly-wrapped MR-NPRK method will always exceed that of its underlying explicit integrator due to the additional implicit solves. However, stiffness introduced by the nonlinear diffusion causes instability for all explicit methods at coarse timesteps. We remark that this phenomenon will become more pronounced as the spatial mesh is refined or the diffusion strength is increased.		 
		
		\item {\em Efficiency gains over single-rate IMEX-NPRK.} For a fixed timestep, MR-NPRK methods are consistently able to achieve improved error compared to single-rate NPRK. The gain in accuracy arises from treating the second component more accurately. The gain in efficiency is possible since the cost of explicit function evaluations is less than an implicit solve (in our experiment, a solve is four times more expensive than an evaluation -- this ratio can be larger on higher dimensional and more complex problems). Additionally, MR-NPRK methods are able to maintain stability at significantly larger timesteps than NPRK methods.
	
	\end{itemize}
	
	We also discuss the degree to which multirating improves accuracy. For the two Gaussian initial condition, all implicitly wrapped MR-NPRK methods are able to achieve almost exactly the same accuracy as the underlying explicit integrator. However, when we introduced a third Gaussian, this effect was significantly reduced. These initial conditions were intentionally selected to vary the coupling strength between diffusion and advection. 
	\begin{itemize}
		\item {\em Two Gaussian IC:} on the time interval $t\in[0,4]$ the two shocks travel through a nearly diffusion-free region. Therefore, for short times, the accuracy is primarily dependent on the explicit integrator. To obtain an accurate final answer the final time $t=6$, it is necessary to first accurately integrate the shock before it enters the diffusive region. This makes multirating the advective term highly beneficial. 
		\item {\em Three Gaussian IC:} The additional Gaussian was intentionally placed inside the diffusive region. In contrast to the previous initial condition, it is now equally necessary to also accurately resolve diffusion. However, since diffusion is stepped using a coarse (macro) timestep, the error does not improve by significantly by multirating.
	\end{itemize}
	Despite the fact that our MR-NPRK methods do not achieve significantly improved error on the three Gaussian initial condition, multirating improves stability in both cases. Specifically, if one requires only moderate precision, then MR-NPRK3-1[ssp-8x] and MR-NPRK2m[ssp2-16x] methods are able to obtain a stable solution about four times faster than the single-rate NPRK methods.

\section{Conclusion}\label{sec:conclusion}

We developed a framework for multirate NPRK methods that can be applied to solve evolution equations with nonlinear splittings. We studied order conditions, linear stability, and introduced new IMEX nonlinear-partitioned multirate methods of order two and three. Our new methods are flexible and allow the user to specify any explicit RK method to solve the fast timescale. Moreover, the methods have a simple structure that allows them to be easily incorporated into existing codes. More broadly, this work generalizes the existing multirate literature that is based on additive and component partitioned systems, thus creating the potential to multirate complex multiscale physical systems whose scales are not linearly separated.

In future work we will evaluate MR-NPRK schemes on hyperbolic equations with stiff relaxation and coupling, including radiation hydrodynamic equations that model the interaction of photons with compressible flow in high energy density regimes.

\appendix

\section{Stability proofs for IMEX MR-NPRK methods}\label{appendix:stability-analysis}
This appendix provides further discussion of the linear stability analysis of MR-NPRK methods and contains proofs of the results stated in \Cref{sec:linear-stability}.

\subsection{The stability function - separating scales}\label{sec:constructing:stab:scales}

	We present a strategy for rewriting the stability function that makes it easier to study the stability of MR-NPRK methods. For multirate methods with an arbitrary number of substeps, the stability function \eqref{eq:ark-stability-region} is challenging to analyze directly since the matrices $\AcmpB{r}$ and $\mathbf{e}(\bcmpB{r})^T$ become high-dimensional. However, it is possible to rewrite the stability function using permutation matrices that separate the coefficients of the two underlying methods. This allows us to more easily obtain some general results regarding linear stability.

Before introducing the permutation matrices, we first describe the structure of the underlying method tableaux for fully-decoupled and $1\to2$ coupled IMEX MR-NPRK methods ($2\to1$ coupled is identical to $1\to2$ except with arguments flipped). Recall that we assume $\Scmp{1} \cup \Scmp{2} = \mathcal{S}$ (\Cref{assump:stage-union}); furthermore, by \Cref{cor:disjoint-not-fully-coupled}, the index sets have trivial intersection, $\Scmp{1} \cap \Scmp{2} = \emptyset$. This allows us to separate all stages in terms of those that belong to $\Mcmp{1}$ and those that belong to $\Mcmp{2}$.

\begin{itemize}
	\item {\em Underlying implicit method}. $\Mcmp{1}$ is diagonally-implicit, which implies that $\AcmpB{1}$ is lower triangular. The method has $\cmp{s}{1}$ irreducible stages whose coefficients are contained in rows of $\AcmpB{1}$ with indices $i \in \Scmp{1}$. The only remaining potentially non-zero entries are located in the columns with index $j\in \Scmp{1}$. The structure of $\AcmpB{1}$ is identical for both fully-decoupled or $1\to 2$ coupled methods. 

	\item {\em Underlying explicit method}. $\Mcmp{2}$ is explicit ,which implies a strictly lower triangular $\AcmpB{2}$. The underlying explicit method $\Mcmp{2}$ has $\cmp{s}{2}$ irreducible stages whose coefficients are located in rows of  $\AcmpB{2}$ with index $i \in \Scmp{2}$. All remaining rows with index $i \in \Scmp{1}$ are all reducible. If the method is fully-decoupled then the columns with index $j \in \Scmp{1}$ are also all zero.
\end{itemize}

\begin{figure}[h!]
    \centering
	\begin{subfigure}[t]{0.95\textwidth}
		\label{subfig:underlying-method-matrix-structure-a}
		\caption{Underlying method matrices $\AcmpB{r}$ and vectors $\bcmpB{r}$ for an example method with $\scmp{1}=3$.}
		\vspace{1em}
		\includegraphics[width=\textwidth]{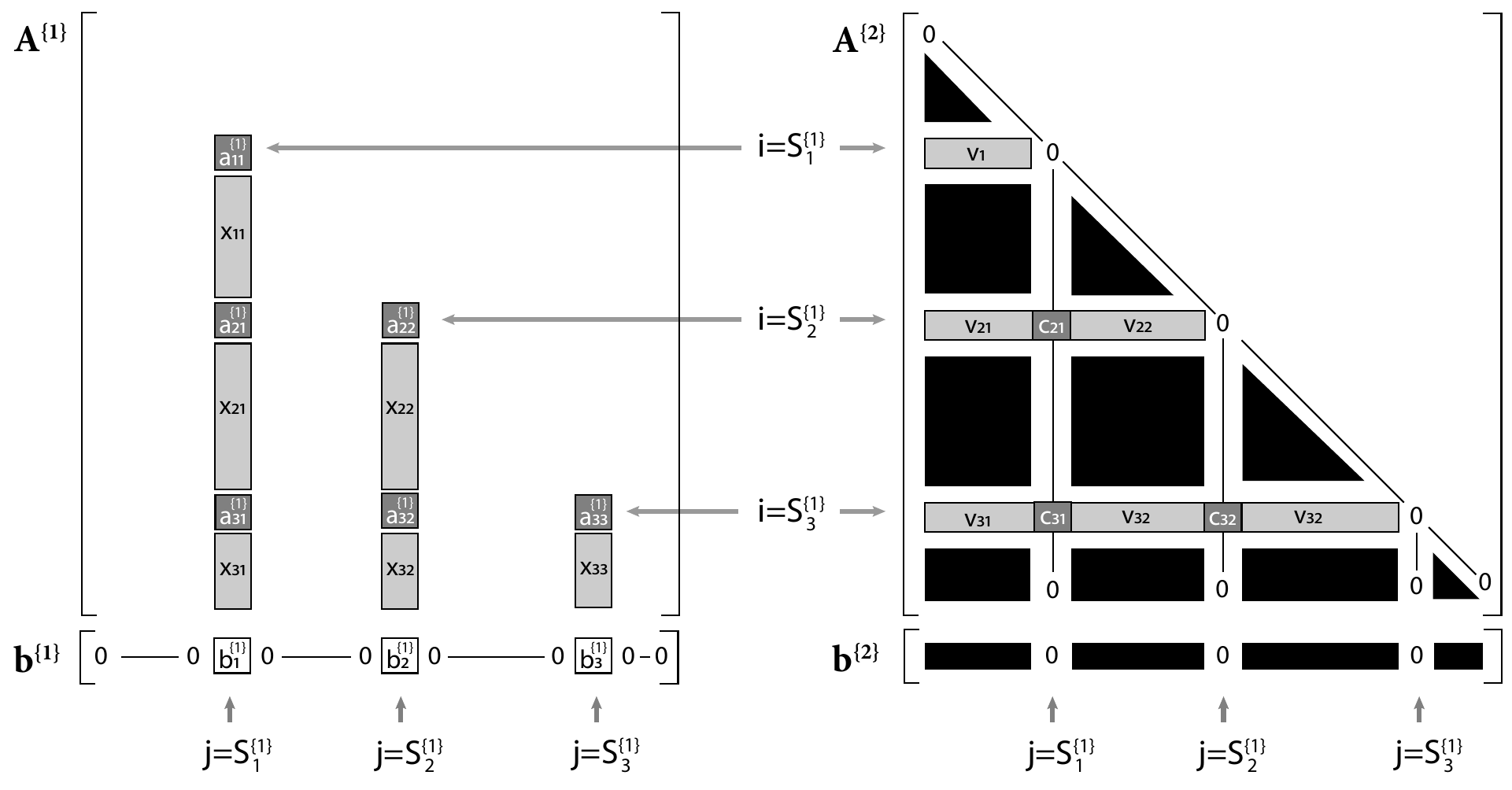}	
	\end{subfigure}

	\vspace{1em}
		
	\begin{subfigure}[t]{0.95\textwidth}
		\label{subfig:underlying-method-matrix-structure-b}
		\caption{Permuted $\AcmpB{r}$ and $\bcmpB{r}$ from (a); the matrix $\mathbf{P}$ is defined in \eqref{eq:underlying-permutation-matrix}}
		\vspace{1em}
		\includegraphics[width=\textwidth]{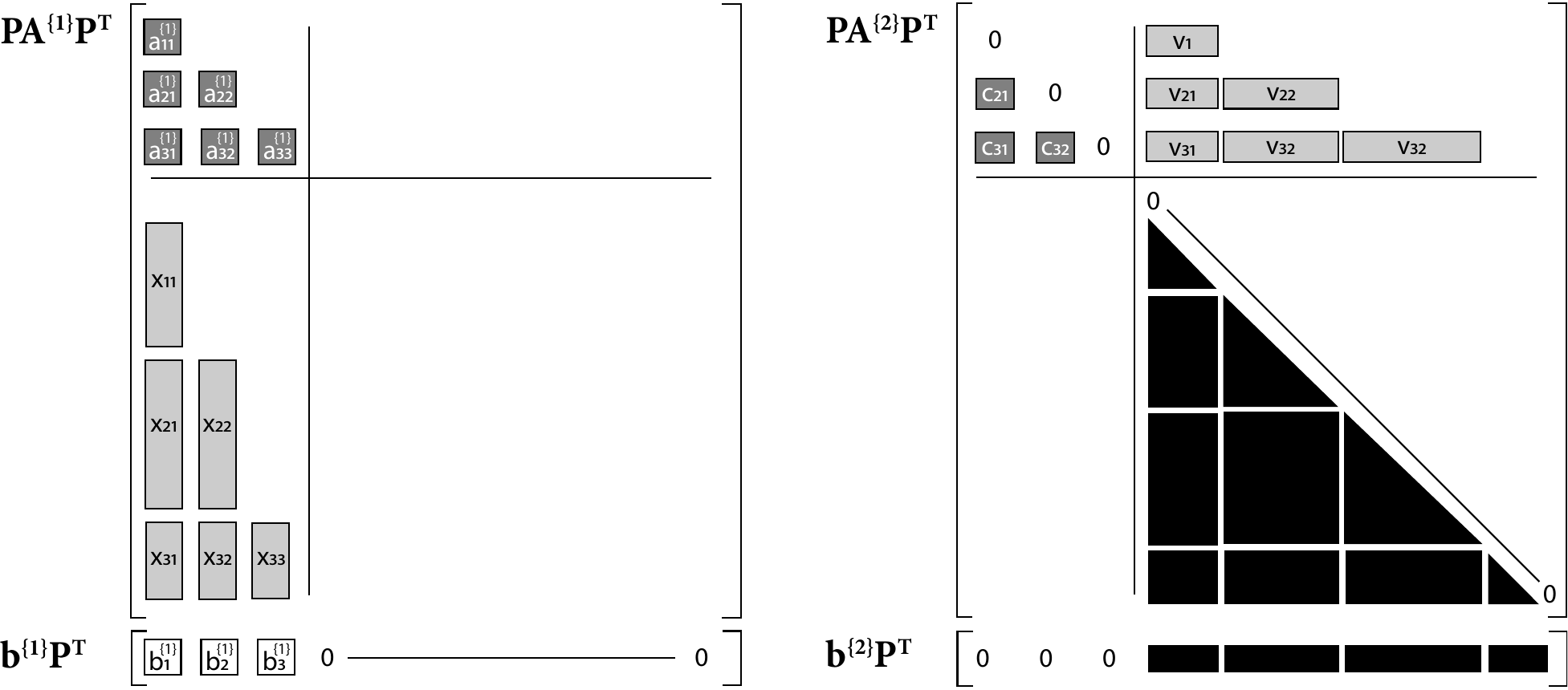}	
	\end{subfigure}
	
	\caption{\small{Butcher matrices $\AcmpB{r}$ pertaining to the underlying methods of a fully-decoupled or $1 \to 2$ coupled IMEX MR-NPRK method. A fully-decoupled method has constants $c_{jk}=0$ in $\AcmpB{2}$. All labeled indices belong to $\Scmp{1}$ while all non-labeled indices belong to $\Scmp{2}$ since $\Scmp{1} \cup \Scmp{2} = \mathcal{S}$ (see \Cref{prop:intersection-fully-coupled})}}.
	\label{fig:underlying-method-matrix-structure}
	\vspace{-3em}
\end{figure}

In \cref{fig:underlying-method-matrix-structure}(a), we illustrate the matrices for a generic method with $\scmp{1}=3$. By introducing the appropriate row and column permutations, one can always separate reducible and irreducible stages in both $\AcmpB{1}$ and $\AcmpB{2}$. Specifically, let $\Scmp{1}_k$ denote the $k^{th}$ element of $\Scmp{1}$, then let $\delta_j$ represent the number of irreducible stages of $\Mcmp{2}$ between sequential pairs of irreducible stages of $\Mcmp{1}$, defined by 
\begin{align}
	\delta_0 := \Scmp{1}_1 - 1, &&
	\delta_j := \Scmp{1}_{j+1} - \Scmp{1}_j - 1, \quad 1\le j < \scmp{1}, &&
	\delta_{\cmp{s}{1}} := s - \Scmp{1}_{\scmp{1}}.
\end{align}
Furthermore, let $\mathbf{e}_n$ be the $n^{th}$ unit vector in $\mathbb{R}^{\cmp{s}{1}}$, $\mathbf{I}(n) \in \mathbb{R}^{n \times n}$ be the identity, and introduce the permutation matrix $\mathbf{P} \in \mathbb{R}^{s \times s}$ defined by 
\begin{align}
	\mathbf{P} := \left[
	\renewcommand*{\arraystretch}{1.1}
	\begin{array}{c:c:c:c}
		& \mathbf{E}(1) & \hdots & \mathbf{E}(\cmp{s}{1}) \\ \hdashline
		\mathbf{I}(\delta_0) & 				&		 & \\ \hdashline
					& \mathbf{Z}(\delta_1) 	&		 & \\ \hdashline
					&				& \ddots & \\ \hdashline
					&				&		 & \mathbf{Z}(\delta_{\cmp{s}{1}})
	\end{array}
	\right], &&
	\begin{aligned}
		\mathbf{Z}(n) &:= 
		\left[
		\begin{array}{c|ccc}
	 			0		& 1 &  \\
	 			\vdots	& & \ddots &  \\
	 			0		& & & 1
	 	\end{array}
	 	\right] \in \mathbb{R}^{n,n+1}, \\
	 	\mathbf{E}(n) &:= 
		\left[
		\begin{array}{c|ccc}
	 							& 0 & \ldots & 0 \\
	 			\mathbf{e}_n		& \vdots & & \vdots \\
	 							& 0 & \ldots & 0
	 	\end{array}
	 	\right] \in \mathbb{R}^{\cmp{s}{1},n+1}.
	\end{aligned}
	\label{eq:underlying-permutation-matrix}
\end{align}
The matrices $\mathbf{P}\AcmpB{r}\mathbf{P}^T$ separate reducible and irreducible stages. Specifically, rows and columns $1, \ldots, \cmp{s}{1}$ of $\mathbf{P} \AcmpB{1} \mathbf{P}^T$ contain the reduced implicit method, while rows and columns $\cmp{s}{1}+1, \ldots, s$ of $\mathbf{P} \AcmpB{2} \mathbf{P}^T$ contain the reduced explicit method. The remaining entries in both matrices correspond to reducible coupling terms; \Cref{fig:underlying-method-matrix-structure}(b) contains  an illustration for a generic method with $\scmp{1}=3$.

\subsection{Proof of L-stability  \Cref{thm:z1-L-stable-ansatz}}\label{sec:constructing:stab:deg}

Here, we prove the main stability result presented in the paper, \Cref{thm:z1-L-stable-ansatz}. To achieve L-stability in the stiff $z_1$ limit, we require the numerator of the stability function to have a lower degree in $z_1$ than the denominator. Additionally, we prefer that the numerator is a polynomial of low degree in $z_2$ to avoid growth for large $|z_2|$ when $|z_1|$ is finite. The permutation matrix \eqref{eq:underlying-permutation-matrix} allows us to make some general claims about the degree of the polynomials in the numerator and denominator of the stability function \eqref{eq:ark-stability-function}. Specifically, we show how a fully-decoupled structure can help achieve L-stability, and the challenges that arise when one introduces coupling. We begin by providing several technical results, that will yield an immediate proof to \Cref{thm:z1-L-stable-ansatz}.

We start with two simple remarks about the determinant of a matrix with coefficients that depend linearly on one or two variables.

\begin{remark}
	Let $\mathbf{M}(z) \in \mathbb{R}^{n \times n}$, where every coefficient is at most a linear function of $z$. If then follows that $\det (\mathbf{M}(z))	$ is a polynomial in $z$ of degree at most $n$ (this follows from induction and expansion by cofactors).
	\label{rmk:det-z-powers}
\end{remark}

\begin{remark}
	Let $\mathbf{M}_1, \mathbf{M}_2 \in \mathbb{R}^{n \times n}$, then $\det (\mathbf{I} - z_1\mathbf{M}_1 - z_2 \mathbf{M}_2)	$ is a polynomial in $z_1$ and $z_2$ of degree at most $n$ in both variables (this follows from induction and recursive definition of determinant expanded by minors). 
	\label{rmk:det-z-powers-nprk}
\end{remark}

Next, we discuss the numerator and denominator of the stability function.  

\begin{proposition}\label{prop:stability-denominator}
	The stability function denominator $R(z_1,z_2)$ for a diagonally-implicit IMEX MR-NPRK method is a polynomial of degree $\scmp{1}$ in $z_1$ given by $R(z_1,z_2) = \prod_{i=1}^{s} \left(I - z_1 \acmp{1}_{ii}\right).$ Moreover, if the method's implicit solves only occur in stages with index $\nu \in S^1$ (i.e. stages pertaining to $\AhcmpB{1}$) then, $R(z_1,z_2) = \prod_{i=1}^{\scmp{1}} \left(I - z_1 \cmp{\widehat{a}}{1}_{ii}\right).$ 

\begin{proof}
	For a diagonally implicit IMEX NPRK integrator, $\AcmpB{1}$ is lower triangular and $\AcmpB{2}$ is strictly lower triangular. If the only implicit stages have indices belonging to $\Scmp{1}$, then $a^{\{1\}}_{ii} \ne 0 \implies  \in \Scmp{1}$. Therefore, both statements follow immediately from \eqref{eq:ark-stability-region} since the determinant is simply the product of the diagonals. \qed
\end{proof}
\end{proposition}

\begin{proposition}\label{prop:stability-numerator}
	\label{thm:mr-num-dem-degree}
	The numerator of the stability function of an $s$-stage MR-NPRK method is a bivariate polynomial in $z_1$, $z_2$ whose degree in $z_1$ and $z_2$ is bounded above as follows:
	\begin{center}
		\begin{tabular}{c|c|c|c|c}
			& \text{fully-decoupled:} & \text{$1 \to 2$} & \text{$2 \to 1$} & \text{fully-coupled:}	\\ \hline
			max degree in $z_1$ & $\scmp{1}$ & $\scmp{1}$ & $s$ & $s$\\
			max degree in $z_2$ & $\scmp{2}$ & $s$ & $\scmp{2}$ & $s$
		\end{tabular}
	\end{center}
\begin{proof}
After permuting rows and columns with \eqref{eq:underlying-permutation-matrix}, the matrix in the numerator of the stability function has the form
\begin{align}
	\left[
	\begin{array}{c|c}
		I + z_1 \mathbf{M}_{11} + z_2 \mathbf{C}_{11} & z_2 \mathbf{M}_{12} \\ \hline
		z_1 \mathbf{M}_{21} & I + z_2 \mathbf{M}_{22} + z_1 \mathbf{C}_{22}
	\end{array}
	\right],
\end{align}
where $\mathbf{M}_{ij}, \mathbf{C}_{ij} \in \mathbb{R}^{\scmp{i} \times \scmp{j}}$; $\mathbf{C}_{11}$ and $\mathbf{C}_{22}$ may be zero depending on the coupling:
\begin{align*}	
	\begin{tabular}{c|c|c}
		\text{$1 \to 2$} & \text{$2 \to 1$} & \text{fully-decoupled: }	\\ \hline
		$\mathbf{C}_{22} = 0$ & $\mathbf{C}_{11} = 0$ & $\mathbf{C}_{22} = 0 \text{ and } \mathbf{C}_{11}=0$
	\end{tabular}
\end{align*}
The degree in $z_1$ and $z_2$ can be immediately bounded above by $s$ using \Cref{rmk:det-z-powers-nprk} applied to the full matrix. Bounds less than $s$ follow from combining the following block determinant formulas with  \Cref{rmk:det-z-powers-nprk}:
\begin{align}
		\text{bounding degree in $z_1$:} && \det \left( \begin{bmatrix} A & B \\ C & D \end{bmatrix} \right) = 
		\det(D) \det \left(A - BD^{-1}C \right), 
		\label{eq:block-det-formula-z1}\\
		\text{bounding degree in $z_2$:} && \det \left( \begin{bmatrix} A & B \\ C & D \end{bmatrix} \right) = 
		\det(A) \det \left(D - CA^{-1}B \right)	.
		\label{eq:block-det-formula-z2}
	\end{align} \qed
\end{proof}
\end{proposition}

In summary, by selecting an IMEX MR-NPRK that is fully-decoupled or $1\to 2$ couplied, we can guarantee that $\lim_{|z_1|\to \infty} |R(z_1,z_2)| < \infty$; however, we still require additional constraints to guarantee L-stability in the stiff $z_1$ limit \eqref{eq:l-stability-z1-limit}. With one additional technical result, we can then prove L-stability.
\begin{lemma}
	An $s$-stage, diagonally-implicit, stiffly-accurate IMEX MR-NPRK method satisfies $s \in \Scmp{1}$.
	\label{lem:stiff-accurate-final-stage-s1}	
\end{lemma}
\begin{proof}
	Stiff accuracy implies that $y_{n+1} = Y_s$ for both $\Mcmp{1}$ and $\Mcmp{2}$. IMEX implies that all stages of $\Mcmp{2}$ are explicit. Now	suppose that $s \in \Scmp{2}$, then by definition the stage $Y_s$ cannot be reducible from $\Mcmp{2}$.  However, this leads to a contradiction; since $\cmp{y}{2}_{n+1}=\cmp{Y}{2}_s$ and $\cmp{Y}{2}_s$ is explicit, then $\cmp{Y}{2}_s$ is always reducible since it is unused by the output:
	\begin{align*}
		\cmp{Y}{2}_s &= y_n + h \sum_{j=1}^{s-1} \acmp{2}_{sj} G(\cmp{Y}{2}_j), &
		\cmp{y}{2}_{n+1} &= y_n + h \sum_{j=1}^{s-1} \acmp{2}_{sj} G(\cmp{Y}{2}_j).
	\end{align*} \qed
\end{proof}

We now prove \Cref{thm:z1-L-stable-ansatz}.
\begin{proof}[Proof of \Cref{thm:z1-L-stable-ansatz}]
	By combining the third property of  \Cref{thm:z1-L-stable-ansatz} and invoking \cref{prop:stability-denominator}, the denominator of the stability function \cref{eq:ark-stability-function} is the polynomial  $Q(z_1,z_2)=\prod_{i=1}^{\scmp{1}}(1-z_1 \ahcmp{1}_{ii}$) whose leading coefficient for $z_1^s$ is $\prod_{i=1}^{\scmp{1}} \ahcmp{1}_{ii} \neq 0$. 
	
	The NPRK method will be L-stable in the stiff $z_1$ limit if the numerator $P(z_1,z_2)$ is a polynomial in $z_1$ of degree less than $\scmp{1}$.  The numerator of the stability function \cref{eq:ark-stability-function} is
	\begin{align}
		P(z_1,z_2) &= \det(\mathbf{I} - (z_1 \mathbf{A}^{\{1\}} + z_2 \mathbf{A}^{\{2\}}) + \mathbf{e}(z_1 \mathbf{b}^{\{1\}} + z_2 \mathbf{b}^{\{2\}})^\text{T} ) \\
		&= \det(\mathbf{I} - \mathbf{P}((z_1 \mathbf{A}^{\{1\}} + z_2 \mathbf{A}^{\{2\}}) + \mathbf{e}(z_1 \mathbf{b}^{\{1\}} + z_2 \mathbf{b}^{\{2\}})^\text{T} ) \mathbf{P}^T ).
	\end{align}
where $\mathbf{P}$ is the permutation matrix defined in \eqref{eq:underlying-permutation-matrix}. Since the second property of the theorem mandates a stiffly accurate method, it follows that: (1) the final rows of the unreduced methods satisfy $\AcmpB{r}_{s,j} = \bcmpB{r}_j$ and (2) the final stage of the NPRK method must belong to $\Scmp{1}$ (\cref{lem:stiff-accurate-final-stage-s1}). The second statement is important because the entries of the row $\bhcmpB{r}_{j}$	will be in the $(\scmp{1})^{\text{th}}$ row of the permuted matrices $\mathbf{P}\mathbf{A}^{\{r\}}\mathbf{P}^T$ (e.g. consider a modification of \cref{fig:underlying-method-matrix-structure} with $\Scmp{1}_3 = s$) and we will see that these terms cancel out when we subtract $\mathbf{P} (\mathbf{e} (\bcmpB{r})^T) \mathbf{P}^*$.	

To show this fact, we define $\cmp{\mathbf{B}}{r}(n) \in \mathbb{R}^{n \times \cmp{s}{r}}$ to be the matrix where each of its $n$ rows is the vector $\bhcmpB{r}$, where $\bhcmpB{r} \in \mathbb{R}^{\scmp{r}}$ is the reduced output vector for the method $\cmp{\widehat{M}}{r}$, $r=1,2$. Next, we consider the effects of the permutation $\mathbf{P}$ from \eqref{eq:underlying-permutation-matrix} applied to matrices $\AhcmpB{r}$, $\mathbf{e} (\bcmpB{r})^T$.
	
\begin{itemize}[leftmargin=2em]
	\item {\em Permuting $\AhcmpB{1}$, $\mathbf{e} (\bcmpB{1})^T$}: applying the permutation matrix \eqref{eq:underlying-permutation-matrix} yields
	\begin{align}
		\mathbf{P} \AcmpB{1} \mathbf{P}^* = 
		\renewcommand*{\arraystretch}{1.25}
		\left[\begin{array}{c|c}
			\cmp{\widehat{\mathbf{A}}}{1} & \mathbf{0} \\ \hline
			\cmp{\mathbf{C}}{1}_D & \mathbf{0}
		\end{array}\right],
		&&
\mathbf{P} \mathbf{e} (\bcmpB{1})^T \mathbf{P}^* = 
		\left[ \begin{array}{c|c}
			\cmp{\mathbf{B}}{1}(\cmp{s}{1}) & \mathbf{0} \\ \hline
			\cmp{\mathbf{B}}{1}(\cmp{s}{2}) & \mathbf{0}
		\end{array}
		\right].
	\end{align}
	Therefore, $\mathbf{P} (\AcmpB{1} - \mathbf{e} (\bcmpB{1})^T) \mathbf{P}^*$ has a zero $(\scmp{1})^{\text{th}}$ row:
	\begin{align} 
		\mathbf{P} (\AcmpB{1} - \mathbf{e} (\bcmpB{1})^T) \mathbf{P}^* = 
		\left[
		\begin{array}{c|c}
			\begin{array}{c}				\cmp{\widehat{\mathbf{A}}}{1}(1:\cmp{s}{1}-1,:) - \cmp{\mathbf{B}}{1}(\cmp{s}{1}-1) \\
				\mathbf{0} \\		
			\end{array} & 
			\mathbf{0} \\ \hline
			\cmp{\mathbf{C}}{1}_D - \cmp{\mathbf{B}}{1}(\cmp{s}{2}) & \mathbf{0} \\
		\end{array}
		\right].
	\end{align}	
	 	
	\item {\em Permuting $\AhcmpB{2}$, $\mathbf{e} (\bcmpB{2})^T$}: applying the permutation matrix \eqref{eq:underlying-permutation-matrix} yields
	\begin{align}
		\mathbf{P} \AcmpB{2} \mathbf{P}^* &= 
		\left[ \begin{array}{c|c}
			\cmp{\mathbf{C}}{2}_{1 \to 2} & \cmp{\mathbf{C}}{2}_{D} \\ \hline
			\mathbf{0} & \cmp{\widehat{\mathbf{A}}}{2}	
		\end{array}
		\right], &&
		\mathbf{P} \mathbf{e} (\bcmpB{2})^T \mathbf{P}^* = 
		\left[ \begin{array}{c|c}
			\mathbf{0} & \cmp{\mathbf{B}}{2}(\cmp{s}{1}) \\ \hline
			\mathbf{0} & \cmp{\mathbf{B}}{2}(\cmp{s}{2})
		\end{array}
		\right].
	\end{align}
	Stiff accuracy of the MR-NPRK method and $s\in\Scmp{1}$, implies that the final row of $\cmp{\mathbf{C}}{2}_{D}$ is equal to $\bhcmpB{2}$:
	\begin{align}
		\cmp{\mathbf{C}}{2}_{D} &= 
		\left[
		\begin{array}{c}
 			\cmp{\mathbf{C}}{2}_{D}(1:\cmp{s}{1}-1,:) \\ \hline
 			\bhcmpB{2}
 		\end{array}
 		\right].
	\end{align}
	Therefore, $\mathbf{P} (\AcmpB{2} - \mathbf{e} (\bcmpB{2})^T) \mathbf{P}^*$ is given by
	\begin{align} 
		\left[
		\begin{array}{c|c}
			\cmp{\mathbf{C}}{2}_{1 \to 2} & 
			\begin{array}{c}
				\cmp{\mathbf{C}}{2}_{D}(1:\cmp{s}{1}-1,:) - \cmp{\mathbf{B}}{2}(\cmp{s}{1}-1) \\
				\mathbf{0} \\		
			\end{array} \\ \hline
			\mathbf{0} & \cmp{\widehat{\mathbf{A}}}{2} - \cmp{\mathbf{B}}{2}(\cmp{s}{2}) \\
		\end{array}
		\right]
	\end{align}	
	\item {\em Combining results:} 	The key observation so far is that the $(\cmp{s}{1})^{th}$ row of the matrix $\mathbf{P} (\AcmpB{1} - \mathbf{e} (\bcmpB{1})^T) \mathbf{P}^*$ is zero and the entries $\scmp{1}+1, \ldots, s$ of the $(\cmp{s}{1})^{th}$ row of the matrix $\mathbf{P} (\AcmpB{2} - \mathbf{e} (\bcmpB{2})^T) \mathbf{P}^*$ are zero.
	
	The numerator is determined from the determinant of the matrix  
		\begin{align}
			&\mathbf{I} + z_1 \mathbf{P} \left[ \AcmpB{2} - \mathbf{e} (\bcmpB{1})^T \right] \mathbf{P}^* + z_2 \mathbf{P} \left[ \AcmpB{2} - \mathbf{e} (\bcmpB{2})^T \right] \mathbf{P}^*.
\end{align}
Using the previous results, we express the matrix in block form to highlight the $z_1$ and $z_2$ dependence: 
\begin{align}
			&\mathbf{I} +   
				\left[ \begin{array}{c|c}
					z_1 \cmp{\mathbf{N}}{1} + z_2 \cmp{\mathbf{N}}{2} & z_2 \mathbf{M} \\ 
					z_2 \mathbf{g} & \mathbf{0} \\ \hline
					z_1\cmp{\mathbf{C}}{1}_D & z_2 \cmp{\widehat{\mathbf{A}}}{2}
				\end{array}
				\right], \quad 
				\begin{aligned}
					\cmp{\mathbf{N}}{j} &\in \mathbb{R}^{(\scmp{1}-1) \times \scmp{1}}, \\
					\mathbf{M} &\in \mathbb{R}^{(\scmp{1}-1) \times \scmp{2}}, \\
					\mathbf{g} &\in \mathbb{R}^{1 \times \scmp{1}}.
				\end{aligned}
		\end{align}
		To conclude the proof, we use the block formula for the determinant \cref{eq:block-det-formula-z1}		along with the following observations:
		\begin{itemize}
			\item {\em Regarding $\det(D)$:} $D$ only depends on $z_2$, therefore $\det(D)$ is a polynomial in $z_2$ of degree at most $\scmp{2}$.
			\item {\em Regarding $\det(A - BD^{-1}D)$:} 
			\begin{itemize}
				\item Since $B$ and $D$ only depend on $z_2$ the matrix $BD^{-1}C$ has entries that are linear in $z_1$. 
				\item Since $B$ has a zero bottom row, the matrix $B D^{-1} C$ also has a zero bottom row.
				\item If the method is fully-decoupled, then the bottom row of the matrix $A$ is zero. If the method is $1 \to 2$ coupled, then the bottom row of $A$ only depends on $z_2$.
				\item In summary the coefficients of the matrix $A-BD^{-1}C$ depend linearly on $z_1$ and the bottom row does not depend on $z_1$. Therefore $\det(A-BD^{-1}C)$ is a polynomial in $z_1$ of degree at most $\scmp{1}-1$
			\end{itemize}	
		\end{itemize}
	In summary, we have shown that the assumptions required by \cref{thm:z1-L-stable-ansatz} guarantee that numerator and denominators of the stability function \cref{eq:ark-stability-function} are polynomials in $z_1$ of degree $\scmp{1}-1$ and $\scmp{1}$ respectively. Therefore, $\lim_{|z_1| \to \infty} S(z_1,z_2) = 0$ and the method is L-stable in the stiff $z_1$ limit.

\end{itemize}
\end{proof}

\bibliographystyle{siamplain}
\bibliography{references}



\def\siamprelabel{SM}
\setcounter{section}{0}
\setcounter{page}{1}


\begin{center}
	{\bfseries{\MakeUppercase{Supplementary Materials: Multirate Runge-Kutta for Nonlinearly Partitioned Systems}}}
	
	\vspace{1em}
	TOMMASO BUVOLI {\scriptsize AND} BRIAN K. TRAN {\scriptsize AND} BEN S. SOUTHWORTH	
	\vspace{2em}	
\end{center}


\section{Some facts regarding MR-NPRK methods}
Here, we provide details on some facts about MR-NPRK methods stated in \Cref{sec:mr-nprk}.

\subsection{Example of an NPRK method violating \Cref{assump:stage-union}}
\label{sup:nprk-example-violating-stage-union}

Here, we show a simple first-order NPRK method where $\Scmp{1} \cup \Scmp{2} \subset \mathcal{S}$ (i.e., violating \cref{assump:stage-union}). Specifically, the method
\begin{align}
	\begin{aligned}
		Y_1 &= y_n, \\
		Y_2 &= y_n + h \alpha F(Y_2,Y_1), \\
		y_{n+1} &= y_n + h F(Y_1,Y_2) + h F(Y_2,Y_1) - h F(Y_2,Y_2),	
	\end{aligned}
\end{align}
has stage sets $\mathcal{S} = \{1, 2\}$, and $\Scmp{1} = \Scmp{2} = \{1\}$. Although stage $Y_2$ is not reducible for the NPRK, it can be harmlessly removed in both the underlying ARK method and the individual integrators:  
\begin{align}
	 \text{Underlying ARK} &
	\begin{cases}
		Y_1 &= y_n, \\
		Y_2 &= y_n + h \alpha F^{\{1\}}(Y_2) + h \alpha F^{\{2\}}(Y_1), \\
		y_{n+1} &= y_n + h F^{\{1\}}(Y_1) + h F^{\{2\}}(Y_1).
	\end{cases} \\[1em]
	\Mcmp{1} &
	\begin{cases}
		Y_1 &= y_n, \\
		Y_2 &= y_n + h \alpha F^{\{1\}}(Y_2), \\
		y_{n+1} &= y_n + h F^{\{1\}}(Y_1).
	\end{cases} \\[1em]
		\Mcmp{2} &
	\begin{cases}
		Y_1 &= y_n, \\
		Y_2 &= y_n +  h \alpha F^{\{1\}}(Y_1), \\
		y_{n+1} &= y_n + h F^{\{2\}}(Y_1).
	\end{cases}
\end{align}

\subsection{Implied sparsity for decoupled, $1 \to 2$, and $2 \to 1$ coupled methods}\label{subsec:implied-sparsity}

Corresponding to \Cref{fig:timescale-coupling-diagrams}, we can precisely characterize each timescale coupling by specifying the allowable non-zero coefficients of the MR-NPRK tensors. For each of the timescale couplings, we define a corresponding \emph{allowable set}, which is a subset of $\mathcal{S} \times \mathcal{S} \times \mathcal{S}$, describing which stage coefficients $a_{ijk}$ \eqref{eq:nprk-general-a} are allowed to be non-zero:
	\begin{subequations} 
		\label{eq:stage-coefficient-sets}
		\begin{align}
	      	    \text{($\emptyset$) allowable set}: \quad &\mathcal{A}_{\emptyset} :=  \left\{i \in \mathcal{S}, j \in \Scmp{1}, k \in \Scmp{2} \right\},
	      		\label{eq:fully-decoupled-stage-conditions} \\
	      		\text{(i) allowable set}:  \quad &\mathcal{A}_{\emptyset} \cup \mathcal{A}_{1 \to 2} \text{ where } \mathcal{A}_{1 \to 2} := \left\{i,j,k \in \Scmp{1} \right\} \label{eq:12-decoupled-stage-conditions}, \\
	      		\text{(ii) allowable set}: \quad &\mathcal{A}_{\emptyset} \cup \mathcal{A}_{2 \to 1} \text{ where } \mathcal{A}_{2 \to 1} :=  \left\{i,j,k \in \Scmp{2}\right\},
	      		\label{eq:21-decoupled-stage-conditions} \\
                \text{(iii) allowable set}: \quad &\mathcal{A}_{\textup{full}} := \mathcal{S} \times \mathcal{S} \times \mathcal{S}.
		\end{align}
	\end{subequations}
	Similarly, for each of the timescale couplings, we define the corresponding \emph{output allowable set}, which is a subset of $\mathcal{S} \times \mathcal{S}$, describing which output coefficients $b_{ij}$ \eqref{eq:nprk-general-b} are allowed to be non-zero:
    	\begin{subequations}\label{eq:all-output-conditions}
	    	\begin{align}
				\text{($\emptyset$), (i), (ii) output allowable sets:} \quad &\mathcal{B}_{\emptyset} := \left\{i \in \Scmp{1}, j \in \Scmp{2} \right\},
		      		\label{eq:fully-decoupled-output-conditions} \\
                    \text{(iii) output allowable set:} \quad &\mathcal{B}_{\textup{full}} := \mathcal{S} \times \mathcal{S}.
			\end{align}
		\end{subequations}

\begin{proposition}\label{prop:sparsity}

MR-NPRK methods satisfy the sparsity structure shown in \Cref{tab:sparsity}.

\begin{proof}
The allowable sets $\mathcal{A}_\emptyset$, $\mathcal{A}_{1\to2}$, $\mathcal{A}_{2\to1}$ and $\mathcal{B}_\emptyset$ in \eqref{eq:stage-coefficient-sets} and \eqref{eq:all-output-conditions} define the indices for which the NPRK coefficients are allowed to be non-zero in decoupled, $1\to2$ and $2\to1$ methods. Thus, the coefficients which \textit{must} be zero, which we refer to as ``implied coefficient conditions" in \Cref{tab:sparsity}, have indices in the complement of the allowable sets.

For all three non fully-coupled methods, we have the complement of the output allowable set
$$ \mathcal{B}_{\emptyset}^c = (\Scmp{1} \times \Scmp{2})^c = (\mathcal{S}^{\{1\}c} \times \mathcal{S} ) \cup (\mathcal{S} \times \mathcal{S}^{\{2\}c}), $$
which yields the conditions $b_{\gamma k} = 0 = b_{j \nu}$ for all $\gamma \in \Scmp{2}$, $k,j \in \mathcal{S}$, $\nu \in \Scmp{1}$.

For fully-decoupled, $1\to 2$ and $2\to1$ methods, the complement of their stage allowable sets are, respectively,
\begin{align*}
        \mathcal{A}_{\emptyset}^c &= (\mathcal{S} \times \Scmp{1} \times \Scmp{2})^c = (\mathcal{S} \times \Scmp{2} \times \mathcal{S}) \cup (\mathcal{S} \times \mathcal{S} \times \Scmp{1}), \\
        (\mathcal{A}_{\emptyset} \cup \mathcal{A}_{1\to2})^c &= \mathcal{A}_{\emptyset}^c \cap \mathcal{A}_{1 \to 2}^c = (\mathcal{S} \times \Scmp{2} \times \mathcal{S}) \cup (\Scmp{2} \times \mathcal{S} \times \Scmp{1}), \\
    (\mathcal{A}_{\emptyset} \cup \mathcal{A}_{2\to1})^c &= \mathcal{A}_{\emptyset}^c \cap \mathcal{A}_{2 \to 1}^c = (\mathcal{S} \times \mathcal{S} \times \Scmp{1}) \cup (\Scmp{1} \times \Scmp{2} \times \mathcal{S}) ,
\end{align*}
which yields the remaining implied coefficient conditions on $a_{ijk}$ in \Cref{tab:sparsity}. \qed
\end{proof}
\end{proposition}

\subsection{Intersection of $\Scmp{1}$ and $\Scmp{2}$}

\begin{proposition}\label{prop:intersection-fully-coupled}
    If $S^{\{1\}}\cap S^{\{2\}} \neq \emptyset$, then $M$ is fully-coupled. 
    \begin{proof}
        This trivially follows from the definitions: if there is an element $k \in S^{\{1\}}\cap S^{\{2\}}$, then, in the NPRK method, $Y_k$ appears both in the first argument and second argument without being able to reduced out of either underlying method and hence, $M$ cannot be strictly $1 \to 2$, $2 \to 1$, nor fully-decoupled, i.e., $M$ is fully-coupled.
    \end{proof}
\end{proposition}
The converse of this statement is not true, there can be a fully-coupled $M$ whose irreducible stage sets have empty intersection. However, the contrapositive of \Cref{prop:intersection-fully-coupled} is naturally true, and gives a characterization of the intersection of the irreducible stage sets for methods which are not fully-coupled.
\begin{corollary}\label{cor:disjoint-not-fully-coupled}
    If $M$ is an MR-NPRK method with fully-decoupled timescales, $1 \to 2$ coupling, or $2 \to 1$ coupling, then the intersection of its irreducible stage sets is empty, $\Scmp{1}\cap \Scmp{2} =\emptyset$.
    \begin{proof}
        This follows from the fact that the negation of being fully-coupled, by definition, is being either $1 \to 2$ coupled, $2 \to 1$ coupled, or fully-decoupled.
    \end{proof} \qed
\end{corollary}

\section{Order conditions reductions}
\label{sup:order-conditions}

The order conditions of NPRK methods can be broadly divided into two types  \cite{nprk2}: those which are identical to the order conditions of the underlying ARK method, and those involving the NPRK tensors that cannot be written as ARK order conditions.

All nonlinear order conditions can reduced by simply applying the sparsity conditions from \Cref{tab:sparsity}. For the order conditions that are identical to those of the underlying ARK method, it is useful to also extend the sparsity to the underlying methods. Using \Cref{tab:sparsity}, the following simplifications follow for a fully-decoupled method: let $\nu \in \mathcal{S}^{\{1\}}$, $\gamma \in \mathcal{S}^{\{2\}}$, $i \in \mathcal{S}$, then \begin{subequations}
	\begin{align}
	    \acmp{1}_{i\gamma} &= \sum_{k=1}^s a_{i\gamma k} = 0, && \bcmp{1}_\gamma = \sum_{j=1}^s b_{\gamma j} = 0, \label{eq:fd-a1-sparsity}
	    \\
	    \acmp{2}_{i\nu} &= \sum_{j=1}^s a_{ij\nu} = 0, && \bcmp{2}_\nu = \sum_{k=1}^s b_{k \nu} = 0. \label{eq:fd-a2-sparsity}
	\end{align}
\end{subequations}
Methods with $1 \to 2$ coupling only satisfy \eqref{eq:fd-a1-sparsity} while methods with $2 \to 1$ coupling only satisfy \eqref{eq:fd-a2-sparsity}. Therefore, the following simplifications are available:
\begin{subequations}
\begin{align}
	&  \sum_{j \in \mathcal{S}} \acmp{1}_{ij} = \sum_{j \in \Scmp{1}} \acmp{1}_{ij}, &&
	  \sum_{i \in \mathcal{S}} \bcmp{1}_i = \sum_{i \in \Scmp{1}} \bcmp{1}_i, &&
	  c_i = \sum_{j,k \in \mathcal{S}} {a_{ijk}} =\sum_{\subalign{\nu &\in \Scmp{1}, \\ \gamma &\in \mathcal{S}}} {a_{i\nu \gamma}}, \label{eq:fd-a1-simp} \\
	&  \sum_{j \in \mathcal{S}} \acmp{2}_{ij} = \sum_{j \in \Scmp{2}} \acmp{2}_{ij}, &&
	  \sum_{i \in \mathcal{S}} \bcmp{2}_i = \sum_{i \in \Scmp{2}} \bcmp{2}_i, &&
	  c_i = \sum_{j,k \in \mathcal{S}} {a_{ijk}} =\sum_{\subalign{\nu &\in \mathcal{S}, \\ \gamma &\in \Scmp{2}}} {a_{i\nu \gamma}}, \label{eq:fd-a2-simp}
\end{align}
\end{subequations}
where $1 \to 2$ coupled methods satisfy \eqref{eq:fd-a1-simp} and $2 \to 1$ coupled methods satisfy \eqref{eq:fd-a2-simp}. Fully-decoupled methods satisfy both \eqref{eq:fd-a1-simp} and \eqref{eq:fd-a2-simp}, along with the further simplification 
	\begin{align} 
	    c_i = \textstyle \sum_{\substack{\nu \in \Scmp{1} \\ \gamma  \in \Scmp{2}}} {a_{i\nu \gamma}}.
	\end{align}

\section{Detailed method derivations}
Here, we provide details of the method derivation for the families of second- and third-order MR-NPRK methods presented in \Cref{subsec:method-ansatze}. 

\subsection{Derivation of second-order MR-NPRK methods}\label{sup:second-order-derivation}
	Here, we describe the derivation of the second-order MR-NPRK method from \Cref{sec:second-order-schemes}.
	We begin by proposng a second-order method ansatz that pads an $\scmp{2}$-stage explicit RK method ($\Ahcmp{2}$, $\bhcmp{2}$) with one implicit stages at the beginning and one at end of the method. The method ansatz is
\begin{align}
	\begin{aligned}
		Y_1 &= y_n, \\
		{\color{magenta} Y_2} &= y_n + ha_{2,2,1} F({\color{magenta} Y_2},Y_1), \\
		Y_i &= y_n + h \sum_{k=1}^{i-1} a_{i,2,k} F({\color{magenta} Y_2},Y_k), \hspace{4em} i = 3,\ldots, \scmp{2}+1, \\
		{\color{Cerulean} Y_{s}} &= y_n + h\sum^{s-1}_{k=0} a_{s,2,k} F({\color{magenta} Y_2},Y_k) + h\sum^{s-1}_{k=0} a_{s,s,k} F({\color{Cerulean} Y_s},Y_k), \\
		y_{n+1} &= {\color{Cerulean} Y_{s}},
	\end{aligned}
\end{align}
with $s=\scmp{2}+2$ and coefficients that satisfy the constraints
\begin{align}
	& \text{\em explicit stages} && \text{\em implicit stages} \nonumber \\ 
	& \left\{
	\begin{aligned}
		& a_{i,2,1} = \ahcmp{2}_{i-1,1} \\
		& a_{i,2,2} = 0 \\
		& a_{i,2,k} = \ahcmp{2}_{i-1,k-1}		
	\end{aligned} 
	\right.
	&&\left\{
	\begin{aligned}
		& a_{s,2,1} + a_{s,s,1} = \bhcmp{2}_1 \\
		& a_{s,2,2} = a_{s,s,2} = 0 \\
		& a_{s,2,k} + a_{s,s,k}	= \bhcmp{2}_{k-1}
	\end{aligned}
	\right.
	&& \text{for} &
	\begin{aligned}
		i &= 3,\ldots, \scmp{2}+1, \\
		k &= 3,\ldots, \scmp{2}+1.
	\end{aligned}
\end{align}
The resulting fully-decoupled, stiffly accurate method has $\scmp{2}+2$ total stages with only the second and final stage requiring implicit solves. The explicit stages have no free parameters, while there are a total of $s$ free parameters in the implicit stages. The constraints guarantee that $\Scmp{1} = \{ 2, s \}$, $\Scmp{2}=\{ 1, 3, \ldots, s-1 \}$ and the explicit reduced underlying method $\Mhcmp{2} = (\Ahcmp{2}, \bhcmp{2})$. 

Incorporating the constraints leads to the $s=\scmp{2}+2$ stage method
\begin{footnotesize}
\begin{align}
	\begin{aligned}
		Y_1 &= y_n, \\
		{\color{magenta} Y_2} &= y_n + ha_{2,2,1} G({\color{magenta} Y_2},Y_0), \\
		Y_i &= y_n + h \acmp{2}_{i-1,1} G({\color{magenta} Y_2},Y_0) + h \sum_{k=3}^{i-1} \acmp{2}_{i-1,k} G({\color{magenta} Y_2},Y_k), \quad i = 2,\ldots, \scmp{2}+1,\\
		{\color{Cerulean} Y_{s}} &= y_n 
			+ h\left[ (\bcmp{2}_1 - a_{s,s,1}) G({\color{magenta} Y_2},Y_1) + (a_{s,s,1}) G({\color{Cerulean} Y_{s}},Y_1) \right]
			\\
            & \hspace{2.5em} + h\sum^{\hat{s}}_{k=3} \left[ (\bcmp{2}_{k-1} - a_{s,s,k}) G({\color{magenta} Y_2},Y_k) + (a_{s,s,k}) G({\color{Cerulean} Y_{s}},Y_k) \right], \\
		y_{n+1} &= {\color{Cerulean} Y_{s}}.
	\end{aligned}
	\label{eq:implicit-padding-stiffly-accurate}
\end{align}
\end{footnotesize}\ignorespaces
The resulting reduced underlying implicit method $\Mhcmp{1}$ is
\begin{align}
	\widehat{M}^{\{1\}} \quad
	\begin{tabular}{c|cc}
		$a_{2,2,1}$ & $a_{2,2,1}$ \\
		$1$ 		& $1 - \gamma$ & $\gamma$  \\ \hline
		    		& $1 - \gamma$ & $\gamma$ 
	\end{tabular},	
	&&
	\gamma = \sum_{\substack{k=1 \\ k \ne 2}}^{s-1} a_{s,s,k}.
	\label{eq:implicitly-wrappedp-second-order-m1}
\end{align}
Note that in \eqref{eq:implicitly-wrappedp-second-order-m1}, we have assumed that the method $\Mhcmp{2}$ is at least first-order accurate such that $\sum_{i=1}^{\scmp{2}} \bhcmp{2}_i = 1$. For a second-order MR-NPRK method, it is necessary and sufficient that both underlying methods are second-order (this is a general property of NPRK methods, see \cite{nprk2}). Therefore, the order conditions are
\begin{align}
	\sum_i \bhcmp{1}_i &= 1, &
	\sum_i \bhcmp{2}_i &= 1, &
	\sum_i \bhcmp{1}_i \chcmp{1}_i &= \tfrac{1}{2}, &
	\sum_i \bhcmp{2}_i \chcmp{2}_i &= \tfrac{1}{2},
	\label{eq:second-order-implicit-padding}
\end{align}
where $(\bhcmp{1}, \Ahcmp{1},\chcmp{1})$ and $(\bhcmp{2},\Ahcmp{2},\chcmp{2})$ are the reduced tableaus of the underlying implicit and explicit integrators. If the method $(\bhcmp{2},\Ahcmp{2},\chcmp{2})$ is second-order accurate, then all but the third conditions are automatically satisfied. Using $\sum_i \bcmp{1}_i = \sum_i \bcmp{2}_i = 1$, we have
\begin{align}
	\begin{aligned}
		\bhcmp{2} &= 
		\begin{bmatrix}
			(\bhcmp{2}_1 - a_{s,s,1}) + \sum_{k=3}^{\hat{s}} (\bhcmp{2}_k - a_{s,s,k}) \\
			(a_{s,s,1}) + \sum_{k=3}^{\hat{s}} (a_{s,s,k}) \\
		\end{bmatrix} =
		\begin{bmatrix}
			1 - \gamma \\
			\gamma \\
		\end{bmatrix}, && \gamma = \sum_{\substack{k=1 \\ k \ne 2}}^{\hat{s}} a_{s,s,k}, \\	
		\chcmp{2} &= 
		\begin{bmatrix}
			a_{2,2,1} \\
			1
		\end{bmatrix}. &&
	\end{aligned}
	\label{eq:second-order-implicit-padding-implicit-bc}
\end{align}
Assuming that $\cmp{\widehat{M}}{2}$ is second-order, the only remaining order condition is
\begin{align}
	\sum_i \bcmp{1}_i \cmp{c}{1}_i &= (1 - \gamma)a_{2,2,1} + \gamma = \tfrac{1}{2}
	\quad \implies \quad a_{2,2,1} = \frac{1 - 2\gamma}{2-2\gamma}.
	\label{eq:second-order-implicit-padding}
\end{align}
The resulting method family has $s-2$ free parameters -- see the definition of $\gamma$ in \eqref{eq:implicitly-wrappedp-second-order-m1}. One possible choice that leads to a simple implicit solve is $a_{s,s,j} = 0$, $j=1,\ldots, s-2$, which implies $\gamma = a_{s,s,s-1}$. Furthermore, imposing a SDIRK structure for $\cmp{\widehat{M}}{1}$ leads to $a_{s,s,s-1} = a_{2,2,1} = (2 \pm \sqrt{2}) / 2$. This then yields the class of second-order MR-NPRK methods \eqref{eq:mr-nprk-second-order-minus} presented in \Cref{sec:second-order-schemes}.

\subsection{Derivation of third-order MR-NPRK methods}\label{sup:third-order-derivation}
Here, we provide a detailed derivation of the third-order MR-NPRK methods presented in \Cref{sec:third-order-schemes}. We propose a third-order method ansatz that pads an $\scmp{2}$-stage explicit RK method ($\cmp{\widehat{A}}{2}$, $\cmp{\widehat{b}}{2}$) with two implicit stages at the beginning and one at end of the method. The method ansatz is
\begin{align}
	\begin{aligned}
		Y_1 &= y_n, \\
		{\color{magenta} Y_2} &= y_n + ha_{2,2,1} F({\color{magenta} Y_2},Y_1), \\
		{\color{ForestGreen} Y_3} &= y_n + h \left[ a_{3,2,1} F({\color{magenta} Y_2},Y_1) + a_{3,2,2} F({\color{magenta} Y_2},{\color{magenta} Y_2}) + a_{3,3,1} F({\color{ForestGreen} Y_3},Y_1) + a_{3,3,2} F({\color{ForestGreen} Y_3},{\color{magenta} Y_2}) \right], \\
		Y_i &= y_n + h \sum_{k=1}^{i-1} \left[ a_{i,2,k} F({\color{magenta} Y_2},Y_k) + a_{i,3,k} F({\color{ForestGreen} Y_3},Y_k) \right], \hspace{4em} i = 4,\ldots, \scmp{2}+2, \\
		{\color{Cerulean} Y_{s}} &= y_n +  h\sum^{s-1}_{k=1} \left[a_{s,2,k} F({\color{magenta} Y_2},Y_k) + a_{s,3,k} F({\color{ForestGreen} Y_3},Y_k) + a_{s,s,k} F({\color{Cerulean} Y_s},Y_k) \right], \\
		y_{n+1} &= {\color{Cerulean} Y_{s}},
	\end{aligned}
\end{align}
with $s=\scmp{2}+3$ and coefficients that satisfy the constraints
\begin{align}
	& \text{\em explicit stages} && \text{\em implicit stages} \nonumber \\ 
	& \left\{
	\begin{aligned}
		& a_{i,2,1} + a_{i,3,1} = \ahcmp{2}_{i-3,1} \\
		& a_{i,2,2} = a_{i,3,2} = a_{i,3,3} = 0 \\
		& a_{i,2,k} + a_{i,3,k} = \ahcmp{2}_{i-3,k-2}		
	\end{aligned} 
	\right.
	&&\left\{
	\begin{aligned}
		& a_{s,2,1} + a_{s,3,1} + a_{s,s,1} = \bhcmp{2}_1 \\
		& a_{s,2,2} = a_{s,2,3} = a_{s,3,2} = a_{s,3,3} = a_{s,s,2} = a_{s,s,3} = 0 \\
		& a_{s,2,k} + a_{s,3,k} + a_{s,s,k}	= \bhcmp{2}_{k-2}
	\end{aligned}
	\right.
\end{align}
for $i,k = 4,\ldots, \scmp{2}+2$. This $1\to2$ coupled, stiffly accurate method has $\scmp{2}+3$ total stages with only the second, third, and final stage requiring implicit solves. The $i^{th}$ explicit stage $2(i-3)$ free parameters, while there are a total of $s$ free parameters in the implicit stages. The constraints guarantee that $\Scmp{1} = \{ 2, 3, s \}$, ${\Scmp{2}=\{ 1, 4, \ldots, s-1 \}}$, and that the underlying methods have tableaux: $\Mhcmp{2} = (\Ahcmp{2}, \bhcmp{2})$ and
\begin{align}
	\widehat{M}^{\{1\}} =
	\begin{tabular}{c|ccc}
		$\widehat{c}_1$ 			& $a_{2,2,1}$ \\
		$\widehat{c}_2$ 	& $(a_{3,2,1}+a_{3,2,2})$ & $(a_{3,3,1}+a_{3,3,2})$ \\ 
		$1$ & $1 - \gamma_3 - \gamma_s $ & $\gamma_3$ & $\gamma_s$ \\  \hline
		$1$ & $1 - \gamma_3 - \gamma_s$ & $\gamma_3$ & $\gamma_s$ 
	\end{tabular},
	&&
	\gamma_i = \sum_{\substack{k=1 \\ k \ne 2,3}}^{s-1} a_{s,i,k}.
	\label{eq:implicitly-wrapped-third-order-m1}
\end{align}
Note that in \eqref{eq:implicitly-wrapped-third-order-m1} we have implicitly assumed that the method $\Mcmp{2}$ is at least first-order accurate such that $\sum_{i} \bcmp{2}_i = 1$ (see the following remark).

\begin{remark}
	Note that the condition $\sum_{k=0}^{s-1} a_{2,3,k} = \acmp{1}_{31} = 1 - \acmp{1}_{32} - \acmp{1}_{33}$ is automatically satisfied if the explicit method is second order. Specifically, using
		\begin{align}
			\quad \sum^{s-1}_{\substack{k=1 \\ k \ne 2,3}} \bcmp{2}_{i-2,\max(k-2,1)} = \sum_{i=1}^{\scmp{2}} \bcmp{2}_i = 1		\end{align}
		it follows that			
		\begin{align}
			\begin{aligned}
				\acmp{1}_{31} &= \sum_{k=0}^{s-1} a_{2,3,k} = \sum^{s-1}_{\substack{k=1 \\ k \ne 2,3}}(\bcmp{2}_{i-2,\max(k-2,1)} - a_{i,3,k} - \delta_{i,\omega} \acmp{1}_{33}) \\ 
			&=1-\sum^{s-1}_{\substack{k=1 \\ k \ne 2,3}}a_{i,3,k} - \acmp{1}_{33} = 1 - \acmp{1}_{32} - \acmp{1}_{33}.	
			\end{aligned}
		\end{align} 
\end{remark}

 This MR-NPRK method has many free parameters; to reduce the number of free parameters, we make several immediate simplifications:
\begin{itemize}
	\item {\em Simple implicit solves}. We choose to sparsify the coefficients associated with the underlying implicit method so that the third stage has a solve of the form $Y_3 = b + G(Y_3, Y_{\omega})$, and the $s^{th}$ stage has a solve $Y_s = b + G(Y_s, Y_{\nu})$; the associated conditions are:
		\begin{align}
			a_{3,3,j} &= 0 \quad \text{for} \quad j \ne \omega, && \omega \in \{1,2\}, \\
			a_{s,s,j} &= 0 \quad \text{for} \quad j \ne \nu	, &&
				\nu \in \{1,4,\ldots, s-1\}.
		\end{align}
	\item {\em Force the underlying implicit method to be L-stable}. The only L-stable, three-stage, stiffly-accurate SDIRK method \cite{butcher2016numerical}[p.274] is 
	\begin{align}
			\begin{tabular}{l|lll}
				$\lambda$					& $\lambda$ \\
				$\tfrac{1}{2}(1 + \lambda)$	& $\tfrac{1}{2}(1 - \lambda)$ & $\lambda$ \\
				$1$  						& $\tfrac{1}{4}(-6 \lambda^2 + 16 \lambda - 1)$ & $\tfrac{1}{4}(6 \lambda^2 - 20 \lambda + 5)$ & $\lambda$ \\ \hline
				 							& $\tfrac{1}{4}(-6 \lambda^2 + 16 \lambda - 1)$ & $\tfrac{1}{4}(6 \lambda^2 - 20 \lambda + 5)$ & $\lambda$
			\end{tabular}
			\label{eq:stiffly-accurate-three-stage-SDIRK}
	\end{align}
	where $\lambda$ is the root of the polynomial $p(z) = \tfrac{1}{6} - \tfrac{3}{2} \lambda + 3 \lambda^2 - \lambda^3$, which to 32 digits of precision is $\lambda=0.43586652150845899941601945119356$. The resulting conditions are
	\begin{align*}
		a_{2,2,1} = a_{3,3,\omega} = a_{s,s,\nu} = \lambda, && a_{3,2,1} + a_{3,2,2} = \tfrac{1}{2}(1+\lambda), \\
		\sum_{\substack{k=1 \\ k \ne 2,3}}^{s-1} a_{s,3,k}=\tfrac{1}{4}(6 \lambda^2 - 20 \lambda + 5).
	\end{align*}
\end{itemize}
Denoting the entires of \cref{eq:stiffly-accurate-three-stage-SDIRK} using $\ahcmp{1}_{i,j}$, the resulting method ansatz after imposing both the conditions in the previous two bullets is:
\begin{align}
	\begin{aligned}
		Y_1 &= y_n, \\
		{\color{magenta} Y_2} &= y_n + h\ahcmp{1}_{11} F({\color{magenta} Y_2},Y_1), \\
		{\color{ForestGreen} Y_3} &= y_n + h \left[ (\ahcmp{1}_{21} - a_{3,2,2}) F({\color{magenta} Y_2},Y_1) + a_{3,2,2} F({\color{magenta} Y_2},{\color{magenta} Y_2}) + \ahcmp{1}_{22} F({\color{ForestGreen} Y_3},Y_\omega)  \right], \\
		Y_i &= y_n + h \sum_{\substack{k=1 \\ k \ne 2,3}}^{i-1} \left[ (\ahcmp{2}_{i-2,\max(k-2,1)} - a_{i,3,k}) F({\color{magenta} Y_2},Y_k) + a_{i,3,k} F({\color{ForestGreen} Y_3},Y_k) \right], \\
        & \qquad \qquad i = 4,\ldots, \scmp{2}+2, \\
		{\color{Cerulean} Y_{s}} &= y_n +  h\sum^{s-1}_{\substack{k=1 \\ k \ne 2,3}} \left[ (\bhcmp{2}_{\max(k-2,1)} - a_{s,3,k} - \delta_{s,\nu} \ahcmp{1}_{33}) F({\color{magenta} Y_2},Y_k) + a_{s,3,k} F({\color{ForestGreen} Y_3},Y_k)   \right] \\
        & \hspace{2.4em} + h \ahcmp{1}_{33} F({\color{Cerulean} Y_s},Y_\nu), \\
		y_{n+1} &= {\color{Cerulean} Y_{s}}.
	\end{aligned}
\end{align}
with $s = \scmp{2}+3$, $\delta_{i,j} = \begin{cases}0 & i \ne j \\ 1 & i = j \end{cases}$ where the coefficeints $a_{s,3,k}$ must be chosen to satisfy
\begin{align}
	\displaystyle \sum_{\substack{k=1 \\ k\ne2,3}}^{s-1} a_{s,3,k} = \ahcmp{2}_{32}.	
	\label{eq:remaining-condition-for-m1}
\end{align}

With the assumptions placed thus far, the remaining free parameters are:
\begin{center}
	\begin{tabular}{lll}
		stages 						& num free variables & free variables \\
		{\color{ForestGreen} $Y_3$} & $1$					 & $a_{321}$ \\
		$Y_4,\ldots,Y_{s-1}$ & $\scmp{2}(\scmp{2}-1)/2$					 & $a_{i3k}$, $i=4,\ldots,s-1$, $k = 1,\ldots,i-1$ \\
		{\color{Cerulean} $Y_s$} & $\scmp{2} - 1$			& $a_{s3k}$ all $k \neq 1$; $a_{s,3,1}$ used to guarantee \\
		& & condition \cref{eq:remaining-condition-for-m1}.
	\end{tabular}
\end{center}
If both the reduced implicit method ($\Ahcmp{1}$, $\bhcmp{1}$) and reduced explicit method ($\Ahcmp{2}$, $\bhcmp{2})$ are third-order, then only three remaining order-conditions need to be satisfied:
\begin{align}
	\bcmp{1}_i \Acmp{2}_{ij} c_j = \frac{1}{6}, && 
	\bcmp{2}_i \Acmp{2}_{ij} c_j = \frac{1}{6}, &&
	c_i b_{ij} c_j = \frac{1}{3}.
\end{align}
Since the underlying implicit and explicit methods are fully specified, the vectors $c_i$, $\bcmp{1}$, $\bcmp{2}$ are completely determined and have no free parameters. Therefore, to satisfy the order conditions we require least one degree of freedom in each of the matrices $\Acmp{1}$, $\Acmp{2}$, and $b_{ij}$. By writing out the matrices one obtains the following requirements
\begin{subequations}
	\begin{align}
		\bcmp{1}_i \Acmp{2}_{ij} c_j &= \frac{1}{6} && \text{requires a specific value for $a_{322}$}, \label{eq:third-oc-12} \\
		\bcmp{2}_i \Acmp{1}_{ij} c_j &= \frac{1}{6} && \text{requires at least one free variable in explicit stages}, \label{eq:third-oc-21} \\
		c_i b_{ij} c_j &= \frac{1}{3} && \text{requires at least one free variable in final stage}. \label{eq:third-oc-nonlinear}
	\end{align}
\end{subequations}
Furthermore, each order condition can be solved independently, meaning that there is no interdependence between the unknowns. We now describe our proposed strategy for satisfying these three order conditions.

\textbf{Order condition \eqref{eq:third-oc-12}}. This order condition simplifies to
\begin{align}
	&\omega = 1: & \chcmp{1} \ahcmp{1}_{32} a_{3,2,2} + \frac{\ahcmp{1}_{33}}{2} &= \frac{1}{6} \quad \implies \quad
	a_{3,2,2} = \frac{1 - 3\ahcmp{1}_{33}}{6 \chcmp{1} \ahcmp{1}_{32}}, \\
	&\omega = 2: & \chcmp{1} \ahcmp{1}_{32}	(\ahcmp{1}_{22} + a_{3,2,2}) + \frac{\ahcmp{1}_{33}}{2} &= \frac{1}{6} \quad \implies \quad
	a_{3,2,2} = \frac{1 - 3\ahcmp{1}_{33}}{6 \chcmp{1} \ahcmp{1}_{32}} - \ahcmp{1}_{22}.
\end{align}
Notice that $a_{3,2,1}$ exclusively depends on the reduced underlying implicit method. When $\Mhcmp{2}$ is given by  \cref{eq:stiffly-accurate-three-stage-SDIRK}, we have
\begin{align}
	& \omega = 1: & a_{3,2,2} &= \frac{2-6 \lambda }{18 \lambda ^3-60 \lambda ^2+15 \lambda } \approx 0.1825367204468751798095174531383,\\
	& \omega = 2: & a_{3,2,2} &=\frac{2-6 \lambda }{18 \lambda ^3-60 \lambda ^2+15 \lambda }-\lambda \approx -0.2533298010615838196065019980553.
\end{align}

\textbf{Order condition \eqref{eq:third-oc-21}}. This order condition requires at least one free variable in the explicit stage. Our aim is to ensure that to the largest degree possible the explicit stage computation equivalent to time-stepping a modified ODE $w'=G(w)$ with the underlying explicit method. There are three possible strategies for achieving this. Specifically, for $i=3,\ldots, s-1$ we can impose:
    \begin{subequations}
    \label{conditions-third-a-sparisfy}
	\begin{align}
		& \text{option 1:} \quad a_{i,3,k} = 0, \\
		& \text{option 2:} \quad a_{i,3,k} = \ahcmp{2}_{i-2,\max(k-2,1)} ~\implies a_{i,2,k}=0, \\
		& \text{option 3:} \quad a_{i,3,k} = \delta \ahcmp{2}_{i-2,\max(k-2,1)} \implies a_{i,2,k} = (1-\delta) \ahcmp{2}_{i-2,\max(k-2,1)}.
	\end{align}
\end{subequations}
  Corresponding to each option is the choice of the function $G$:
\begin{subequations}
    \begin{align}
		& \text{option 1:} \quad G(w) = F( {\color{magenta} Y_2}, w),\\
		& \text{option 2:} \quad G(w) = F({\color{ForestGreen} Y_3}, w),\\
		& \text{option 3:} \quad G(w) = (1 - \delta) F( {\color{magenta} Y_2}, w) + \delta F( {\color{ForestGreen} Y_3}, w).
    \end{align}
\end{subequations}
	Note that for options 1 and 2, it is necessary to violate the assumptions in at least one stage or it will not be possible to satisfy the order conditions since $\Acmp{1}$ will have no free parameters. In the final option, the parameter $\delta$ can be used to satisfy the order condition.
	
	\begin{itemize}
		\item {\em Option 1:} leave a single $a_{i,3,k}$ free in an explicit stage free and zero-out all others:
		\begin{align}
			a_{i,3,k} = \begin{cases}
					\eta & i=u, k=w \\
 					0 & \text{otherwise}
 				\end{cases}
		\end{align}
		for $u \in \{4, \ldots, \scmp{2}+2\}$ and $w \in \{1\} \cup \{ 4, \ldots, u - 1\}$. The resulting order condition is
		\begin{align}
			& \chcmp{1}_1 \left(\frac{1}{2} - \bhcmp{2}_{u-2} \eta \right) + \chcmp{1}_2 \left( \bhcmp{2}_{u-2} \eta \right) = \frac{1}{6} \\
				& \implies \quad \eta = \frac{1 - 3\chcmp{1}_1}{6	\bhcmp{2}_{u-2}(\chcmp{1}_2 - \chcmp{1}_1)}
		\end{align}

		\item {\em Option 2:} leave a single $a_{i,3,k}$ free in an explicit stage free and zero-out other $a_{i,2,k}$:
			\begin{align}
			a_{i,3,k} &= 
				\begin{cases}
					\eta & i=u, k=w \\
 					\ahcmp{2}_{i-2,\max(k-2,1)} & \text{otherwise}
 				\end{cases}  \nonumber
 				\\ & \implies 
 				a_{i,2,k} = 
				\begin{cases}
					\ahcmp{2}_{i-2,\max(k-2,1)} - \eta & i=u, k=w \\
 					0 & \text{otherwise}
 				\end{cases}
			\end{align}
			The resulting order condition is
			\begin{align}
				\frac{1}{6} &= \chcmp{1}_1 \left(\bhcmp{2}_{\max(1,\kappa-2)} (\ahcmp{2}_{u-2,\max(w-2,1)} - \eta) \right) \nonumber \\ & \qquad + \chcmp{1}_2 \left( \tfrac{1}{2} - \bhcmp{2}_{\max(1,\kappa-2)} (\ahcmp{2}_{u-2,\max(w-2,1)} - \eta) \right),  \\
					& \implies \quad \eta = \frac{1 - 3\chcmp{1}_2}{6	\bhcmp{2}_{u-2}(\chcmp{1}_2 - \chcmp{1}_1)} + \ahcmp{2}_{u-2,\max(w-2,1)}.
			\end{align}
 
 		\item {\em Option 3:} convex combination: 
			\begin{align}
			a_{i,3,k} = \delta \ahcmp{2}_{i-2,\max(k-2,1)},
			&&
 				a_{i,2,k} = (1-\delta) \ahcmp{2}_{i-2,\max(k-2,1)}.
			\end{align}
			The resulting order condition is 
			\begin{align}
				(1-\delta) (\tfrac{1}{2}) \chcmp{1}_1+\delta(\tfrac{1}{2}) \chcmp{1}_2 = \frac{1}{6} 
				\quad \implies \quad 
				\delta = \frac{1-3\chcmp{1}_1}{3(\chcmp{1}_2 - \chcmp{1}_1)}.
			\end{align}
			Once again this order condition only depends on the reduced underlying implicit method. When $\Mhcmp{2}$ is given by \eqref{eq:stiffly-accurate-three-stage-SDIRK}, we have
				\begin{align}
					\delta = \frac{2 (3 \lambda -1)}{3 (\lambda -1)} \approx -0.36350683689006809456061097715969.
				\end{align}
	\end{itemize}

\textbf{Order condition \eqref{eq:third-oc-nonlinear}}. This order condition requires at least one free variable in the NPRK output weight matrix $b_{ij}$. One possibility is to again preserve, to the largest degree possible, that the explicit pieces of the final implicit stage are computed as if explicitly timestepping. This leads to conditions on the final stage that are very similar to those for the internal stages $a_{ijk}$ tensor \cref{conditions-third-a-sparisfy}, namely:
    \begin{subequations}
	    \label{conditions-third-as-sparisfy}
		\begin{align}
			& \text{option 1:} \quad a_{s,3,k} = 0, \\
			& \text{option 2:} \quad a_{s,3,k} = \ahcmp{2}_{i-2,\max(k-2,1)} ~\implies a_{s,2,k}=0, \\
			& \text{option 3:} \quad a_{s,3,k} = \delta \ahcmp{2}_{i-2,\max(k-2,1)} \implies a_{s,2,k} = (1-\delta) \ahcmp{2}_{s-2,\max(k-2,1)}.
		\end{align}
	\end{subequations}
    Corresponding to each option is the choice of the function $G$:
\begin{subequations}
    \begin{align}
		& \text{option 1:} \quad G(w) = F( {\color{magenta} Y_2}, w),\\
		& \text{option 2:} \quad G(w) = F({\color{ForestGreen} Y_3}, w),\\
		& \text{option 3:} \quad G(w) = (1 - \delta) F( {\color{magenta} Y_2}, w) + \delta F( {\color{ForestGreen} Y_3}, w).
    \end{align}
\end{subequations}
	In order to specify $\Mhcmp{1}$ one of the final output coefficients must be chosen to satisfy \cref{eq:remaining-condition-for-m1}. This implies that we need to violate the conditions (options 1-3 listed above) at least twice or there will be no degrees of freedom in $b_{ij}$ and the order condition cannot be satisfied. Moreover, if the free variables has index $k \ne \nu$ ($\nu$ is used in implicit solve $Y_s=b+F(Y_s,Y_\nu)$) then the assumptions will be violated three times. 
	\begin{remark}
		If ($\Ahcmp{2}$, $\bhcmp{2}$) retains low storage cost as $\scmp{2} \to \infty$ (e.g. a repeated RK tableau), then one main advantage of forcing the explicit stages to resemble explicit time-stepping is that the MR-NPRK method retains the low storage property. We note that one can still retain low-storage for the MR-NPRK method even if the output stage does not satisfy the conditions listed above. Specifically, a method can still be implemented in a low-storage way (at the cost of one extra register) even if the structure of the final output stage does not match with the explicit stages. Specifically, if we let $b=y_n$ then as soon as the stage $Y_k$ is computed we incrementally add to b according to the rule:
			\begin{align}
				b \mathrel{+}= 	\left[ (\bhcmp{2}_{\max(k-2,1)} - a_{s,3,k} - \delta_{s,\nu} \ahcmp{1}_{33}) F({\color{magenta} Y_2},Y_k) + a_{s,3,k} F({\color{ForestGreen} Y_3},Y_k)   \right].
			\end{align}
			Then to compute the output we solve the implicit system $Y_s = b + \ahcmp{1}_{33} F(Y_s,Y_\nu)$. Nevertheless, to simplify the implementation we only consider the three options listed above.
	\end{remark}
	\begin{itemize}
		\item {\em Option 1.} Leave a single $b_{3,k}=a_{s,3,k}$ free in an explicit stage, select another to enforce $\Mhcmp{1}$ and then zero-out all others:
			\begin{align}
				a_{s,3,k} = 
					\begin{cases}
						\eta & k=\alpha \\
						\ahcmp{1}_{32} - a_{s,3,\alpha} & k=\beta \\
	 					0 & \text{otherwise}
	 				\end{cases}
			\end{align}
			The resulting order condition is
			\begin{align}
				\begin{aligned}
					\chcmp{1}_1 \left(\frac{1}{2} - a_{s,3,\alpha} \chcmp{2}_{\max(1,\alpha-2)} - (\ahcmp{1}_{32}-a_{3,s,\alpha}) \chcmp{2}_{\max(1,\beta-2)} - \ahcmp{2}_{33} \chcmp{2}_{\max(1,\nu-2)} \right) \\ +
				\chcmp{1}_2 \left( a_{s,3,\alpha} \chcmp{2}_{\max(1,\alpha-2)} + (\ahcmp{1}_{32}-a_{3,s,\alpha}) \chcmp{2}_{\max(1,\beta-2)} \right) \\ +
				\acmp{2}_{33} \chcmp{2}_{\max(1,\nu-2)}
				= \tfrac{1}{3}
				\end{aligned}
			\end{align}
			this leads to
			\begin{align}				
				\begin{tabular}{c}
						$a_{s,3,\alpha} = \frac{\tfrac{1}{3} - \ahcmp{1}_{33} \chcmp{2}_{\max(1,\nu-2)} + \chcmp{1}_{1}(\ahcmp{1}_{33} \chcmp{2}_{\max(1,\nu-2)} - \tfrac{1}{2}) - \ahcmp{1}_{32}(\chcmp{1}_2 - \chcmp{1}_1)\chcmp{2}_{\max(1,\beta-2)}}{(\chcmp{1}_2-\chcmp{1}_1)(\chcmp{2}_{\max(1,\alpha-2)}-\chcmp{2}_{\max(1,\beta-2)})}$
				\end{tabular}				               
			\end{align}

		\item {\em Option 2.} Leave a single $b_{3,k}=a_{s,3,k}$ free in an explicit stage, select another to enforce $\Mhcmp{1}$ and the rest to zero-out all $a_{s,2,k}$ except for $k = \nu$:
			\begin{align}
				a_{s,3,k} = 
					\begin{cases}
						\eta & k=\alpha \\
						\ahcmp{1}_{32} - a_{s,3,\alpha} - (1 - \bhcmp{2}_{\max(1,\alpha-2)} - \bhcmp{2}_{\max(1,\beta-2)}) & k=\beta \\
	 					\bhcmp{2}_{\max(1,k-2)} & \text{otherwise}
	 				\end{cases}
			\end{align}
			The resulting order condition is
			\begin{align*}
				\frac{1}{3} &= \chcmp{2}_1 \Big( (\bhcmp{2}_{\max(1,\alpha-2)} - a_{s,3,\alpha}) \chcmp{2}_{\max(1,\alpha-2)} \\
                & \qquad \qquad + (\bcmp{2}_{\max(1, \beta-2)} - D)\chcmp{2}_{\max(1,\beta-2)} - \ahcmp{1}_{33} \chcmp{2}_{\max(1,\nu-2)} \Big) \\
				& \qquad + \chcmp{2} \Big( \tfrac{1}{2} - \bhcmp{2}_{\max(1,\alpha-2)}\chcmp{2}_{\max(1,\alpha-2)} - \bhcmp{2}_{\max(1,\beta-2)}\chcmp{2}_{\max(1,\beta-2)} \\ 
                & \qquad \qquad + a_{s,3,\alpha} \chcmp{2}_{\max(1,\alpha-2)} + D \chcmp{2}_{\max(1,\beta-2)} \Big) + \ahcmp{1}_{33} \chcmp{2}_{\max(1,\nu-2)},
			\end{align*}
			for $D = \left(\ahcmp{1}_{32} - a_{s,3,\gamma} - \left(1 - \bhcmp{2}_{\max(1, \alpha-2)} - \bhcmp{2}_{\max(1,\beta-2)}\right)\right)$. This leads to to the choice
			\begin{align*}
				a_{s,3,\alpha} &= \frac{1}{(\chcmp{1}_2 - \chcmp{1}_1)(\chcmp{2}_{\max(1,\alpha-2)} - \chcmp{2}_{\max(1,\beta-2)})} \left( \frac{1}{3} - T_1 - T_2\right) \\
                & \qquad \qquad - \frac{T_3}{(\chcmp{2}_{\max(1,\alpha-2)} - \chcmp{2}_{\max(1,\beta-2)})},
            \end{align*}
            where $T_1, T_2, T_3$ are defined as
            \begin{footnotesize}
            \begin{align*}
				T_1 &:= \ahcmp{1}_{33} \chcmp{2}_{\max(1,\nu-2)} \\
                & \qquad - \chcmp{1}_1 \left( \bcmp{2}_{\max(1,\alpha-2)} \chcmp{2}_{\max(1,\alpha-2)} + \bcmp{2}_{\max(1,\beta-2)}  \chcmp{2}_{\max(1,\beta-2)} - \ahcmp{1}_{33}  \chcmp{2}_{\max(1,\nu-1)}\right),\\
				T_2 &:= \chcmp{2} \left( \tfrac{1}{2} - \bhcmp{2}_{\max(1,\alpha-2)}\chcmp{2}_{\max(1,\alpha-2)} - \bhcmp{2}_{\max(1,\beta-2)}\chcmp{2}_{\max(1,\beta-2)} \right),\\
				T_3 &:= \left(\ahcmp{1}_{32} - \left(1 - \bhcmp{2}_{\max(1, \alpha-2)} - \bhcmp{2}_{\max(1,\beta-2)}\right)\right) \chcmp{2}_{\max{(1,\beta-2)}}.
			\end{align*}
			\end{footnotesize}%
		\item {\em Option 3.} Leave a single $b_{3,k}=a_{s,3,k}$ free in an explicit stage, select another to enforce $\Mhcmp{1}$ and then zero-out all others:
			\begin{align}
				a_{s,3,k} = 
					\begin{cases}
						\eta, & k=\alpha, \\
						\ahcmp{1}_{32} - a_{s,3,\alpha} - \delta \left(1-\bhcmp{2}_{\max(1,\alpha-2)} - \bhcmp{2}_{\max(1,\beta-2)} \right), & k=\beta, \\
	 					\delta \bhcmp{2}_{\max(1,k-2)}, & \text{otherwise.}
	 				\end{cases}
			\end{align}
			Note that here we consider $\delta$ as a known value.  This choice of $a_{s,3,k}$ implicitly lead to
			\begin{align}
				a_{s,2,k} = 
					\begin{cases}
						\bhcmp{2}_{\max(1,\alpha-2)} - \eta - \Ahcmp{1}_{3,3} \delta_{\nu,k}, & k=\alpha,   \\
						\bhcmp{2}_{\max(1,\beta-2)} + a_{s,3,\alpha} - \Ahcmp{1}_{3,3} \delta_{\nu,k} \\
                         \qquad -\left[ \ahcmp{1}_{32} - \delta \left(1-\bhcmp{2}_{\max(1,\alpha-2)} - \bhcmp{2}_{\max(1,\beta-2)} \right) \right], & k=\beta, \\
	 					(1 - \delta) \bhcmp{2}_{\max(1,k-2)} - \Ahcmp{1}_{3,3} \delta_{\nu,k}, & \text{otherwise.}
	 				\end{cases}
			\end{align}
			Therefore with the exception of a maximum of three cases (i.e. $k \in \{\alpha, \beta, \nu\}$) the final stage has evaluations of the form $G(w) = (1 - \delta) F( {\color{magenta} Y_2}, w) + \delta F( {\color{ForestGreen} Y_3}, w)$ which matches with option 3 for the explicit stages. 
			
			If we define the following constants
			\begin{align*}
				X_1 &= \ahcmp{1}_{32} - \delta \left(1-\bhcmp{2}_{\max(1,\alpha-2)}-\bhcmp{2}_{\max(1,\beta-2)}\right), \\ 
				X_2 &= \bhcmp{2}_{\max(1,\alpha-2)}, \\
				X_3 &= \bhcmp{2}_{\max(1,\beta-2)} - X_1, \\
				X_4 &= \frac{1}{2} - \left( \bhcmp{2}_{\max(1,\alpha-2)} \chcmp{2}_{\max(1,\alpha-2)} + \bhcmp{2}_{\max(1,\beta-2)} \chcmp{2}_{\max(1,\beta-2)} \right),
			\end{align*}
			then, the order condition is
			\begin{align*}
				\frac{1}{3} &= \chcmp{1}_1 \Big( (1 - \delta)X_4 + \chcmp{2}_{\max(1,\alpha-2)} \left( -a_{s,3,\alpha} + X_2\right) \\
                & \qquad \qquad + \chcmp{2}_{\max(1,\beta-2)} \left( a_{s,3,\alpha} + X_3 \right) - \ahcmp{1}_{33} \chcmp{2}_{\max(1,\nu-2)} \Big) \\ 
                & \qquad  + \chcmp{1}_2 \left(  \delta X_4 + \chcmp{2}_{\max(1,\alpha-2)} a_{s,3,\alpha} +\chcmp{2}_{\beta-2} (-a_{s,3,\alpha} + X_1)\right) \\ 
                & \qquad +	\acmp{2}_{33} \chcmp{2}_{\max(1,\nu-2)},
			\end{align*}
			which has the solution
            \begin{subequations}
			\begin{align}
				a_{s,3,\alpha} &= \frac{\tfrac{1}{3} - \ahcmp{1}_{33} \chcmp{2}_{\max(1,\nu-2)} - \chcmp{1}_{1}N_1 - \chcmp{1}_{2}N_2}{(\chcmp{1}_2-\chcmp{1}_1)(\chcmp{2}_{\max(1,\alpha-2)}-\chcmp{2}_{\max(1,\beta-2)})}, \\
				N_1 &= (1-\delta)X_4+\chcmp{2}_{\max(1,\alpha-2)}X_2 + \chcmp{2}_{\max(1,\beta-2)}X_3 - \ahcmp{1}_{33} \chcmp{2}_{\max(1, \nu-2)}, \\
				N_2 &= \delta X_4 + \chcmp{2}_{\max(1,\beta-2)} X_1.
			\end{align}
            \end{subequations}
	\end{itemize}

This derivation allows for a large number of possible methods depending on the choice of constants. We summarize the free parameters below:
\begin{itemize}
	\item {\em Explicit stages:} the following only apply to options 2, and 3.
	\begin{center}
		\begin{tabular}{lll}
			param & restriction & description \\ \hline
			$u$ & $u \in \{4, \ldots, \scmp{2}+2\}$ & Coupling coefficient is in stage $Y_u$. \\
			$w$ & $w \in \{1\} \cup \{ 4, \ldots, u - 1\}$ & Coupling coefficient is $a_{u,3,w}$.
		\end{tabular}
	\end{center}
	\item {\em Output stages:} the following only apply to all options. Recall that since the method is stiffly accuracy $b_{jk} = a_{s,j,k}$.
	\begin{center}
		\begin{tabular}{lll}
			param & restriction & description \\ \hline
			$\alpha$ & $\alpha \in \{1\} \cup \{ 4, \ldots, s-1\}$ & Coupling coefficients include $a_{s,3,\alpha}$. \\
			$\beta$ & $\beta \in \{1\} \cup \{ 4, \ldots, s - 1\}$ & Coupling coefficients include $a_{s,3,\beta}$. \\
			$\nu$ & $\nu \in \{1\} \cup \{ 4, \ldots, s - 1\}$ & Implicit coefficient is $a_{s,s,\nu}$, which \\
			&& also forces a condition on $a_{3,s,\nu}$.
		\end{tabular}
	\end{center}
\end{itemize}

\subsection{Method selection:}

In the main paper text, the method listed in \cref{tab:mr-nprk-third-order-v1} is obtained by selecting parameters so that sub-stepping on is done on $w'=G(w)$ for $G(w) = F( {\color{magenta} Y_2}, w)$ with requirement to save initial and final stage. Specifically:
	\begin{itemize}
		\item {\bf Configuration 1:} substep $w'=G(w)$ for $G(w) = F( {\color{magenta} Y_2}, w)$.
		\begin{itemize}
			\item {\em  Implicit coupling:} 		$\omega=2$, $\nu = s-1 = \scmp{2}+2$ 
			\item {\em Explicit stages:} 		Option 1: $u = s-1 = \scmp{2} + 2$ and $w = s-2 = \scmp{2} + 1$
			\item {\em Output stages:} 		Option 1: $\alpha = 1$, $\beta = s-1 = \scmp{2}+2$ (requires that one save first and final during sub-stepping).	
		\end{itemize}
	\end{itemize}
Since this method can have unbounded coupling coefficients (See main paper text \cref{remark:mr3-v1-unbounded-coeff}), we propose an additional method here for anyone who wishes to use MR-NPRK methods with large $\scmp{2}$. The complete formula and implementation of this method are listed in \cref{subsec:mr3-alternative}.
	\begin{itemize}
		\setcounter{enumi}{1}
		\item {\bf  Configuration 2:} substep $w'=G(w)$ for $G(w) = (1 - \delta) F( {\color{magenta} Y_2}, w) + \delta F( {\color{ForestGreen} Y_3}, w)$ 
		\begin{itemize}
			\item {\em Implicit coupling:} 		$\omega=2$, $\nu = s-1 = \scmp{2}+2$ 
			\item {\em Explicit stages:} 		Option 3
			\item {\em Output stages:} 		Option 3  with $\alpha = 1$, $\beta = s-1 = \scmp{2}+2$	 (requires that one save first and final during sub-stepping).
		\end{itemize}
		{\bf Important:} stable coefficients as $\scmp{2} \to \infty$ however twice the cost in terms of explicit stages.	
	\end{itemize}

	\begin{remark}
		The coupling coefficients shown in \cref{tab:mr-nprk-third-order-v1} (main paper) and \cref{tab:mr-nprk-third-order-v2} (Supplemental) are simplification of the previous equations using the facts:
			\begin{itemize}
				\item For $\alpha = 1$ it follows that $\chcmp{2}_{\max(1,\alpha-2)} = c_1 = 0$
				\item For $\beta = s-1 = \scmp{2}+2$ it follows that $\chcmp{2}_{\max({1,\beta-2})} = \chcmp{2}_{\scmp{2}}$
			\end{itemize}	
	\end{remark}

\subsection{Alternative 3rd Order MP-NPRK with stable coefficients}
\label{subsec:mr3-alternative}

In \Cref{tab:mr-nprk-third-order-v2} we propose a third-order method ansatz that pads an $\scmp{2}$-stage explicit RK method ($\AhcmpB{2}$, $\bhcmpB{2}$) with two implicit stages at the beginning and one at end of the method. 
A detailed derivation of this method is contained in \Cref{sup:third-order-derivation}.

\subsubsection*{Practical Implementation for MR-NPRK method \eqref{eq:mr-nprk-third-order-coupling-coeff-v2}:}
	Given the initial condition $y_n$ and an explicit RK method $(\AhcmpB{2}, \bhcmpB{2})$ of at least third order:
		\begin{enumerate}
			\item Solve the implicit equations for ${\color{magenta} Y_2}$ and ${\color{ForestGreen} Y_3}$.
			\item Consider the ODE $w' = (1 - \delta) F({\color{magenta} Y_2}, w) + \delta F({\color{ForestGreen} Y_3}, w)$, $w_0=y_n$  and compute one timestep with stepsize $h$ using the explicit method $(\ahcmp{2}_{ij}, \bhcmp{2}_i)$. In addition to the output $w_1$, also store the final stage value $W_{\scmp{2}}$, which is equal to $Y_{s-1}$.
			\item To obtain {\color{Cerulean} $Y_s$} solve the implicit system ${\color{Cerulean} Y_{s}} = x + h\ahcmp{1}_{33} F({\color{Cerulean} Y_s},W_{\scmp{2}})$ with
			\begin{align} 
					x = w_1 & + h\Big[  (a_{s,3,1} - \delta \bhcmp{2}_1) H(y_n) \label{eq:b-vector-mr-nprk-third-order-coupling-coeff-v2} \\
                    & \qquad + (a_{s,3,s-1} - \delta \bhcmp{2}_{\scmp{2}}) H(W_{\scmp{2}})  - \ahcmp{1}_{33}F({\color{magenta} Y_2},W_{\scmp{2}} )	\Big], \nonumber
			\end{align}
			where $H(w)$ is defined in \eqref{eq:h-fun-mr-nprk-third-order-coupling-coeff-v1}.
		\end{enumerate}

	\begin{center}
		\begin{tabular}{|p{0.95\textwidth}|}
		
			\hline
			\vspace{-0.25em}
			
			{\bf Pseudocode for MR-NPRK method \eqref{eq:mr-nprk-third-order-coupling-coeff-v2}} \hfill $\omega \in \{1,2\}$ and remaining coefficients in \Cref{tab:mr-nprk-third-order-v2} \\[.25em] \hline         
			
			\vspace{-0.25em}
			Solve ${\color{magenta} Y_2} = y_n + h\ahcmp{1}_{11} F({\color{magenta} Y_2},y_n)$ \\[0.5em]
			
			Solve ${\color{ForestGreen} Y_3} = y_n + h \left[ (\ahcmp{1}_{21} - a_{3,2,2}) F({\color{magenta} Y_2},y_n) + a_{3,2,2} F({\color{magenta} Y_2},{\color{magenta} Y_2}) + \ahcmp{1}_{22} F({\color{ForestGreen} Y_3},Y_\omega)  \right]$ \\[0.5em]
			
			[$w_1$, $W_{\scmp{2}}$] = rkStep($w' = (1 - \delta) F({\color{magenta} Y_2}, w) + \delta F({\color{ForestGreen} Y_3}, w)$, $w_0 = y_n$, $\AhcmpB{2}$, $\bhcmpB{2}$, $h$) \\[.5em]
			
			Solve ${\color{Cerulean} Y_s} = x + h F({\color{Cerulean}Y_s}, W_{\scmp{2}})$ \hfill ($x$ defined in \eqref{eq:b-vector-mr-nprk-third-order-coupling-coeff-v2}) \\[.5em] \hline 
		        
		\end{tabular}
	\end{center}

\begin{table}[h!]
		
		\begin{tabular}{|p{0.95\textwidth}|}
		\hline \\[-0.5em]
		{ 
			Select either $\omega=1$ or $\omega = 2$. The third-order MR-NPRK method ansatz is:
		\begin{align}
			\begin{aligned}
				Y_1 &= y_n, \\
				{\color{magenta} Y_2} &= y_n + h\ahcmp{1}_{11} F({\color{magenta} Y_2},Y_1), \\
				{\color{ForestGreen} Y_3} &= y_n + h \left[ (\ahcmp{1}_{21} - a_{3,2,2}) F({\color{magenta} Y_2},Y_1) + a_{3,2,2} F({\color{magenta} Y_2},{\color{magenta} Y_2}) + \ahcmp{1}_{22} F({\color{ForestGreen} Y_3},Y_\omega),  \right] \\
				Y_i &= y_n + h \sum_{\substack{k=1 \\ k \ne 2,3}}^{i-1} \ahcmp{2}_{i-2,\max(k-2,1)} \left[ (1 - \delta) F({\color{magenta} Y_2},Y_k) + \delta F({\color{ForestGreen} Y_3},Y_k) \right], \\
                & \hspace{2 em} \text{for } i = 4,\ldots, \scmp{2}+2, \\
				{\color{Cerulean} Y_{s}} &= y_n +  
					h \left[ \bhcmp{2}_{1} - a_{s,3,1} \right]F({\color{magenta} Y_2}, Y_1) + h a_{s,3,1} F({\color{ForestGreen} Y_3}, Y_1) \\				
					& \hspace{2.4em}
					+ h\sum^{s-2}_{k=4} \bhcmp{2}_{k-2} \left[ (1 - \delta) F({\color{magenta} Y_2},Y_k) + \delta F({\color{ForestGreen} Y_3},Y_k) \right]  \\
					& \hspace{2.4em}
					+ h \left[ \bhcmp{2}_{\scmp{2}} - a_{s,3,s-1} - \ahcmp{1}_{33} \right]F({\color{magenta} Y_2}, Y_{s-1}) \\
                    & \hspace{2.4em} + h a_{s,3,s-1} F({\color{ForestGreen} Y_3}, Y_{s-1}) + h\ahcmp{1}_{33} F({\color{Cerulean} Y_s},Y_{s-1}), \\					
				y_{n+1} &= {\color{Cerulean} Y_{s}},
			\end{aligned}
			\label{eq:mr-nprk-third-order-coupling-coeff-v2}
		\end{align}
		where $s = \scmp{2}+3$. 	The coupling coefficients are
			\begin{align*}
				a_{3,2,2} &= 
				\begin{cases}
			 		\frac{2-6 \lambda }{18 \lambda ^3-60 \lambda ^2+15 \lambda } \approx 0.182536720446875179809517453138 & \omega = 1 \\
			 		\frac{2-6 \lambda }{18 \lambda ^3-60 \lambda ^2+15 \lambda }-\lambda \approx -0.253329801061583819606501998055 & \omega = 2
			 	\end{cases} \\
				\delta &= \frac{2 (3 \lambda -1)}{3 (\lambda -1)} \approx -0.36350683689006809456061097715969 	\\
				a_{s,3,1} &= \frac{\tfrac{1}{3} - \ahcmp{1}_{33} \chcmp{2}_{\scmp{2}} - \chcmp{1}_{1}N_1 - \chcmp{1}_{2}N_2}{(\chcmp{1}_1-\chcmp{1}_2)\chcmp{2}_{\scmp{2}}} \\
				a_{s,3,s-1} &= \ahcmp{1}_{32} - a_{s,3,1} - \delta \left(1-\bhcmp{2}_{1} - \bhcmp{2}_{\scmp{2}} \right)
			\end{align*}
		}
		where the constants $N_j$ are given by
		{ 
			\begin{align*}
				X_1 &= \ahcmp{1}_{32} - \delta \left(1-\bhcmp{2}_{1}-\bhcmp{2}_{\scmp{2}}\right) &
				X_3 &= \bhcmp{2}_{\scmp{2}} - X_1 &
				X_4 &= \frac{1}{2} - \bhcmp{2}_{\scmp{2}} \chcmp{2}_{\scmp{2}}, \\
				N_1 &= (1-\delta)X_4 + \chcmp{2}_{\scmp{2}}(X_3 - \ahcmp{1}_{33}) &
				N_2 &= \delta X_4 + \chcmp{2}_{\scmp{2}} X_1.
			\end{align*}
		}
		\\[-0.5em] \hline
		
		\end{tabular}
		
		\caption{Third-order MR-NPRK \eqref{eq:mr-nprk-third-order-coupling-coeff-v2} with bounded coefficients as $\scmp{2} \to \infty$.}
		\label{tab:mr-nprk-third-order-v2}
		\end{table}


\vfill
\begin{center}
	\includegraphics[width=1em]{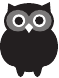}	
\end{center}

\end{document}